\documentclass[10pt]{article}
\pdfoutput=1

\usepackage[top=1in, bottom=1in, left=1.25in, right=1.25in]{geometry}
\usepackage[T1]{fontenc}

\usepackage{amsmath,amssymb,amsthm,graphicx, mathtools}
\newtheorem{theorem}{Theorem}
\newtheorem{lemma}{Lemma}

\usepackage{graphics,graphicx,color}         	% Allows inclusion of eps files.
\graphicspath{{/EPSF/}{../Figures/}{Figures/}}
\usepackage{url}		% Allows good typesetting of web URLs.
\usepackage{hyperref}
\usepackage[normalem]{ulem}

\usepackage{floatrow}
\newfloatcommand{capbtabbox}{table}[][0.5\textwidth]

\usepackage{xr} % cross references
\usepackage{appendix}
\usepackage{color}  % Added by Andrew
\usepackage{fancyhdr}

\input{setupmacros.sty} % Added by Sarah
\usepackage[boxed]{algorithm2e} % for algorithms

\newtheorem{prop}{Proposition}

\newcommand{\xcondh}{\hat{\boldsymbol{x}}_{\text{cond}}}
\newcommand{\condh}{\widehat{\boldsymbol{\Gamma}}_{\text{cond}}}
\newcommand{\xcond}{\boldsymbol{x}_{\text{cond}}}
\newcommand{\cond}{\boldsymbol{\Gamma}_{\text{cond}}}
\newcommand{\prior}{\boldsymbol{\Gamma}_{\text{pr}}}
\newcommand{\normtwo}[1]{\| #1\|_2}

\title{ Low Rank Independence Samplers in Bayesian Inverse Problems}

\author{
	D. Andrew Brown\thanks{
		Department of Mathematical Sciences,
		Clemson University, Clemson, SC 29634;
		\href{mailto:ab7@clemson.edu}{ab7@clemson.edu};
	}
	\and
	Arvind Saibaba\thanks{
		Department of Mathematics,
		North Carolina State University, Raleigh, NC 27695;
		\href{mailto:asaibab@ncsu.edu}{asaibab@ncsu.edu};
	}
	\and
	Sarah Vall\'elian\thanks{
		Statistical and Applied Mathematical Sciences Institute,
		Research Triangle Park, NC 27709;
		\href{mailto:svallelian@samsi.info}{svallelian@samsi.info};
	}
}

\begin{document}

\maketitle

\begin{abstract}
In Bayesian inverse problems, the posterior distribution is used to quantify uncertainty about the reconstructed solution. In practice,  Markov chain Monte Carlo algorithms often are used to draw samples from the posterior distribution. However, implementations of such algorithms can be computationally expensive. We present a computationally efficient scheme for sampling high-dimensional Gaussian distributions in ill-posed Bayesian linear inverse problems. Our approach uses Metropolis-Hastings independence sampling with a proposal distribution based on a low-rank approximation of the prior-preconditioned Hessian. We show the dependence of the acceptance rate on the number of eigenvalues retained and discuss conditions under which the acceptance rate is high. We demonstrate our proposed sampler by using it with Metropolis-Hastings-within-Gibbs sampling in numerical experiments in image deblurring, computerized tomography, and NMR relaxometry.
\end{abstract}

%\begin{keywords}  % AB - I was taught to list key words alphabetically, but I'm not sure what the convention is in other communities

%%%

%%Hierarchical Bayesian approach
%%Sampling schemes -
%%Block Gibbs, Metropolis-within-Gibbs (Hastings independence chain)
%%Acceptance Ratio and why it should be high
%%Introduce our sampler and why it will be high.
%%Metrics for assessing diagnostic.
%%Applications: 1D, CT and EEG.
%%Conclusions

% Discussion points:
%       1. Our approach is not dependent on the prior distributions for the precision parameters
%       2. High autocorrelation in the realizations of the precision parameters, which has been addressed using partially collapsed Gibbs sampling and MTC sampling. Our procedure would
%           still be helpful here, though, since PC-Gibbs still involves sampling from the full conditional distribution of x
%       3. Making a fully Bayesian approach computationally feasible makes it easier to include more prior information when solving the inverse problem (e.g., information from fMRI %   %               activation studies for source localization) while accounting for all sources of uncertainty
%       4. MAP estimator is interpreted as the posterior mode in the Bayesian paradigm. With the full posterior distribution, virtually any estimator we want is available, depending on %          criteria (0-1 loss (classification) <-> MAP,  absolute error loss <-> posterior median, squared error loss (most common) <->  posterior mean, etc.)

% intro section

\section{Introduction}\label{SEC:Intro}
Inverse problems aim to recover quantities that cannot be directly observed, but can only be measured indirectly and in the presence of measurement error. Such problems arise in many applications in science and engineering, including medical imaging~\cite{Hauk2004}, earth sciences~\cite{AndersenEtAl03}, and particle physics~\cite{NakhlehEtAl13}. The deterministic approach to inverse problems involves minimizing an objective function to obtain a point estimate of the unknown parameter. Inverse problems also admit a Bayesian interpretation, facilitating the use of prior information and allowing full quantification of uncertainty about the solutions in the form of a posterior probability distribution. An overview of Bayesian approaches to inverse problems is available in~\cite{KaiSom05, Mosegaard2002, Stuart2010}. A recent special issue of {\em Inverse Problems} also highlights the advances in the Bayesian approach and the broad impacts of its applicability~\cite{Calvetti2014}.

In the Bayesian statistical framework, the parameters of interest, $\B{x}$, and the observed data, $\B{b}$, are modeled as random variables. {\em A priori} uncertainty about the parameters is quantified in the prior distribution, $\pi(\B{x})$. Bayesian inference then proceeds by updating the information about these parameters given the observed data. The updated information is quantified in the posterior distribution obtained via Bayes' rule, $\pi(\B{x} | \B{b}) \propto f(\B{b} | \B{x})\pi(\B{x})$, where $f(\cdot | \B{x})$ is the assumed data generating model determined by the forward operator and $\B{x}$, called the likelihood function. Rather than providing a single solution to the inverse problem, the Bayesian approach provides a distribution of plausible solutions. Thus, sampling from the posterior distribution allows for simultaneous estimation of quantities of interest and quantifying the associated uncertainty.
A challenge of the fully Bayesian approach is that the posterior distribution will usually not have a closed form, in which case approximation techniques become necessary. In light of this, an indirect sampling-based approximation often is used to explore the posterior distribution. Since the seminal work of Gelfand and Smith \cite{GelfandSmith90}, Markov chain Monte Carlo (MCMC), particularly Gibbs sampling \cite{GemanGeman84}, has become the predominant technique for Bayesian computation. Several MCMC methods for sampling the posterior distributions obtained from inverse problems have been proposed in the literature \cite{bardsley2012mcmc,FoxNorton16,bardsley2016partially,agapiou2014analysis,bardsley2016metropolis}. However, these methods can be computationally expensive on large-scale problems due to the need to factorize a large covariance matrix at each iteration; though there are cases in which the choice of the prior and the forward operator lead to a reduction in computational cost \cite{bardsley2013efficient}. Approximating complex, non-Gaussian posteriors without the computational intensity of MCMC is still an ongoing area of research; e.g., variational Bayes \cite{JordanEtAl99} and integrated nested Laplace approximation \cite{RueEtAl09}. Each approach has features and caveats, a full exposition of which is beyond the scope of this paper. In this work, we assume that a researcher has already decided that they will use MCMC to access the posterior distribution. %In this case, Monte Carlo approximations can be used to approximate the quantities of interest. Even Monte Carlo, however, requires samples from the posterior distribution, which are generally unavailable so that indirect sampling approaches are necessary.
Our aim in this work is to address the computational burden posed by repeatedly sampling high dimensional Gaussian random variables as part of a larger MCMC routine, e.g., block Gibbs \cite{LiuEtAl94} or one-block \cite{RueHeld05}. We do so by leveraging the low-rank structure of forward models typically encountered in linear inverse problems. Specifically, we propose a Metropolis-Hastings independence sampler in which the proposal distribution, based on a low-rank approximation to the prior-preconditioned Hessian, is easy to construct and to sample. We also develop a proposal distribution using a randomized approach for computing the low-rank approximation when doing so directly is computationally expensive. We derive explicit formulas for the acceptance rates of our proposed approaches and analyze their statistical properties. We provide a detailed description of the computational costs. Numerical experiments support the theoretical properties of our proposed approaches and demonstrate the computational benefits over standard block Gibbs sampling.

%Our aim in this work is to address the computational burden incurred by a standard block Gibbs algorithm by leveraging the mathematical structure of Bayesian inverse problems, in particular the low-rank structure of the forward model. We use a Metropolis-Hastings-within-Gibbs (MHwG) algorithm for sampling from the posterior distribution. The computationally expensive part of a block Gibbs algorithm is drawing repeated samples from the high-dimensional Gaussian distribution of the estimand, conditional on the remaining parameters and the data. To approximate this step, we sample from a proposal distribution based on a low-rank approximation to the prior-preconditioned Hessian which is easy to construct and to sample from, while noting that the proposal distribution can be used also as a candidate distribution in a Rejection algorithm. When the computation of a low-rank approximation is expensive, we also develop a distribution using a randomized approach to compute the low-rank approximation. We derive explicit formulas for the acceptance rate from our proposed approach and analyze its statistical properties. Additionally, we provide a detailed description of the computational costs. Numerical experiments on challenging applications, including image deblurring and computerized tomography (CT), confirm the theoretical properties of our proposed sampling approach and demonstrate the computational benefits over standard approaches.

The rest of the paper is organized as follows. In Section~\ref{SEC:BayesIP}, we formulate a general linear inverse problem in the hierarchical Bayesian framework with particular attention paid to the computational bottleneck arising in standard MCMC samplers. In Section~\ref{SEC:Proposal}, we present our proposed approach of using low rank approximation as the basis of an independence sampler to accelerate drawing realizations from high-dimensional Gaussian distributions. In Section~\ref{sec:Simulation}, we demonstrate the performance of our approach on simulated examples in image deblurring and CT reconstruction via Metropolis-Hastings-within-Gibbs sampling. The paper concludes with a discussion in Section \ref{SEC:Discussion} and proofs of stated results in an Appendix. Further numerical studies, including convergence and alternative parameterizations for MCMC, are presented in Supplementary Material to this paper.

% Hierarchical Bayesian inverse problems approach

\section{The Bayesian Statistical Inverse Problem}\label{SEC:BayesIP}
%In this section, we formulate our inverse problem in the hierarchical Bayesian framework and define the posterior distribution we wish to sample from. We

Assume that the observed data are corrupted by additive noise so that the stochastic model for the forward problem is
	\begin{equation}\label{EQ:StochasticModel}
	\B{b} = \B{Ax} + \boldsymbol{\epsilon},
	\end{equation}
where $\B{A} \in \mathbb{R}^{m \times n}$ is the forward operator, or the parameter-to-observation map, $\boldsymbol{\epsilon}$ is the measurement error, and $\B{x}$ is the underlying quantity that we wish to reconstruct. We suppose that $\boldsymbol{\epsilon}$ is a Gaussian random variable with mean zero and covariance $\mu^{-1}\B{I}$, independent of the unknown $\B{x}$. %This illustrates how $\mu$ weighs the certainty or uncertainty, i.e. the variance, we have in our measured data $\B{b}$.
In some applications, $\mu$ may be known. Quite often, however, it is unknown and we assume that is the case here. Under this model, $\B{b} \mid \B{x},\mu \sim \mc{N}(\B{Ax}, \mu^{-1}\B{I})$ so that the likelihood is
	\begin{equation}\label{EQ:Likelihood}
		f(\B{b} \mid \B{x}, \mu) \ \propto \ \mu^{m/2} \exp \left( -\frac{\mu}{2} (\B{b} - \B{Ax})^\top (\B{b} - \B{Ax} )\right), ~\B{b} \in \mathbb{R}^m.
	\end{equation}
The prior distribution for $\B{x}$ encodes the structure we expect or wish to enforce on $\B{x}$ before taking observed data into account. An often reasonable prior for $\B{x}$ is Gaussian with mean zero and covariance $\sigma^{-1} \prior \equiv \sigma^{-1} (\B{L}^\top \B{L})^{-1}$; i.e.,
\begin{equation}\label{EQ:Priors-X}
	\pi(\B{x} \mid \sigma) \ \propto \ \sigma^{n/2} \exp \left( -\frac{\sigma}{2} \B{x}^\top \prior^{-1} \B{x}  \right), ~\B{x} \in \mathbb{R}^n,
\end{equation}
where the covariance matrix $\prior$ is assumed known up to the precision $\sigma$.
%For the prior to be a proper distribution this requires that $\prior$ be symmetric positive definite. By choosing the prior for the unknown $\pi(\B{x} \mid \sigma)$ to be Gaussian we ensure that the posterior distribution is also Gaussian.

Different covariance matrices may be chosen depending on what structure one wishes to enforce on the estimand $\B{x}$. The prior structure we use in our numerical experiments (Section \ref{sec:Simulation}) is motivated by Gaussian Markov random fields (GMRFs) \cite{RueHeld05}. Other popular choices involve Gaussian processes (GP) \cite{RasWilliams06}, which are parameterized in terms of covariance kernels.
%%%%%%%%% Arvind: removed this discussion to save some space
%\sout{Several choices of covariance kernel for the GP are available, with the most common being in the Mat\'ern family.  In this approach, the precision matrix is expressed as the discretized differential operator $(\kappa^2 - \Delta)^{\alpha/2}$, where $\alpha = \nu + d/2$, $\nu$ is the smoothness parameter of the Mat\'ern kernel, and $d$ is the dimension. An explicit connection between GPs and GMRFs is presented in~\cite{lindgren2011explicit}. An advantage of the GMRF representation is that the precision matrix is sparse and thus computationally efficient to handle when $\alpha$ is an even integer. GPs, on the other hand, more easily accommodate stationarity and generally allow for smoother fields. However, GP covariance kernels have infinite support so that the resulting covariance matrices are dense, creating challenges for large-scale applications. Fast Fourier Transform-based approaches on regular grids (e.g.,~\cite{ambikasaran2013fast}) and the $\mathcal{H}$-matrix approach on irregular grids (e.g.,~\cite{saibaba2012efficient,saibaba2015fastc}), though, can make the computational costs comparable to those associated with GMRFs.}

Gaussian process covariance kernels typically depend on parameters other than the multiplicative precision. For example, the covariance matrix of a GP often takes the form $\prior = \sigma\B{R}(\B{\theta})$, where $\B{\theta}$ may be the correlation length or other parameters determining the covariance function. When $\B{\theta}$ is unknown, one can assign it a prior distribution and estimate it along with the other parameters in the model by, e.g., updating it on each iteration of an MCMC algorithm. Such repeated updates are not feasible for extremely high dimensional problems since each factorization $\prior^{-1} = \B{L}^\top \B{L}$ is too expensive. However, it is possible to assign a prior to $\B{\theta}$ and subsequently obtain an empirical Bayes estimate $\widehat{\B{\theta}}$ (e.g., marginal posterior mode as done in \cite{QianWu08}). This estimator can be plugged in to the covariance function so that $\B{R}(\widehat{\B{\theta}})$ is fixed. Thus, we assume in this work that $\prior$ is fixed up to a multiplicative constant, without much loss of generality. This makes available an {\em a priori} factorization that we use to construct a low-rank approximation. See \cite{CarlinLouis09} for more details concerning empirical Bayes estimation.

%The conditional distribution
%\[ \pi(\B{x}\mid \B{b},\mu,\sigma) \propto \]

Conditional on $\mu$ and $\sigma$, the Bayesian inverse problem as formulated in equations~\eqref{EQ:Likelihood} and~\eqref{EQ:Priors-X} yields $\B{x} \mid \B{b}, \mu,\sigma \sim \mc{N}(\xcond, \cond)$, where $\cond = (\mu \B{A}^\top \B{A} + \sigma \prior^{-1})^{-1}$ and $\xcond = \mu \cond \B{A}^\top\B{b}$; i.e.,
\begin{equation}\label{EQ:XFC}
    \pi(\B{x} \mid \B{b}, \mu, \sigma) \propto \exp\left( -\frac{\mu}{2} \| \B{Ax} - \B{b} \|^2_{2} - \frac{\sigma}{2} \| \B{Lx} \|^2_{2}\right).
\end{equation}
The conditional posterior mode, $\hat{\B{x}} = \argmax_{\B{x}\in \mathbb{R}^n} \pi(\B{x} \mid \B{b}, \mu, \sigma)$, is the minimizer of the negative log-likelihood $(\mu/2)\| \B{Ax} - \B{b} \|^2_{2} + (\sigma/2)\| \B{Lx} \|^2_{2}$ and thus corresponds to Tikhonov regularization in the deterministic linear inverse problem. In a fully Bayesian analysis in which $\mu$ and $\sigma$ are unknown, we assign them a prior $\pi(\mu, \sigma)$ so that they can be estimated along with other parameters. In this case, the joint posterior density becomes
\begin{equation}\label{EQ:Joint-Conditional-Explicit}
	\pi(\B{x},\mu,\sigma \mid \B{b}) \propto \ \mu^{m/2} \sigma^{n/2} \exp \left( -\frac{\mu}{2} \| \B{Ax} - \B{b} \|^2_{2} - \frac{\sigma}{2} \| \B{Lx} \|^2_{2}\right)\pi(\mu, \sigma).
\end{equation}
In Section~\ref{sec:Simulation}, we consider two different priors on the precision parameters, conditionally conjugate Gamma distributions and a so-called weakly informative prior.

% Edit this paragraph
With priors on the precision components, the full posterior distribution is no longer Gaussian and generally not available in closed form. Non-Gaussian posteriors can sometimes be approximated by a Gaussian distribution, but such an approximation can be poor, especially with high-dimensional parameter spaces or multi-modal posterior distributions (\cite{GelmanEtAl14}, Chapter 4). Thus we appeal to MCMC for sampling from the posterior distribution. A version of the basic block Gibbs sampler for sampling from \eqref{EQ:Joint-Conditional-Explicit} is given in Algorithm \ref{ALG:Basic-BlockGibbs}. Most often, $\mu$ and $\sigma$ are updated individually (especially when using conditionally conjugate Gamma priors), but this is not necessary. Typically, $\B{x}$ is drawn separate from $(\mu, \sigma)$ to take advantage of its conditionally conjugate Gaussian distribution.

%An early modification of Gibbs sampling is the so-called Metropolis-Hastings-within-Gibbs to deal with challenging full conditional distributions \cite{Muller91}.

\LinesNumbered
\begin{algorithm}[tb]
\SetAlgoLined
\DontPrintSemicolon
\SetKwInput{Input}{Input}
\SetKwInput{Output}{Output}
	\Input{Full conditional distributions of $\B{x} \mid \B{b}, \mu, \sigma$ and $(\mu, \sigma) \mid \B{x}, \B{b}$, sample size $N$, burn-in period $N_b$.}
	\Output{Approximate sample from the posterior distribution \eqref{EQ:Joint-Conditional-Explicit}, $\{ \B{x}_{(t)}, \mu_{(t)}, \sigma_{(t)} \}_{t=N_b +1}^N$.}
	\BlankLine
	Initialize $\B{x}_{(0)}$, $\mu_{(0)}$, and $\sigma_{(0)}$. \;
  	\For{$t=1$ to $N$} {
    	Draw $\B{x}_{(t)} \>\sim\> \mc{N}(\mu_{(t-1)}\cond^{(t)}\B{A}^\top\B{b},\cond^{(t)})$, where \label{EQ:Gibbs-Conditional-x}
     \hfill $\cond^{(t)} = (\mu_{(t-1)} \B{A}^\top \B{A} + \sigma_{(t-1)} \prior^{-1})^{-1}.$ \hfill	\;
		Draw $(\mu_{(t)}, \sigma_{(t)}) \>\sim\> \mu, \sigma \mid \B{x}_{(t)}, \B{b}$. \;%\pi(\mu \mid \B{b}, \B{x}_{(t)}, \sigma_{(t-1)})$. \;
        }
\caption{An outline of the standard block Gibbs algorithm for sampling the posterior density~\eqref{EQ:Joint-Conditional-Explicit}.}
\label{ALG:Basic-BlockGibbs}
\end{algorithm}

For any iterative sampling algorithm in the Bayesian linear inverse problem, the computational cost per iteration is dominated by sampling $\B{x} \mid \B{b}, \mu, \sigma$ in~\eqref{EQ:XFC}. While sampling from this Gaussian distribution is a very straightforward procedure, the fact that it is high-dimensional makes it very computationally intensive. To circumvent the computational burden, we propose so-called Metropolis-Hastings-within-Gibbs sampling \cite{Muller91}. Specifically, we substitute direct sampling with a Metropolis-Hastings independence sampler using a computationally cheap low-rank proposal distribution. We present our proposed approach in Section \ref{SEC:Proposal}.

\section{Independence Sampling with Low-Rank Proposals}\label{SEC:Proposal}
Here we briefly review independence sampling and discuss a proposal distribution that uses a low-rank approximation to efficiently generate samples from~\eqref{EQ:XFC}.
\subsection{Independence Sampling}
 Let $\B\theta \in\mathbb{R}^n$ and denote the (possibly unnormalized) target density by $h(\B{\theta})$. The Metropolis-Hastings (MH) algorithm \cite{MetropolisEtAl53, Hastings70} proceeds iteratively by generating at iteration $t$ a draw, $\B{\theta}_{\ast}$, from an available proposal distribution possibly conditioned on the current state, $\B{\theta}_{(t-1)}$, and setting $\B{\theta}_{(t)} = \B{\theta}_{\ast}$ with probability $\alpha(\B{\theta}_{(t-1)}, \B{\theta}_{\ast}) = h(\B{\theta}_{\ast})q(\B{\theta}_{(t-1)} \mid \B{\theta}_{\ast})/(h(\B{\theta}_{(t-1)})q(\B{\theta}_{\ast} \mid \B{\theta}_{(t-1)})) \wedge 1$,
where $q(\cdot \mid \B{\theta}_{(t-1)})$ is the density of the proposal distribution. This algorithm produces a Markov chain $\{\B{\theta}_{(t)}\}$ with transition kernel
\[
    K(\B{\theta}, \B{\theta_{\ast}}) = \alpha(\B{\theta}, \B{\theta}_{\ast})q(\B{\theta}_{\ast} \mid \B{\theta}) + \delta_{\B{\theta}}(\B{\theta}_{\ast})\left(1 - \int \alpha(\B{\theta}, \B{\theta}^{\prime})q(\B{\theta}^{\prime} \mid \B{\theta})d\B{\theta}^{\prime}\right),
\]
where $\delta_{\B{\theta}}(\cdot)$ is the point mass at $\B{\theta}$. Properties of the MH algorithm, including convergence to the target distribution, may be found in \cite{RobertCasella04} and elsewhere.

An independence Metropolis-Hastings sampler (IMHS) proposes states from a density that is independent of the current state of the chain. The proposal has density $q(\B{\theta}_{\ast} \mid \B{\theta}_{(t-1)}) \equiv  g(\B{\theta}_{\ast})$, and the ratio appearing in $\alpha(\B{\theta}_{(t-1)}, \B{\theta}_\ast)$ can be written as
\begin{equation}\label{EQN:IndepRatio}
    \frac{h(\B{\theta}_{\ast})g(\B{\theta}_{(t-1)})}{h(\B{\theta}_{(t-1)})g(\B{\theta}_{\ast})} = \frac{w(\B{\theta}_{\ast})}{w(\B{\theta}_{(t-1)})},
\end{equation}
where $w(\B{\theta}) \propto h(\B{\theta})/g(\B{\theta})$. The IMHS is similar to the rejection algorithm. The rejection algorithm draws a candidate value $\B{\theta}_{\ast}$ from an available generating distribution with density $g$ such that for some $M \geq 1$, $h(\B\theta) \leq Mg(\B\theta)$, for all $\B\theta$. It then accepts the draw with probability $h(\B\theta_{\ast})/Mg(\B\theta_{\ast})$. Rejection sampling results in an exact draw from the target distribution.

For both the IMHS and the rejection sampler, it is desirable for $g$ to match the target density as closely as possible and, hence, to have an acceptance rate as high as possible. At least, $g$ should generally follow $h$, but with tails that are no lighter than $h$ \cite{GelmanEtAl95, RobertCasella04}. These guidelines are in contrast to those prescribed for the more common random walk MH, in which the best convergence is generally obtained with acceptance rates between 20\% and 50\% \cite{GelmanEtAl95, Rosenthal11}. In the sequel, we discuss our proposed generating distribution, both as an independence sampler as well as its use in a rejection algorithm.

\subsection{Approximating the Target Distribution}
Samples from the conditional distribution $\mc{N}(\xcond, \cond)$ can be generated as $\B{x} = \xcond + \B{G} \B{\varepsilon}$,
where $\B{\varepsilon} \sim \mathcal{N}(\B{0},\B{I})$ and $\B{G}$ satisfies $\cond = \B{G}\B{G}^\top$. Forming the mean $\xcond$ and computing the random vector $\B{G} \B{\varepsilon}$ involve expensive operations with the covariance matrix. By leveraging the low-rank nature of the forward operator $\B{A}$, we can construct a fast proposal distribution for an independence sampler.

%For convenience, we suppress the dependence of covariance matrices on the iteration number.
Consider the covariance matrix $\cond= (\mu \B{A}^\top \B{A} + \sigma \B{L}^\top \B{L} )^{-1}$. Factorizing this matrix so that
\begin{equation}\label{EQ:PriorPreconditionedHessian}
\begin{aligned}
	\cond &= \B{L}^{-1} ( \mu \B{L}^{-\top} \B{A}^\top \B{A} \B{L}^{-1} + \sigma \B{I} )^{-1} \B{L}^{-\top}%(\B{L}^\top \mu \B{L}^{-\top} \B{A}^\top \B{A} \B{L}^{-1} \B{L} + \sigma \B{L}^\top \B{L} )^{-1}
%	&= \left( \B{L}^\top (\mu \B{L}^{-\top} \B{A}^\top \B{A} \B{L}^{-1} + \sigma\B{I} ) \B{L} \right)^{-1} \\
%= \B{L}^{-1} ( \mu \B{L}^{-\top} \B{A}^\top \B{A} \B{L}^{-1} + \sigma \B{I} )^{-1} \B{L}^{-\top},
\end{aligned}
\end{equation}
yields the so-called \emph{prior-preconditioned Hessian transformation} $\B{H} := \B{L}^{-\top} \B{A}^\top \B{A} \B{L}^{-1}$ ~\cite{bui2013computational, flath2011fast, petra2014computational, saibaba2015fastc}. For highly ill-posed inverse problems such as those considered here, $\B{A}$ either has a rapidly decaying spectrum or is rank deficient. The product singular value inequalities~\cite[Theorem 3.3.16 (b)]{HoJ91} ensure that  $\B{AL}^{-1}$ has the same rank as $\B{A}$ and the same rate of decay of singular values. A detailed discussion on the low-rank approximation of the prior-preconditioned Hessian is provided in~\cite[Section 3]{flath2011fast}.

We approximate $\B{H}$ using a truncated eigenvalue decomposition,
\begin{equation}\label{EQ:EigenvalueDecomp}
	\B{L}^{-\top} \B{A}^\top \B{A} \B{L}^{-1} \approx \B{V}_k \boldsymbol{\Lambda}_k \B{V}_k^\top,
\end{equation}
where $\B{V}_k \in \bbR^{n \times k}$ has orthonormal columns and $\boldsymbol{\Lambda}_k \in \bbR^{k \times k}$ is the diagonal matrix containing the $k \leq n$ largest eigenvalues of $\B{H}$. If $\text{rank}(\B{A}) = k$, then exact equality holds. The truncation parameter $k$ controls the tradeoff between accuracy on the one hand and computational and memory costs on the other.

%\paragraph{Approximating $\xcond$ and $\cond$}
We approximate the conditional covariance matrix $\cond$ by substituting~\eqref{EQ:EigenvalueDecomp} into~\eqref{EQ:PriorPreconditionedHessian},
\begin{equation}\label{EQ:ApproximateHatGamma-1}
	\condh \> \equiv \>\B{L}^{-1} \frac{1}{\sigma} ( \B{I} + \frac{\mu}{\sigma} \B{V}_k \boldsymbol{\Lambda}_k \B{V}_k^\top )^{-1} \B{L}^{-\top}.
\end{equation}
Using the Woodbury identity and the fact that $\B{V}_k$ has orthonormal columns, the right-hand side of \eqref{EQ:ApproximateHatGamma-1} becomes
\begin{equation*}%\label{EQ:ApproximateHatGamma-2}
	\condh \>  = \> \frac{1}{\sigma} \B{L}^{-1} ( \B{I} - \B{V}_k \B{D}_k \B{V}_k^\top ) \B{L}^{-\top},\qquad \B{D}_k = \text{diag}\left(\mu \lambda_j(\mu \lambda_j + \sigma)^{-1} : j= 1, \ldots, k\right) \in \bbR^{k \times k},
\end{equation*}
where $\lambda_j, ~j= 1, \ldots, k$, are the diagonals of $\boldsymbol{\Lambda}_k$. To approximate the mean $\xcond$, replace $\cond$ by $\condh$ so that $\xcondh = \mu\condh\B{A}^\top\B{b}$.
%\end{equation*}
With these approximations, the proposal distribution for our proposed independence sampler is $\mc{N}(\xcondh,\condh)$. Optimality of this low-rank approximation was studied in~\cite{spantini2015optimal}.

%\paragraph{Sampling from the proposal distribution}
A factorization of the form $\condh = \B{GG}^\top$ can be used to sample from $\mc{N}(\xcondh, \condh)$. It can be verified that $\B{G} := \sigma^{-1/2} \B{L}^{-1} (\B{I} - \B{V}_k\widehat{\B{D}}_k\B{V}_k^\top)$, with $\widehat{\B{D}}_k = \B{I} \pm \left( \B{I}-\B{D}_k\right)^{1/2}$, satisfies $\condh = \B{GG}^\top$. Since $\widehat{\B{D}}_k$ is diagonal and $k \ll n$, we obtain a computationally cheap way of generating draws from the high-dimensional proposal distribution $\mc{N}(\xcondh,\condh)$. Then we can use a Metropolis-Hastings step to correct for the approximation. This results in our proposed low rank independence sampler (LRIS).

\subsection{Analysis of Acceptance Ratio}
Here, we derive an explicit formula for evaluating the acceptance ratio for our proposed algorithm and provide insight into the conditions under which the proposal distribution closely approximates the target distribution. For simplicity of notation, we suppress the conditioning on $\B{b}, \mu$, and $\sigma$.

The target density is
\begin{equation}\label{EQ:Target-Density}
h(\B{x}) \> := \> \frac{1}{\sqrt{(2\pi)^n\det(\cond)}}\exp\left( -\frac{1}{2}(\B{x}- \xcond)^\top\cond^{-1}(\B{x}- \xcond)\right),
\end{equation}
and the proposal density, $g(\B{x})$, replaces $\xcond$ by $\xcondh$ and $\cond$ by $\condh$ in $h(\B{x})$.
%\begin{equation}\label{EQ:Proposal-Density}
%q(\B{x}) \> := \> \frac{1}{\sqrt{(2\pi)^n\det(\condh)}}\exp\left( -\frac{1}{2}(\B{x}- \xcondh)\condh^{-1}(\B{x}- \xcondh)\right),
%\end{equation}
%It follows that the MH acceptance ratio is  $\eta(\B{z},\B{x}) \> := \> \frac{h(\B{z})q(\B{x})}{h(\B{x})q(\B{z})}$.
%\textcolor{red}{which reduces to a simple and useful form.}
The following result gives a practical way to compute the acceptance ratio. It can be verified with a little algebra, so the proof is omitted.
\begin{prop}\label{p_mh_ratio1}
Let $\B{x}$ be the current state of the LRIS chain and let $\B{z}$ be the proposed state. Then the acceptance ratio can be computed as $\eta(\B{z} , \B{x}) = w(\B{z})/w(\B{x}), $
where $w(\B{x}) = \exp\left(-\B{x}^\top (\cond^{-1} -\condh^{-1}) \B{x} / 2 \right)$.
\end{prop}
An efficient implementation and the cost of computing this ratio is discussed in Subsection~\ref{ss_cost}. The quality of the low-rank approximation to the target distribution can be seen through the acceptance ratio.
\begin{prop}\label{p_mh_ratio2}
Let $\B{x}$ be the current state of the LRIS chain and let $\B{z}$ be the proposed state. Then the LRIS acceptance ratio can be expressed as
\begin{equation}\label{e_mh_ratio2}
 \eta(\B{z},\B{x}) = \exp\left( -\frac{\mu}{2}\sum_{j=k+1}^n \lambda_j \left[\left(\B{v}_j^\top\B{Lz}\right)^2 - \left(\B{v}_j^\top\B{Lx}\right)^2 \right] \right)\,.
\end{equation}
\end{prop}
\begin{proof} See Appendix~\ref{a_proofs}.
\end{proof}
This proposition asserts that the acceptance ratio is high when either $\mu$ is small or the discarded eigenvalues $\{ \lambda_j\}_{j=k+1}^n$ are small. The dependence of the acceptance ratio on the eigenvectors can be seen explicitly by writing $(\B{v}_j^\top\B{Lz})^2 - (\B{v}_j^\top\B{Lx})^2 = [ \B{v}_j^\top\B{L}(\B{z}+\B{x})][ \B{v}_j^\top\B{L}(\B{z}-\B{x})].$ Thus, if $\B{z} \pm \B{x} \perp \B{L}^\top \B{v}_j, ~j= k+1, \ldots, n$, then the acceptance ratio is $1$. %Furthermore, assuming that $\B{z} = \xcondh + \B{G\epsilon}_{\B{z}}$ and $\B{x} = \xcondh + \B{G\epsilon}_{\B{x}}$, this expression simplifies to $ \B{\epsilon}_{\B{z}} - \B{G\epsilon}_{\B{x}} \perp \B{v}_j$.

While Proposition \ref{p_mh_ratio2} provides insight into realizations of the acceptance ratio, the actual acceptance ratio is a random variable. The expected behavior and variability of this quantity can be understood through Theorem~\ref{p_expect}. To this end, define the constants
\begin{equation}\label{e_constant}
 N_\ell := \exp\left( \frac{\mu^2}{2\sigma}\sum_{j=k+1}^n \frac{\ell\mu\lambda_j}{\ell\mu\lambda_j + \sigma}(\B{b}^\top\B{A}\B{L}^{-1}\B{v}_j)^2\right)\prod_{j=k+1}^n \left( 1 + \frac{\ell\mu}{\sigma}\lambda_j\right)^{1/2}, ~~\ell= 1, 2, \ldots.
\end{equation}
%for integers $\ell \geq 1$.

\begin{theorem}\label{p_expect}
Let $\B{x}$ be the current state of the LRIS chain and let $\B{z}$ be the proposed state. Then the expected value and the variance of the acceptance ratio are given by
\begin{equation} \label{e_expect}
    \begin{aligned}
        e_{\eta} := \bbE_{\B{z} |\B{x}} [\eta(\B{z}, \B{x})] &= \frac{1}{ N_1w(\B{x})}\\
        v_{\eta}^2 := \mathbb{V}_{\B{z} |\B{x}} [\eta(\B{z}, \B{x})] &= \frac{1}{ w^2(\B{x})} \left(\frac{1}{N_2} - \frac{1}{N_1^2}\right).
    \end{aligned}
\end{equation}
where $\bbE_{\B{z} \mid \B{x}}(\cdot)$ denotes expectation conditional on $\B{x}$, $\mathbb{V}_{\B{z} \mid \B{x}}(\cdot)$ denotes the variance conditional on $\B{x}$, and $w(\B{x})$ is as defined in Proposition \ref{p_mh_ratio1}.
%where $M$ is defined in Proposition~\ref{p_mh_bound}.
\end{theorem}
\begin{proof}See Appendix~\ref{a_proofs}.
\end{proof}
Using this result, a straightforward application of Chebyshev's inequality~\cite{Resnick99} shows that for any $\epsilon > 0$, $\text{Pr}_{\B{z} \mid \B{x}} \left( \left|  \eta(\B{z},\B{x}) - (N_1 w(\B{x}))^{-1} \right| \geq \epsilon  \right) \leq \left(N_2^{-1} - N_1^{-2}\right)/[\epsilon^2 w^2(\B{x})],$
where $\text{Pr}_{\B{z} \mid \B{x}}(\cdot)$ denotes probability conditional on the current state $\B{x}$. Thus, we can construct conditional prediction intervals about the realized acceptance rate. For instance, at any given state $\B{x}$, $\text{Pr}_{\B{z}\mid \B{x}}(\eta(\B{z}, \B{x}) \in [e_{\eta} \pm 4.47v_{\eta}]) \geq 0.95$. In Appendix~\ref{a_proofs}, we derive expressions for all moments of the acceptance ratio.

It is clear from Theorem~\ref{p_expect} that if the eigenvalues $\{\lambda_j\}_{j=k+1}^n$ are zero, then the acceptance probability is $1$. Likewise, if the eigenvalues are nonzero but small in magnitude, then the acceptance rate is close to $1$. Further, consider the SVD of $\B{A}\B{L}^{-1} = \B{U\Sigma V}^\top$. Then $\B{b}^\top\B{A}\B{L}^{-1}\B{v}_j$ in \eqref{e_constant} is equal to $\sigma_j \B{b}^\top \B{u}_j$, where $\B{u}_j$ is the $j$th singular vector of $\B{A}\B{L}^{-1}$. Thus, the acceptance rate may be close to $1$ even if the components of the measurement $\B{b}$ along the left singular vectors of $\B{AL}^{-1}$ are small. This is closely related to filter factors that are used to analyze deterministic inverse problems~\cite{hansen2010discrete}.

These results establish the moments of the acceptance rate for fixed precision parameters $\mu$ and $\sigma$ and fixed rank $k$ of the proposal distribution. In practice, when running MCMC, $\mu$ and $\sigma$ will change on each iteration, meaning that the actual acceptance rate will vary from one iteration to the next. Thus, it may not be clear {\em a priori} which truncation level to use to achieve an acceptable acceptance rate while minimizing the computational cost. Of course, if the low-rank matrix is obtained from a rank-deficient forward model by discarding only the zero eigenvalues, then the acceptance rate is one for all $\mu$ and $\sigma$. Otherwise, a practitioner can employ an {\em adaptive} LRIS in which the acceptance rate is tracked during an initial burn-in period, adding rank to the distribution every, e.g., 100 iterations if the acceptance rate is too low. This allows finding the minimum number of eigenvalues needed to achieve high acceptance over the high probability region of $\mu$ and $\sigma$. Provided the adaptation stops after a finite number of iterations, convergence to the stationary distribution is still guaranteed \cite{CarlinLouis09}. An outline of the adaptive LRIS approach, along with practical guidelines to determine the target rank $k$, is given in the Supplementary Material.

\paragraph{Convergence and the Rejection Algorithm}

Our proposed candidate generating distribution $g(\B{x})$ bounds the target distribution up to a fixed constant as a function of the remaining eigenvalues in the low-rank approximation, as asserted by the next Proposition.
\begin{prop}\label{p_mh_bound}
The target density $h(\B{x})$ \eqref{EQ:Target-Density} and the proposal density $g(\B{x})$ can be bounded as $h(\B{x}) \leq N_1 g(\B{x})$ for all $\B{x}$, where $N_1 \geq 1$ is given in~\eqref{e_constant}.
\end{prop}

\begin{proof} See Appendix~\ref{a_proofs}.
\end{proof}

Proposition \ref{p_mh_bound} establishes that the subchain produced by our proposed sampler has stationary distribution $\pi(\cdot \mid { \B{b}}, \mu, \sigma)$ and is uniformly ergodic by \cite[Theorem 7.8]{RobertCasella04}; i.e., for $p \in \mathbb{N},$
\begin{equation}\label{eqn:ergodic}
    \|K^p(\B{x}, \cdot) - \pi(\cdot \mid { \B{b}}, \mu, \sigma)\|_{TV} \leq 2\left(1 - N_1^{-1}\right)^p ~~\forall \B{x} \in \text{supp} ~\pi,
\end{equation}
where $K^p(\B{x}, \cdot)$ is the $p$-step LRIS transition kernel starting from $\B{x}$ and $\|\cdot\|_{TV}$ denotes the total variation norm. Thus, if one runs several sub-iterations of the LRIS, the realizations will converge to a draw from the true full conditional distribution at a rate independent of the initial state. Convergence is faster as the remaining eigenvalues from the low-rank approximation become small, and is immediate when the remaining eigenvalues are zero.

Equation \eqref{eqn:ergodic} explicitly quantifies convergence of the subchain to the full conditional distribution as a function of the quality of the approximation to the target, quantified in $N_1$. However, when the LRIS is used inside a larger MCMC algorithm (e.g., Metropolis-Hastings-within-Gibbs), convergence of the entire Markov chain to its stationary distribution is affected not only by the LRIS proposal distribution, but also by modeling choices on the remaining parameters and the manner in which they are updated. There exist results for establishing geometric ergodicity of componentwise Metropolis-Hastings independence samplers and so-called two-stage Metropolis-Hastings-within-Gibbs algorithms \cite{JohnsonEtAl13} for which Proposition \ref{p_mh_bound} could be useful. To the best of our knowledge, though, more general effects of the proposal distribution on the convergence of a MHwG algorithm are unknown. While exploring this issue is beyond the scope of this work, we carry out in the Supplementary Material an empirical study in which we assess convergence of an MHwG chain as a function of the rank of the proposal. We observe that as the number of eigenvalues retained increases, the convergence of the LRIS algorithm becomes more rapid.

Proposition \ref{p_mh_bound} suggests also that the approximating distribution can be used in a rejection algorithm instead of LRIS. The proof of the Proposition shows that $\det (\condh) \geq \det (\cond)$, but each determinant is a generalized variance~\cite{JohnWich07}. When there are non-zero eigenvalues left out of the low-rank approximation, the proposal density will have heavier tails than the target density, a desirable property for a candidate distribution in a rejection algorithm~\cite{GelmanEtAl14}. Otherwise, the approximation is exact. We remark, however, that for a given candidate density $g$, LRIS is more efficient than a rejection algorithm in terms of variances of the concomitant estimators~\cite{Liu96}. Further, the rejection sampler requires knowledge of $N_1$, which depends on eigenvalues that may be unavailable.

\subsection{Generating Low-Rank Approximations}\label{ss_randsvd}
A major cost of our proposed sampler is in the precomputation associated with constructing the low-rank approximation. The standard approach for computing this low-rank approximation is to use a Krylov subspace solver (e.g., Lanczos method~\cite{saad1992numerical}) for computing a partial eigenvalue decomposition. Alternatively, we can compute the rank$-k$ singular value decomposition $\B{A} \B{L}^{-1} \approx \B{U}_k \B{\Sigma}_k \B{V}_k^\top$. Then the approximate low-rank decomposition can be computed as $ \B{H}\approx \B{V}_k\B{\Sigma}_k^2 \B{V}_k^\top$. Here we discuss a computationally efficient alternative. %Another approach is to use the Golub-Kahan bidiagonalization algorithm~\cite{simon2000low}.
%\subsubsection{Randomized SVD approach}\label{sss_randsvd}

Randomized SVD, reviewed in~\cite{halko2011finding}, is a computationally efficient approach for computing a low-rank approximation to the prior-preconditioned Hessian $\B{H}$. The basic idea of the randomized SVD approach is to draw a random matrix $\B\Omega \in \mathbb{R}^{n\times (k+p)}$, where the entries of $\B\Omega$ are i.i.d.\ standard Gaussian random variables. Here, $k$ is the target rank and $p$ is an oversampling parameter. An approximation to the column space of $\B{H}$ is computed by the matrix product $\B{Y} = \B{H\Omega}$. A thin-QR factorization $\B{Y}=\B{QR}$ is computed, and the resulting low-rank approximation to $\B{H}$ is given by
\begin{equation}\label{eqn:randlowrank}
\B{H} \approx \widehat{\B{H}} := \B{QQ}^\top\B{HQQ}^\top.
\end{equation}
This can be postprocessed to obtain an approximate low-rank decomposition of the form~\eqref{EQ:EigenvalueDecomp}. This is summarized in Algorithm~\ref{alg:randsvd}.

\LinesNumbered
\begin{algorithm}[tb]
\SetAlgoLined
\DontPrintSemicolon
\SetKwInput{Input}{Input}
\SetKwInput{Output}{Output}
\Input{Matrix $\B{H}\in \mathbb{R}^{n\times n}$ and random matrix $\B{\Omega} \in  \mathbb{R}^{n\times (k+p)}$.}
\Output{Approximate eigenvectors $\B{V}$ and approximate eigenvalues $\B{\Lambda}$.}
Compute $\B{Y} = \B{H\Omega}$ and thin-QR factorization $\B{Y}=\B{QR}$. \;
Compute $\B{T} = \B{Q}^\top\B{HQ}$ and its eigendecomposition $\B{T} = \B{U}\B{\Lambda}\B{U}^\top$. \;
Compute $\B{V} = \B{QU}$. \;
%Set $\widehat{\B{H}} = \B{V}_k\B{\Lambda}_k\B{V}_k^\top$.
\caption{Randomized SVD algorithm for computing low-rank decomposition.}\label{alg:randsvd}
\end{algorithm}

Similar to Theorem \ref{p_expect}, we can bound the expected value of the acceptance ratio under the randomized SVD approach.
\begin{theorem}\label{p_randsvd}
    Suppose we compute the low-rank approximation $\widehat{\B{H}}$ using Algorithm~\ref{alg:randsvd} with guess $\B\Omega \in  \mathbb{R}^{n\times (k+p)}$. Let $p \geq 2$ be the oversampling parameter. Then
    \[ \mathbb{E}_{\B\Omega \mid \B{z},\B{x}}\, [\eta(\B{z},\B{x})] \geq \exp\left(-\mu\|\B{Lz}\|^2 \left[\alpha \lambda_{k+1} + \beta\left(\sum_{j=k+1}^{n}\lambda_j^2\right)^{1/2}\right]\right),\]
    where $\alpha = 1 + \sqrt{\frac{k}{p-1}}, \beta = \frac{e\sqrt{k+p}}{p}$ and $\mathbb{E}_{\B{\Omega}\mid \B{z},\B{x}}$ denotes expectation w.r.t.\ $\B\Omega$ given the current state $\B{x}$ and the proposed step $\B{z}$.
\end{theorem}

\begin{proof} See Appendix~\ref{a_proofs}.
\end{proof}
The interpretation of this result is similar to Theorem~\ref{p_expect}. That is, if the eigenvalues of the prior-preconditioned Hessian $\B{H}$ are rapidly decaying or zero beyond the index $k$, then the expected acceptance rate, averaged over all random matrices $\B{\Omega}$, is high. %We provide a proof of this result.

% Do we need a citation of recommendation of the sampling parameter?
In practice, an oversampling parameter of $p \lesssim 20$ is recommended~\cite{halko2011finding}. As proposed, Algorithm~\ref{alg:randsvd} requires $2(k+p)$ matrix-vector products (matvecs) with $\B{H}$. The second round of matvecs required in Step 2 can be avoided by using the approximation~\cite[Section 2.3]{saibaba2016randomized}
\[
    \B{T} \approx (\B\Omega^\top \B{Q})^{-1}(\B\Omega^\top \B{Y})(\B{Q}^\top\B\Omega)^{-1}.
\]
This is an example of the so-called single pass algorithm. Other single pass algorithms are discussed in~\cite{tropp2016randomized}. In practice, the target rank $k$ may not be known, in which case a modified approach may be used to adaptively estimate the subspace~\cite[Algorithm 4.2]{halko2011finding}. %Another approach is to first compute a low-rank approximation of $\B{AL}^{-1} \approx \B{U}_k\B{\Sigma}_k\B{V}_k^\top$ using~\cite[Algorithm 5.1]{halko2011finding} and then computing the low-rank approximation $\B{H} \approx \B{V}_k\B{\Sigma}_k^2 \B{V}_k^\top$. The error analysis is similar and the details are omitted.

\subsection{Computational Costs}\label{ss_cost}

Denote the computational cost of a matvec with $\B{A}$ by $T_{\B{A}}$, and the cost of a matvec with $\B{L}$ and $\B{L}^{-1}$ as $T_{\B{L}}$ and $T_{\B{L}^{-1}}$, respectively. For simplicity, we assume the cost of the transpose operations of the respective matrices is the same as that of the original matrix.

It is difficult to accurately estimate the cost of the Krylov subspace method \textit{a priori}, but the cost is roughly 2 sets of matvecs with $\B{A}$ and $\B{L}^{-1}$ and an additional $\mc{O}(nk^2)$ operations.  The quantities $\B{A}^\top \B{b}$ and $\B{L}^{-\top}\B{A}^\top \B{b}$ can also be precomputed at a cost of $T_{\B{A}}$ and $T_{\B{A}}+T_{\B{L}^{-1}}$ flops, respectively. Generally speaking, this is the same asymptotic cost for randomized SVD. In practice, however, randomized SVD can be much cheaper since it only seeks an approximate factorization~\cite{halko2011finding}.

\begin{table}[tb]\centering
\begin{tabular}{c|c|c}
Operation  & Formula & Cost \\ \hline
Precomputation & Equation~\eqref{EQ:EigenvalueDecomp} & $2k(T_{\B{A}} + T_{\B{L}^{-1}}) + \mc{O}(nk^2)$ \\
Computing mean & $\xcondh = \mu\condh\B{A}^\top\B{b}$ & $T_{\B{L}^{-1}} + 4nk$ \\
Generating sample & $\B{x} \sim \mc{N}(\xcondh,\condh)$ & $T_{\B{L}^{-1}} + 4nk$ \\
Acceptance ratio & Proposition~\ref{p_mh_ratio1} & $T_{\B{A}} + T_{\B{L}} + 2n(k+2)$
\end{tabular}
\caption{Summary of computational costs of various steps in the LRIS.}
\label{tab:costs}
\end{table}

The cost of computing the mean $\xcondh$ involves the application of $\B{L}^{-1}$ and $(\B{I}-\B{V}_k\B{D}_k\B{V}_k^\top)$. This costs $T_{\B{L}^{-1}} + 4nk$ flops. Similarly, the cost of $\B{G\varepsilon}$ is also $T_{\B{L}^{-1}} + 4nk$ flops. The important point here is that generating a sample from the proposal distribution does not require a matvec with $\B{A}$. This is useful for applications in which $T_{\B{A}}$ can be extremely high. The computational cost of computing the acceptance ratio can be examined in light of Proposition~\ref{p_mh_ratio1}. On each iteration, the weight $w(\B{x})$ will already be available from the previous iteration, so we only need to compute $w(\B{z})$. We can simplify this expression as $\log w(\B{z}) = - \mu \B{z}^\top ( \B{A}^\top \B{A} - \B{L}^\top \B{V}_k\B{\Lambda}_k\B{V}_k^\top\B{L})\B{z}$, which requires one matvec with $\B{A}$ and $\B{L}$ each, two inner products and $4n$ flops, and an additional $2nk$ flops. Aside from the precomputational cost of the low-rank factorization, only the evaluation of the acceptance ratio requires accessing the forward operator $\B{A}$. The resulting costs are summarized in Table~\ref{tab:costs}.

\section{Illustrations}\label{sec:Simulation}  % Add computation times for cross-validation and MCMC? (The cross-validation took a LONG time, but it is dependent upon the resolution of the grid of values considered)

Here we demonstrate our proposed approach on two simulated examples. The first example is a standard two dimensional deblurring problem in which we compare the performance of our proposed low-rank independence sampler to conventional block Gibbs sampling to demonstrate the competitive solutions and the ability to access the posterior distribution in an efficient manner. The second example is a more challenging application motivated by medical imaging with a rank deficient forward model. We apply our proposed approach there to demonstrate feasibility and to consider a different prior on the precisions than the conventional independent conjugate Gammas.

%%Consider moving this paragraph to the appendix.
To ensure meaningful inferences based on the MCMC output, it is important to assess whether the Markov chain is sufficiently close to its stationary distribution. It is well known that an MCMC procedure will generally not result in an immediate draw from the target distribution, unless the initial distribution is the stationary distribution. Usually it is not possible to prove that a chain has converged to its limiting distribution, except in special cases (e.g., perfect sampling \cite{CraiuMeng11}). However, diagnostic tools %such as trace plots and potential scale reduction factors
can be used to assess whether or not a chain is sufficiently close so that one can safely treat its output as draws from the target distribution.
%Trace plots are graphical displays of particular scalar quantities versus the iteration number.  When multiple chains are run, each starting from different values well dispersed with respect to the target distribution, a sign of convergence is when the realizations of the chains tend to the same part of the parameter space and overlap when superimposed. Potential scale reduction factors (PSRFs), denoted $\widehat{R}$, are quantitative measures which use the {\em within} chain variability and {\em between} chain variability (of three or more chains) to estimate the factor by which the variance of the current distribution of a parameter of interest might be reduced if the simulations were allowed to run indefinitely \cite{GelmanRubin92}. This quantity tends to $1$ as the chain length tends to infinity, and $\widehat{R} \gg 1$ indicates the need to run the chain longer before approximate convergence is attained. A generalization of PSRFs to the multivariate case takes $\widehat{R}_{MV} := (N-1)/N + \lambda_1(m_c+1)/m_c$, where $N$ is the chain length, $\lambda_1$ is the largest eigenvalue of $W_{MV}^{-1}B_{MV}/N$, and $W_{MV}$ and $B_{MV}$ are the within- and between-chain covariance matrices, respectively \cite{BrooksGelman98}. As in the scalar case, $\widehat{R}_{MV} \rightarrow 1$ as $N \rightarrow \infty$. We use these diagnostic tools to assess convergence of the Markov chains in the following applications.
To diagnose convergence, we use (scalar and multivariate) potential scale reduction factors (PSRF/MPSRF) \cite{GelmanRubin92, BrooksGelman98}, trace plots, and autocorrelation plots. The reader is referred to \cite[Ch. 12]{RobertCasella04}, \cite[Ch. 3]{CarlinLouis09},  or \cite[Ch. 11]{GelmanEtAl14}  for further discussions of convergence diagnostics for MCMC.

%Common diagnostic tools used in practice are .  These quantities can either be realizations of particular parameter states in the chain, or of partial sums used to calculate ergodic averages. on the trace plot. One can also use the trace plot of a single chain, looking for when the chain seems to ``settle down" without moving from one part of the space to another. The latter approach is more susceptible to {\em pseudo-convergence} in which the chain gets trapped in a local mode of the target distribution, without fully exploring all of the high probability areas. P The rule of thumb for sufficient convergence \cite{GelmanEtAl14} is $\widehat{R} \leq 1.1$.   See

%All simulations were performed using MATLAB on a desktop running OS X Yosemite (8GB RAM, Intel Core i5 2.66GHz processor).

%\input{5_1D}
%\input{5_eeg}

\subsection{2D Image Deblurring}\label{subsec:2d}
\begin{figure}[tb]
        \centering
    \includegraphics[scale = 0.5]{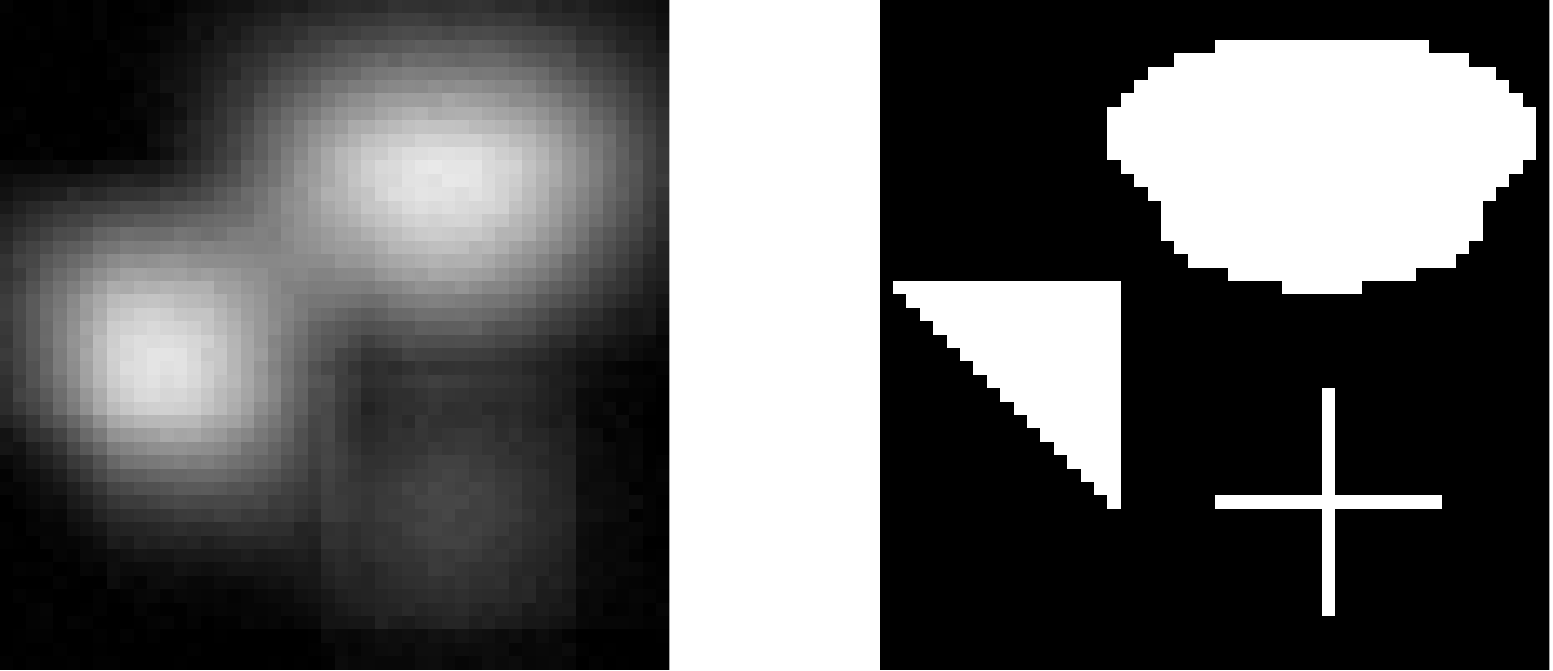}
        \caption{Observed image (left panel) and true image (right panel) in the 2D image deblurring example.}
        \label{FIG:2DBlurTrue}
 \end{figure}

We take as our target image a $50 \times 50$ pixel grayscale image of geometric shapes so that $n = 2500$ in (\ref{EQ:Priors-X}). We blur the image by convolution with a Gaussian point spread function. The forward model $\B{A}$ and true image $\B{x}$ are created using the \texttt{Regularization Tools} package \cite{Hansen2007}. The data are generated by adding Gaussian noise with variance $0.01^2\|\B{Ax}\|_{\infty}^2$. Figure~\ref{FIG:2DBlurTrue} displays the target image and the noisy data.

\begin{figure}[tb]
    \centering
    \includegraphics[scale = 0.3]{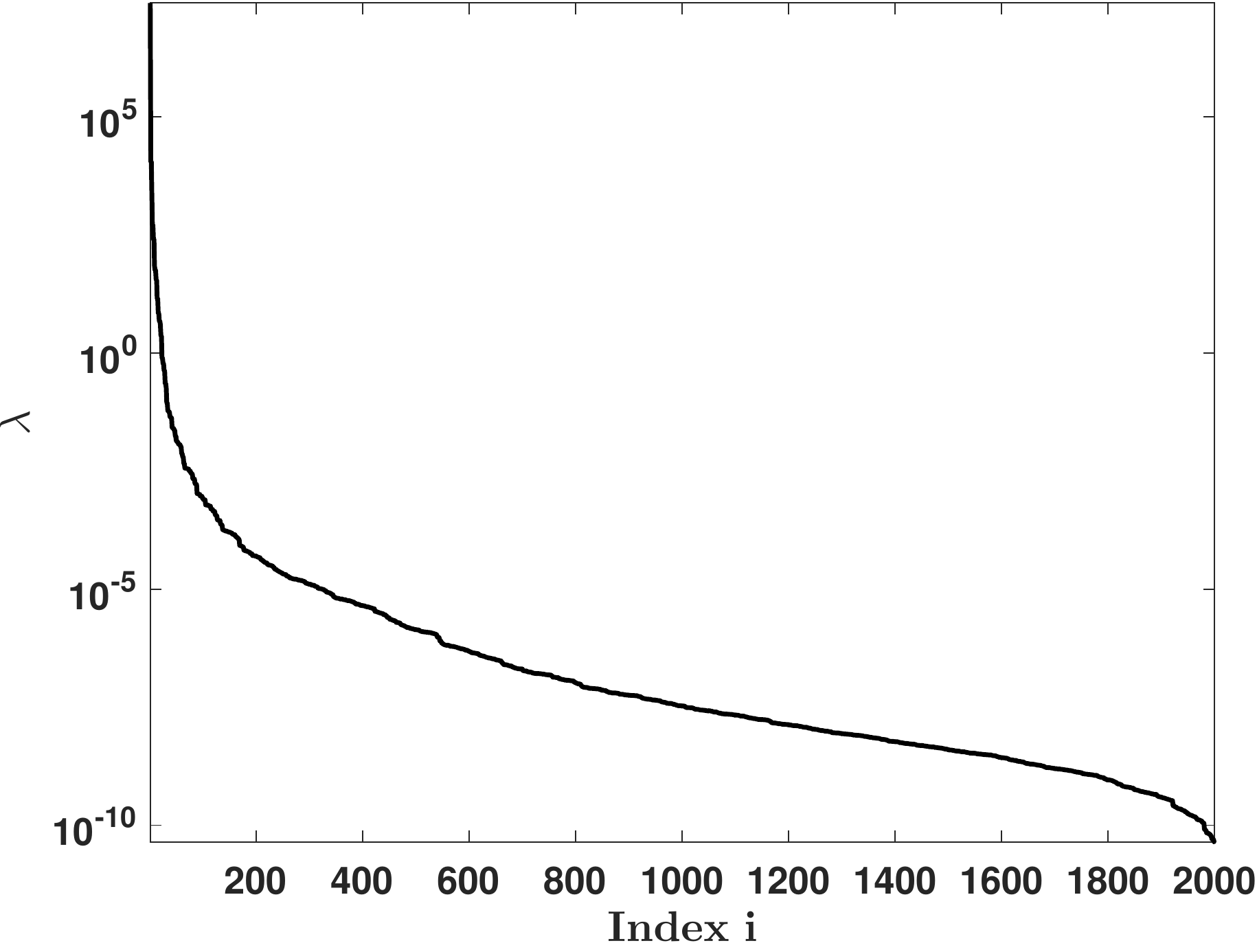}
    \includegraphics[scale = 0.3]{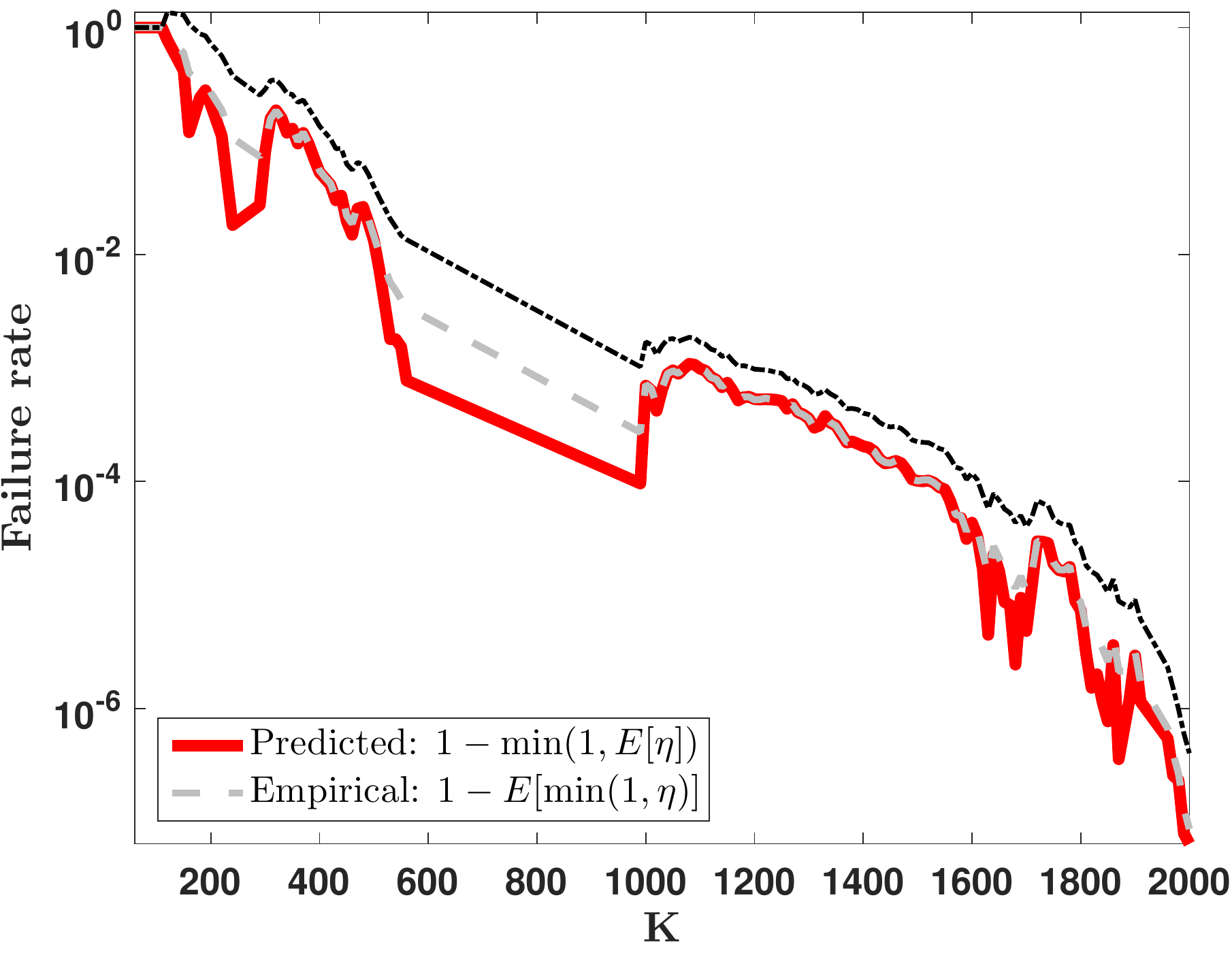}
    \caption{(left) Eigenvalues of $\B{H}$ in the 2D image deblurring example. (right) Predicted and empirical failure rates for different truncation levels $K$. The empirically determined upper $0.975$ quantile is given by the dashed line; the lower quantile often attained zero values, so it is not displayed.}
    \label{FIG:2Dfailures}
\end{figure}

%\paragraph{Choice of priors}
We model smoothness on $\B{x}$ {\em a priori} by taking $\B{L} = -\Delta + \delta \B{I}$ in \eqref{EQ:Priors-X}, where $-\Delta$ is the discrete Laplacian and $\delta$ is a small constant to ensure positive definiteness \cite{KaiSom05}. For the prior and noise precision parameters, we assign a vague Gamma prior, $\text{Gamma}(0.1, 0.1)$, which approximates the scale invariant objective prior while maintaining conditional conjugacy. We compute the eigenvalues of the prior preconditioned Hessian matrix $\B{H}$ via SVD to determine an appropriate cutoff. Figure~\ref{FIG:2Dfailures} (discussed further below) indicates rapid decay within the first few eigenvalues, followed by a smoother decay, and another sharp decay. We use the first $k = 500$ eigenvalues of the matrix to construct our low-rank approximation. We analyze below the effect of truncation level on the acceptance rate of the sampling algorithm. For comparison, we also compute a low-rank approximation using the Randomized SVD approach described in Subsection~\ref{ss_randsvd}. %To ensure the same number $\ell = k + p$ of eigenvalues are retained with this approach, we specify a cutoff of $k=480$  and an oversampling parameter $p=20$.

%\begin{figure}[tb]
%    \centering
%\includegraphics[scale = 0.45]{blur-spectra.eps}
%    \caption{Eigenvalues of $\B{H}$ for different blurring operators.}
%    \label{FIG:2Dspectra}
%\end{figure}

\paragraph{Convergence and UQ metrics}
We implement a Metropolis-Hastings-within-Gibbs algorithm in which Step \ref{EQ:Gibbs-Conditional-x} of Algorithm \ref{ALG:Basic-BlockGibbs} is substituted with our proposed low rank independence sampler (LRIS) presented in Section \ref{SEC:Proposal}. Three different chains are run in parallel, with each chain initialized by drawing $\B{x}, \mu$, and $\sigma$ randomly from their prior distributions. Each chain is run for $N=50,000$ iterations with the first $25,000$ iterations discarded as a burn-in period. For comparison, the Gibbs sampler is implemented identically to the low-rank procedure with three independent chains run in parallel with widely dispersed initial values. All simulations are done in \texttt{MATLAB} running on OS X Yosemite (8GB RAM, Intel Core i5 2.66GHz processor).
\begin{figure}[tb]
    \centering
    \includegraphics[scale= 0.25]{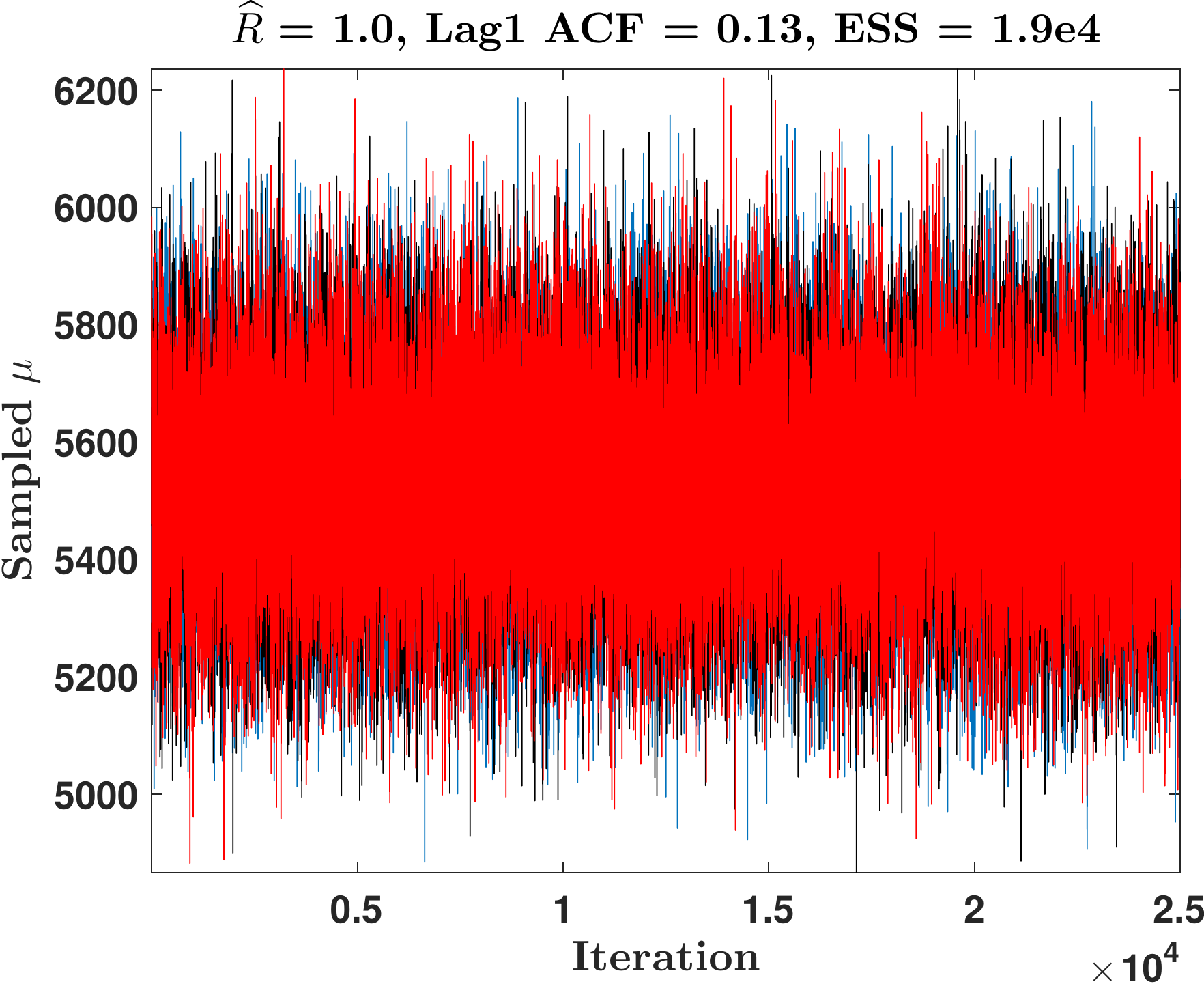}
    \includegraphics[scale= 0.25]{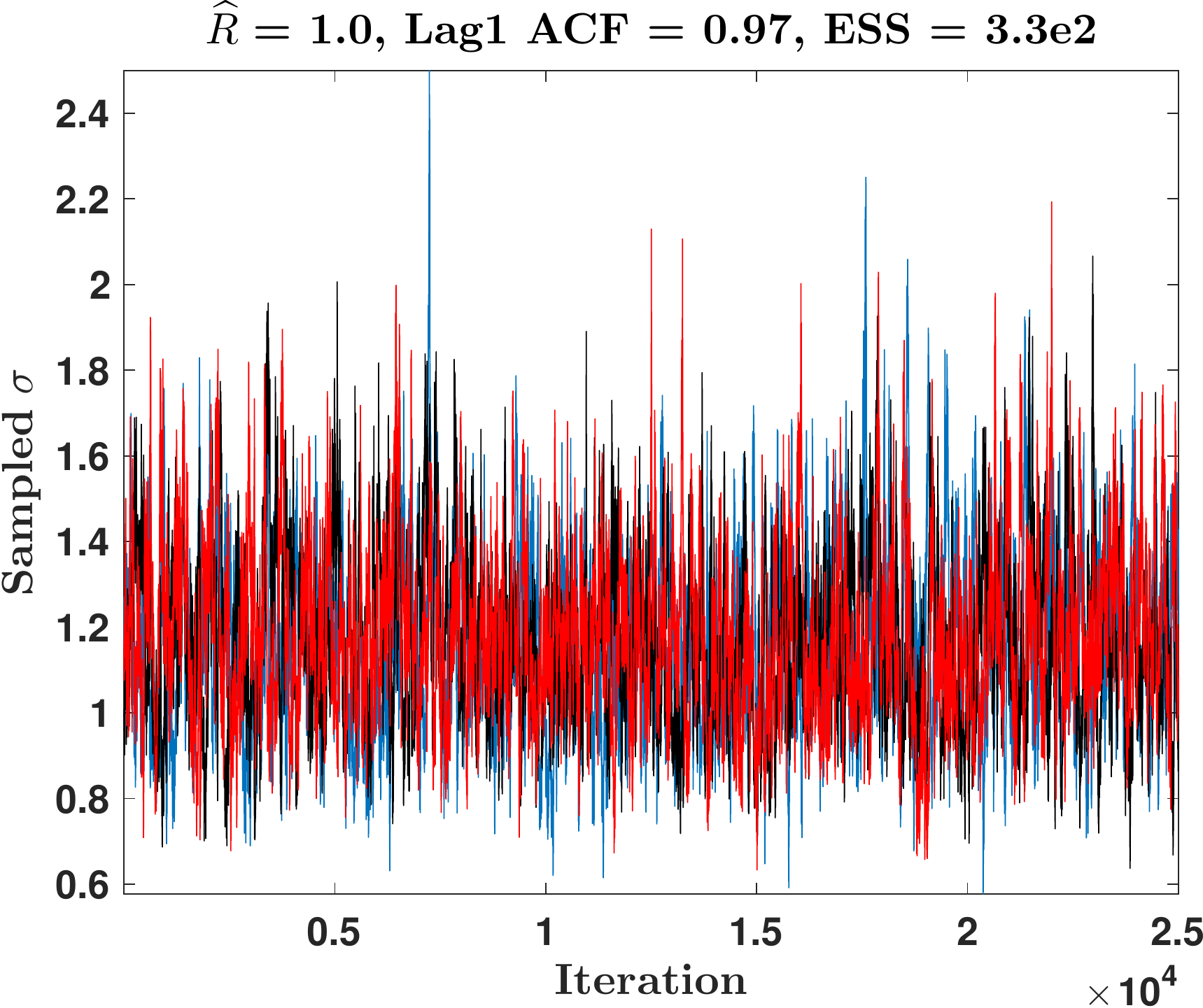}\\
     \includegraphics[scale= 0.25]{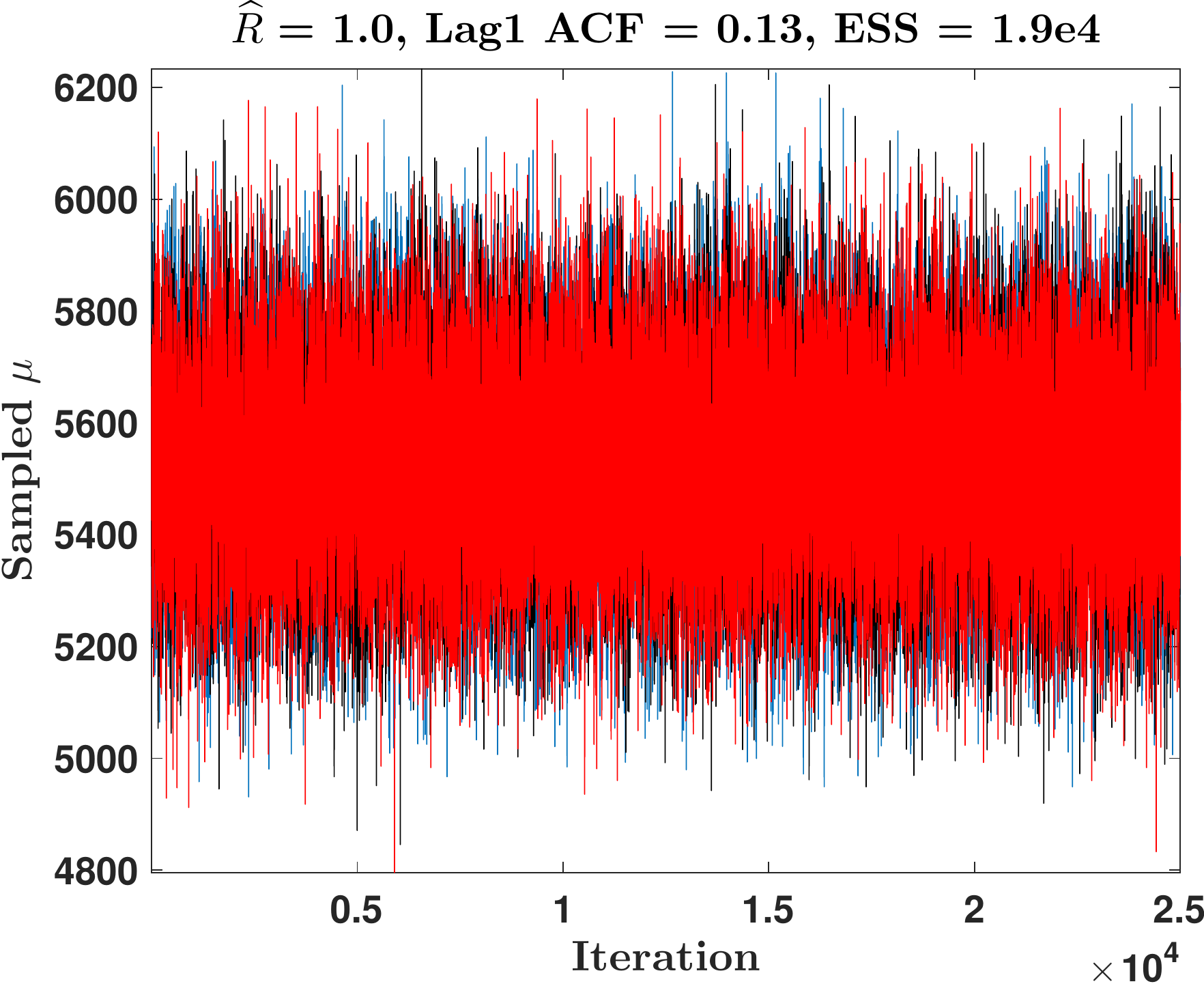}
    \includegraphics[scale= 0.25]{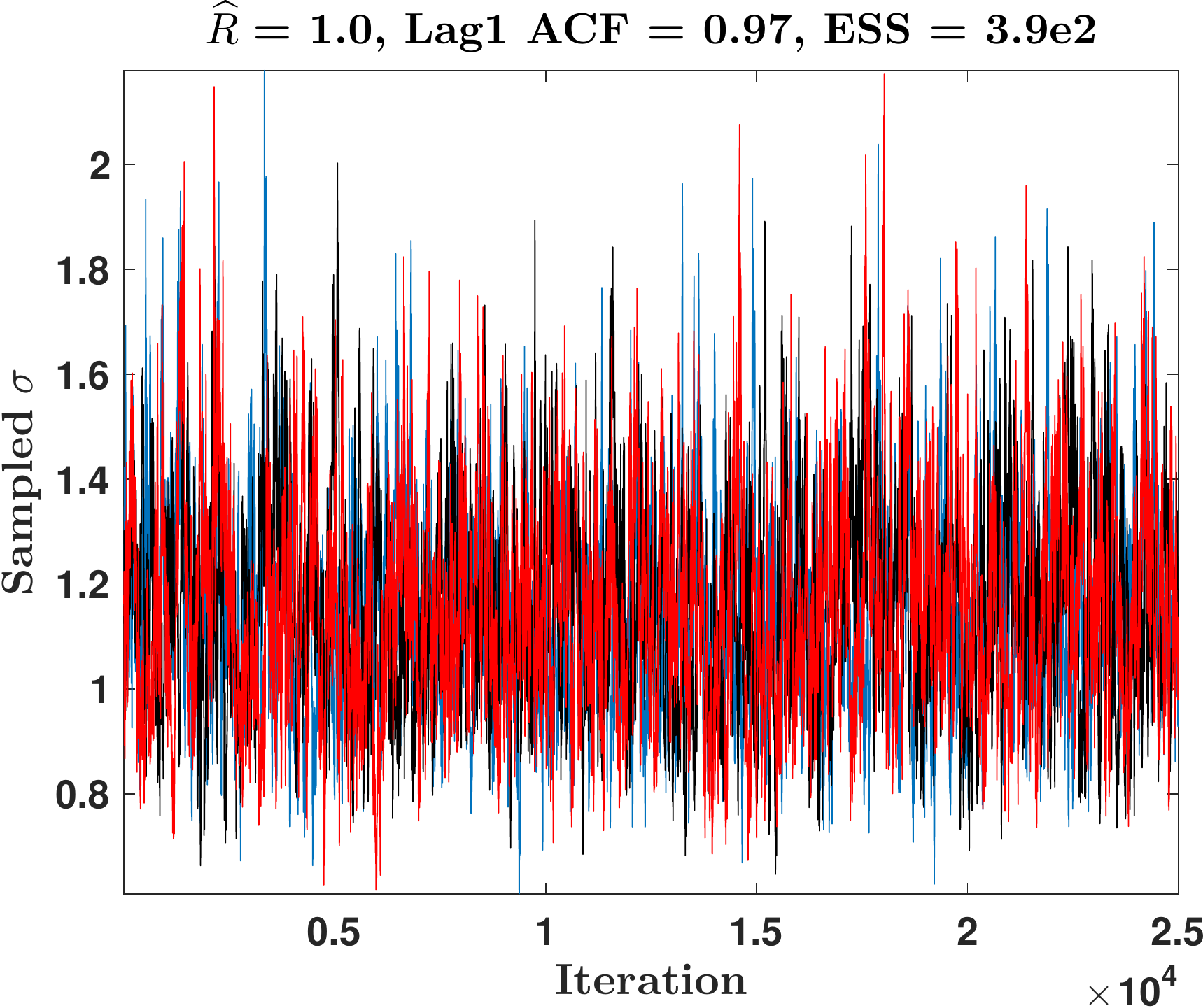}
    \caption{Trace plots of $\mu$ (left) and $\sigma$ (right) obtained from the MCMC output of both block Gibbs sampling (top row) and the LRIS-based algorithm (bottom row), where each of the three colors represents a different chain. The potential scale reduction factors ($\widehat{R}$) are displayed above each plot, along with the lag 1 autocorrelation coefficients and effective sample size (ESS) estimated from one chain each under both approaches.}
    \label{fig:BgMwgTPs}
\end{figure}

\begin{figure}[tb]
    \centering
    \includegraphics[scale= 0.25]{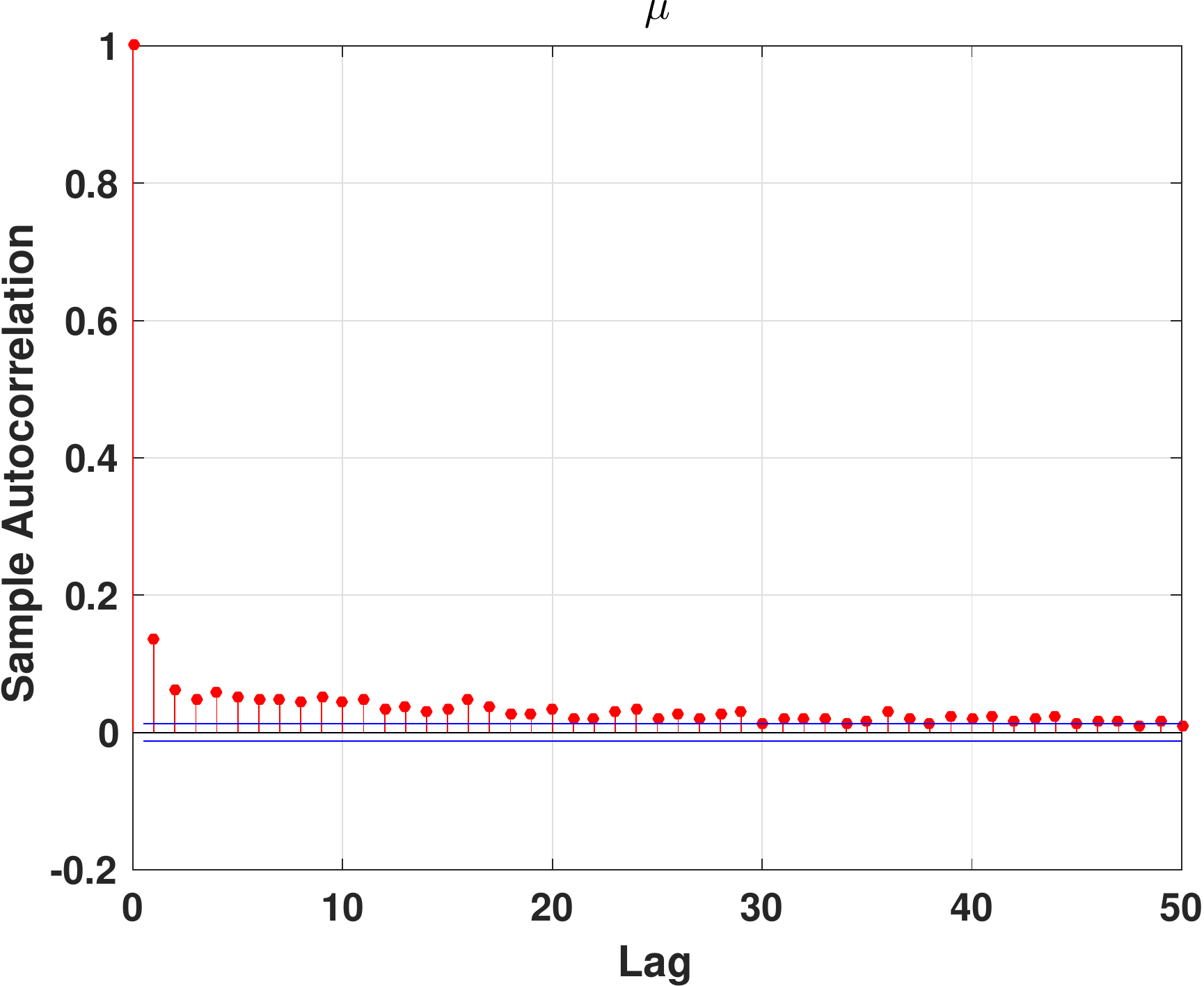}
    \includegraphics[scale= 0.25]{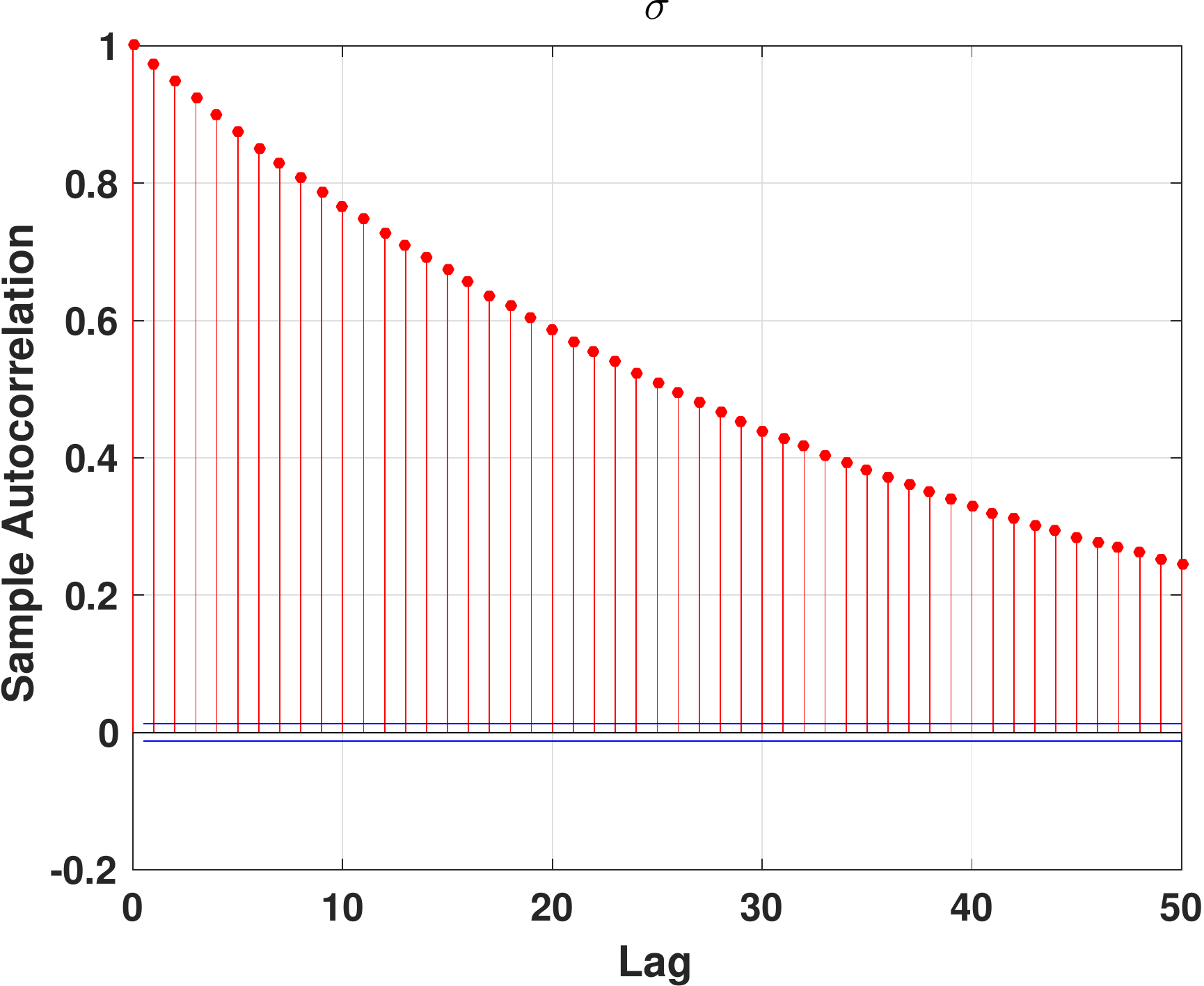}\\
     \includegraphics[scale= 0.25]{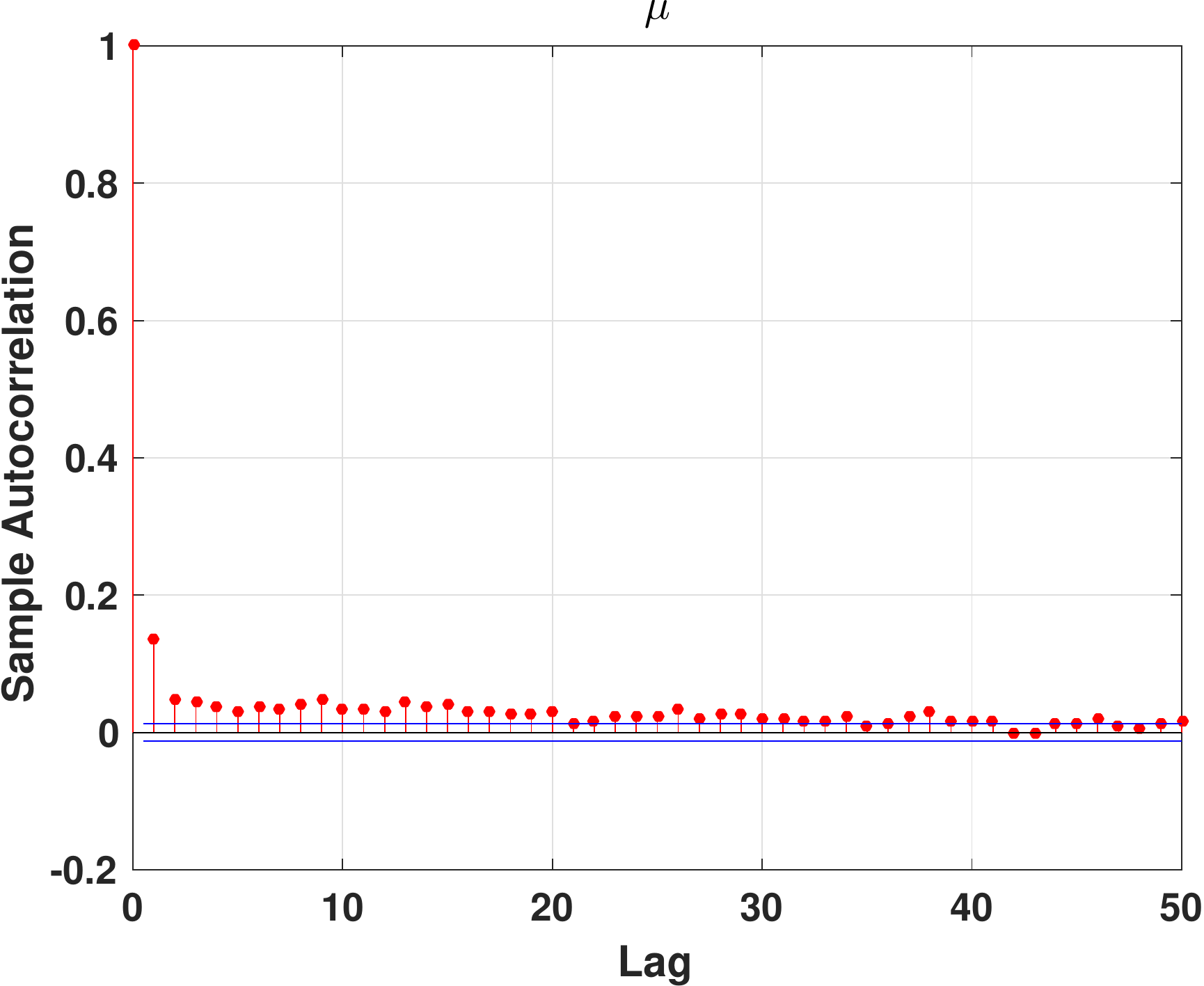}
    \includegraphics[scale= 0.25]{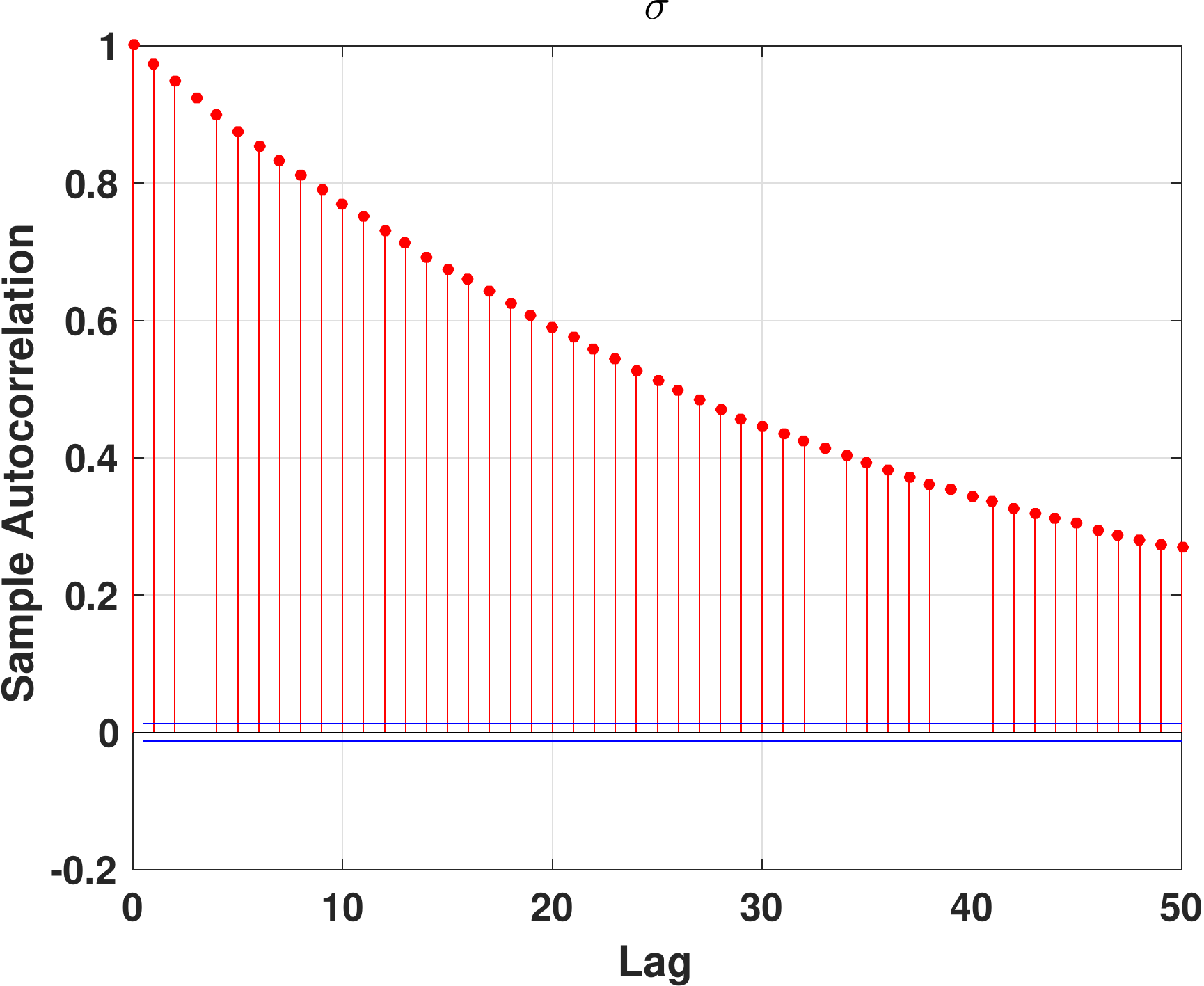}
    \caption{Estimated autocorrelation functions of $\mu$ (left) and $\sigma$ (right) obtained from one of the chains each under block Gibbs sampling (top row) and our proposed LRIS algorithm (bottom row).}
    \label{fig:BgMwgACF}
\end{figure}

\begin{figure}[tb]
    \centering
    \includegraphics[scale= 0.325]{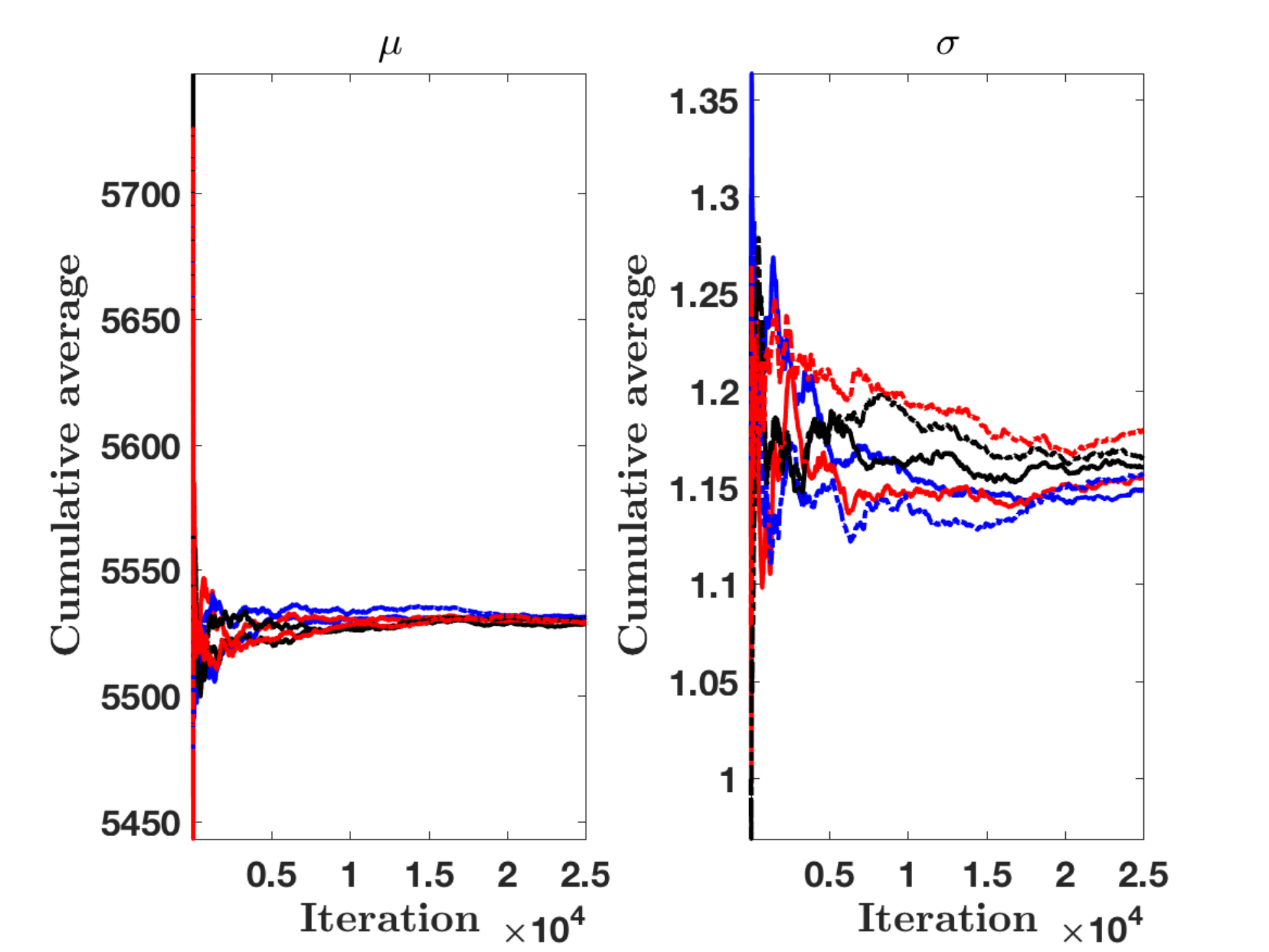} %Arvind: is this the right plot
    \includegraphics[scale= 0.325]{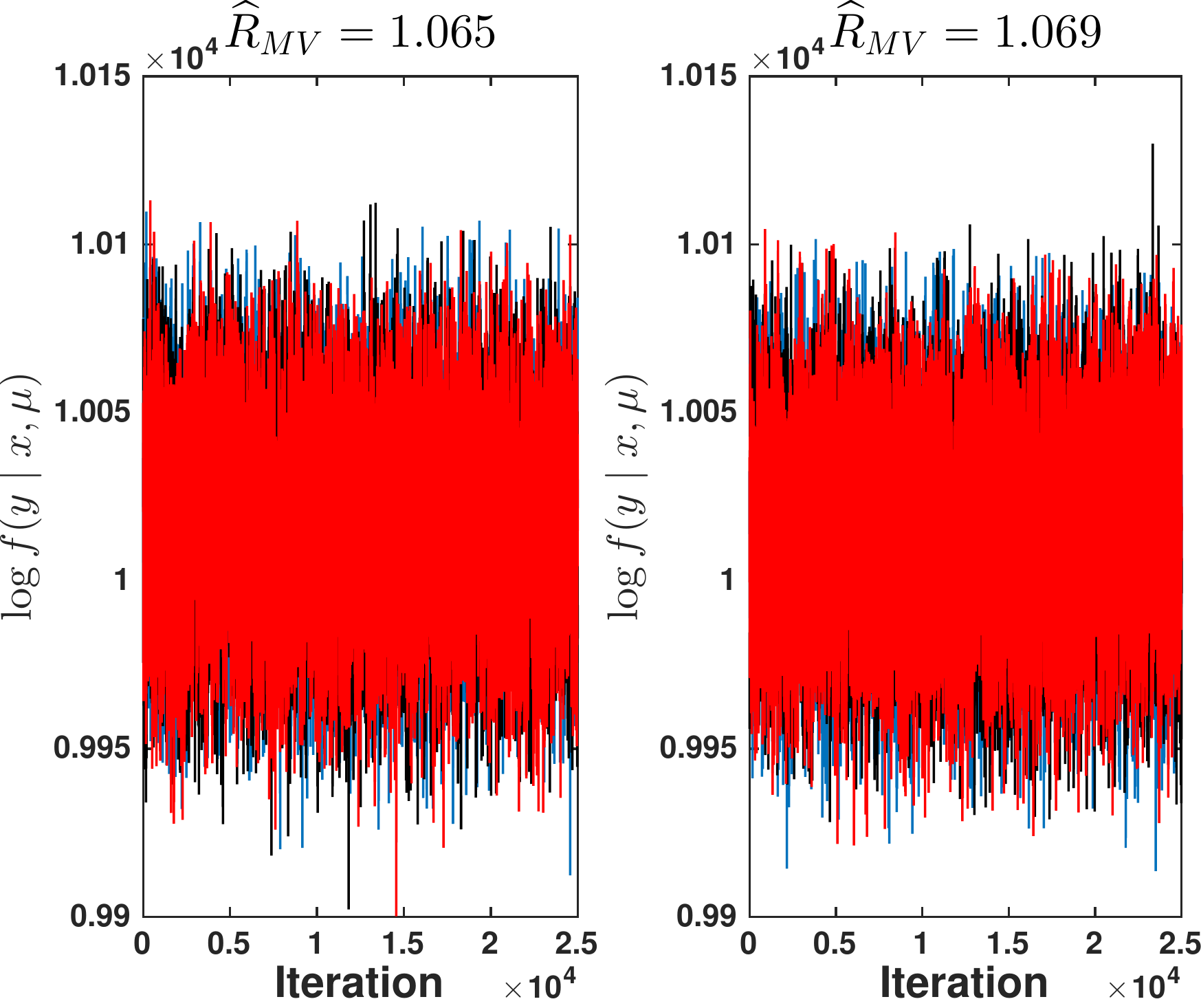}
    \caption{Cumulative averages of the $\mu$ chains (far left) and $\sigma$ chains (middle left) obtained from the MCMC output in the 2D deblurring example. The dotted lines represent the three chains from block Gibbs, the solid lines correspond to the LRIS algorithm output. Plot of the data-misfit part of the log-likelihood values calculated from the MCMC output of both block Gibbs sampling (middle right) and our low-rank independence sampling algorithm (far right), where each of the three colors corresponds to a different chain. The multivariate potential scale reduction factors are displayed above.}
    \label{fig:BgMwgPAvg}
\end{figure}

Figures \ref{fig:BgMwgTPs} and \ref{fig:BgMwgACF} display trace plots and autocorrelation functions, respectively, for the last $25,000$ iterations of the $\mu$ and $\sigma$ chains for both ordinary block Gibbs sampling and our proposed algorithm. As is known to occur with block Gibbs sampling in high-dimensional linear inverse problems \cite{bardsley2012mcmc}, we observe near independence within the $\mu$ chains and strong autocorrelation in the $\sigma$ chains. Despite the high autocorrelation, we still are able to achieve approximate convergence and a sufficient effective sample size (ESS) from the $\sigma$ chains by running each chain long enough. By combining the three independent chains after approximate convergence, we effectively triple the ESS and thus the number of independent pieces of information available about the target posterior. Thinning the chains to, e.g., every $10${th}, $50${{th}}, or $100${th} draw would dramatically reduce the autocorrelation of the chains. However, it was argued by Carlin and Louis \cite{CarlinLouis09} that such thinning is not necessary and does not improve estimates of quantities of interest. Figure \ref{fig:BgMwgPAvg} illustrates the approximate convergence of the ergodic averages $\widehat{\mu}_{(n)} = n^{-1}\sum_{t=1}^n\mu_{(t)}$ and $\widehat{\sigma}_{(n)} = n^{-1}\sum_{t=1}^n\sigma_{(t)}, ~n= 1, \ldots, 25000,$ despite the high autocorrelation of the $\sigma$ chain. The limiting values from both approaches closely agree.

%As pointed out by \cite{BesagEtAl95}, \cite{GelfandSahu99}, and \cite{EberlyCarlin00}, weakly identified subsets of parameters in the full posterior are not necessarily problematic when the parameters of interest are more strongly identified in a so-called embedded posterior distribution. In our case, we are less concerned with $\sigma$ and more concerned with the {\em a posteriori} information we have about the estimand of interest, $\B{x}$.
While assessing convergence of high-dimensional parameters is more difficult than for scalar quantities, we can track realized values of the data-misfit part of the log-likelihood as a proxy for monitoring convergence. These realizations also should settle down as the chain approaches the target distribution. Figure \ref{fig:BgMwgPAvg} displays these plots for both algorithms along with the multivariate PSRFs. Again, we see consistency between ordinary block Gibbs and our own approach, as well as approximate convergence according to the rule of thumb that the PSRF should be less than or equal to approximately 1.1 \cite{GelmanEtAl14}.
%%%Arvind: I removed the loglikelihood plots
%\begin{figure}[tb]
%    \centering
%    \caption{P}
%    \label{fig:BgMwgLkhd}
%\end{figure}

The advantage of using the LRIS approach is clear in Table \ref{tab:CES}, which displays the total wall time to complete the $50,000$ MCMC iterations for both block Gibbs and our proposed low-rank sampling approach. Table \ref{tab:CES} also displays the {\em cost per effective sample} (CES) for one of the $\sigma$ chains obtained under both algorithms as well as the Randomized SVD approach. CES is a measure of the average computational effort required between effectively independent draws. The LRIS approach yields a 76\% reduction in computation time compared to the standard block Gibbs sampler, along with an approximate 80\% reduction in computational effort between independent draws of $\sigma$. The average acceptance rate over the three chains using our low-rank proposals is 98\%, for both ``exact" and Randomized SVD. The acceptance rate versus rank is discussed further below.
\begin{table}[!ht]
    \centering
    \begin{tabular}{l | c c}
        Algorithm & Wall Time (s) & CES for $\sigma$ \\
        \hline
        Block Gibbs & 27907 & 84.30 \\
        LRIS & 5134 & 13.20\\
        LRIS (Randomized SVD) & 5307 & 13.86
    \end{tabular}
    \caption{Total wall time to complete $50,000$ MCMC iterations under block Gibbs sampling and the LRIS approach for the 2D deblurring example, along with the estimated cost per effective sample (CES) for one of the $\sigma$ chains in each case.}
    \label{tab:CES}
\end{table}

\begin{figure}[tb]
    \centering
\includegraphics[scale = 1.75]{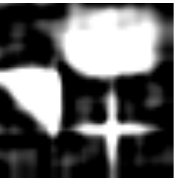}
\includegraphics[scale = 1.75]{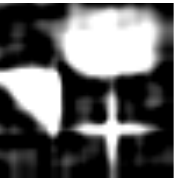}
    \caption{Approximate posterior mean images obtained from block Gibbs sampling (left) and the proposed low-rank sampling algorithm (right) for the 2D deblurring example.}
    \label{FIG:2Dpostmean}
\end{figure}

\begin{figure}[tb]
    \centering
\includegraphics[scale = 0.25]{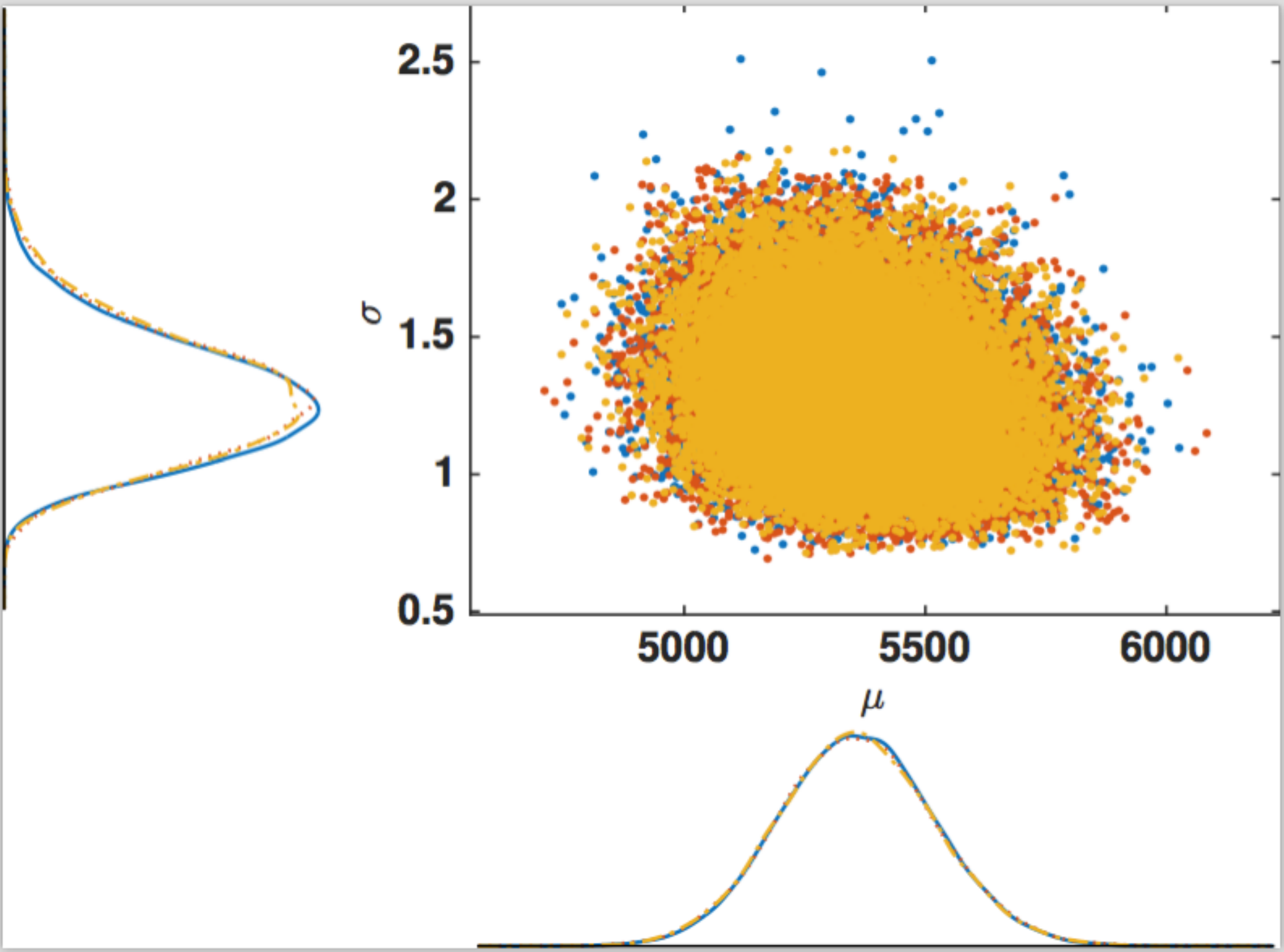} \hspace{8pt}%\hfill
\includegraphics[scale = 0.25]{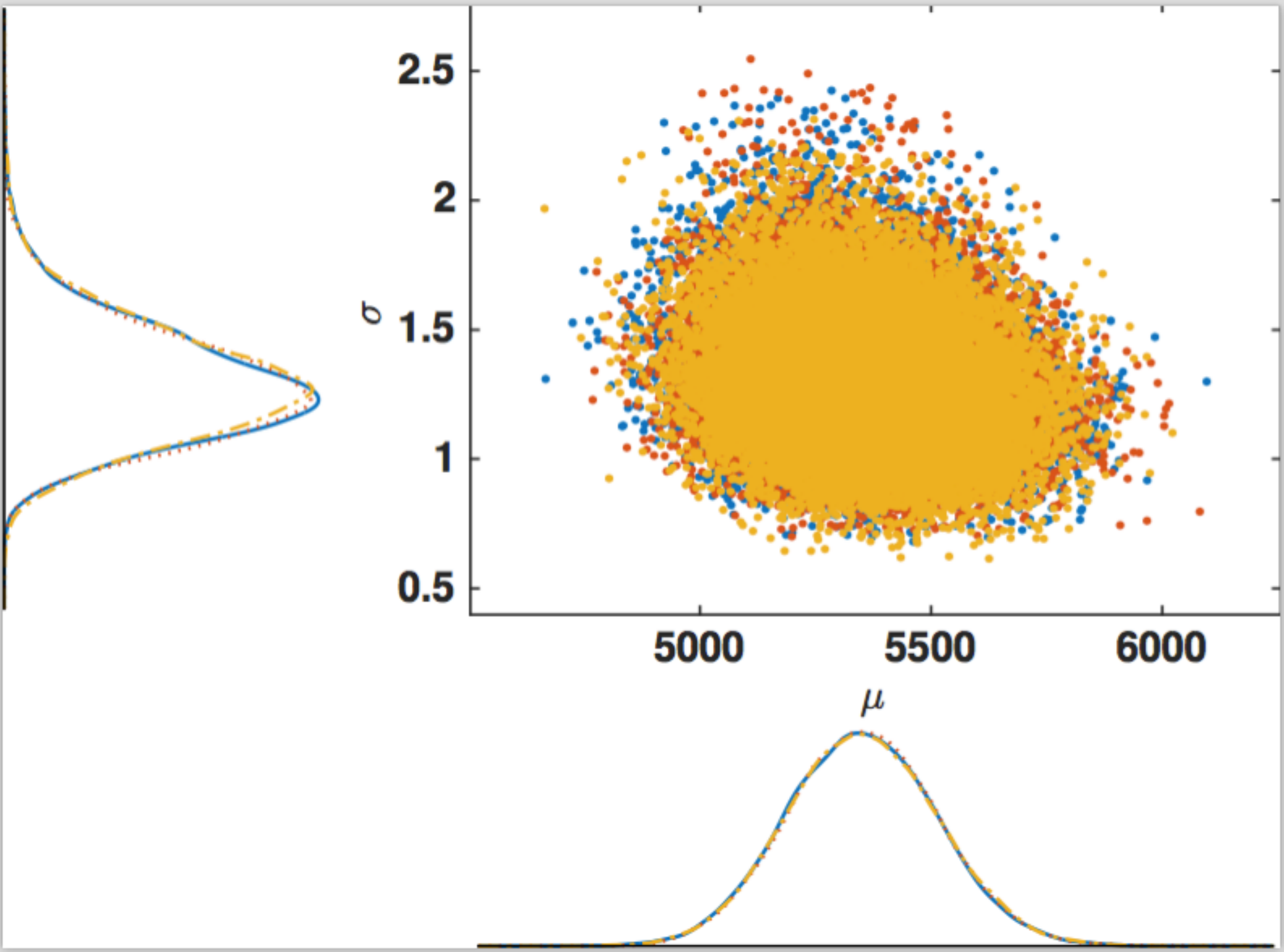}
    \caption{Scatter plots of the $(\mu, \sigma)$ realizations obtained from block Gibbs (left panel) and the low-rank approach (right panel) for the 2D deblurring example, where the three different colors correspond to different chains. The smoothed marginal posterior densities are displayed in the margins of the plots.}
    \label{FIG:2DmuSigDens}
\end{figure}

We attain this dramatic reduction in computational effort without sacrificing the quality of posterior inferences, as evident in Figure \ref{FIG:2Dpostmean}. This Figure displays the approximate posterior means of $\B{x}$ from both block Gibbs and our low-rank approach. The estimators we obtain with Randomized SVD are similar and hence omitted. We show also the $(\mu, \sigma)$ scatterplots and approximate marginal densities obtained from both algorithms in Figure \ref{FIG:2DmuSigDens}, again showing agreement. The strong Bayesian learning that occurred about these parameters is evident in Supplementary Figure \ref{fig:2Dmarginals}. Table \ref{tab:MapMean2Derrors} gives the relative errors to quantify the quality of the reconstructions. We observe nearly identical solutions under both MCMC approaches, both graphically and quantitatively.
\begin{table}[!ht]
    \centering
    \begin{tabular}{l | c c}
        Estimator & RE \\ %& RMSE \\
        \hline
        Posterior Mean (Block Gibbs) & 0.4453 \\ %& 0.4900 \\
        Posterior Mean (LRIS) & 0.4455 %& 0.4907 \\
        % Posterior Mean (Randomized SVD) & 0.4413 \\ %& 0.4861 \\
        %MAP & 0.5285 \\ %& 0.5822
    \end{tabular}
    \caption{Relative error (RE) of the estimates for the 2D image deblurring example.}
    \label{tab:MapMean2Derrors}
\end{table}
%In addition to posterior means, we consider also the posterior mode (MAP), a common point estimator that corresponds to deterministic minimization of the objective function. The MAP estimate, obtained with a block coordinate descent on the negative log-likelihood of the joint posterior distribution, is displayed in Figure \ref{FIG:2Dpostmean} and the reconstruction error quantified in Table \ref{tab:MapMean2Derrors}.  In this case, the posterior mean produces a more faithful reconstruction of the solution than the MAP, which is clearly overly smooth. This is not surprising, since the posterior mean is the Bayes estimator under squared-error loss (the most common loss function for estimation) \cite[Ch. 4]{LehmannCasella98}. This is one reason why Carvalho et al.~\cite{CarvalhoEtAl10} argue in favor of the posterior mean for point estimation instead of the mode. %It was observed also by \cite{FoxNicholls01} that the MAP estimator can often be misleading.

\paragraph{Acceptance rate versus rank} To explore the effect of the retained number of eigenvalues on the acceptance rate for our algorithm, we estimate the predicted and empirical acceptance rates of sampling from the proposal distribution, over a range of truncation levels, at a given state of the chain. We fix the state by initializing $(\B{x}, \mu, \sigma)$ as the last sample from one of the chains obtained from the LRIS. At each truncation level $k$, we compute the expected value of the acceptance ratio using Theorem~\ref{p_expect}. We draw $2,000$ samples from the proposal distribution and compute the acceptance ratio of each using Proposition~\ref{p_mh_ratio2}. From these we estimate the empirical failure rates to compare with their expected values as the truncation level increases. Figure~\ref{FIG:2Dfailures} displays the results. The close agreement between the predicted and empirical acceptance rates support the theoretical results in Section~\ref{SEC:Proposal}.

\subsection{CT Image Reconstruction}
Computed x-ray tomography (CT) is a common medical imaging modality in which x-rays are passed through a body from a source to a sensor along parallel lines indexed by an angle $\omega$ and offset $y$ with respect to a fixed coordinate system and origin. The intensities of the rays are attenuated according to an unknown absorption function as they pass through tissue. The attenuated intensity $I$ is recorded while the lines are rotated around the origin so that $I(S) = I(0)\exp\{-\int_0^S \alpha(x(s)) \, \dx s\}$, where $s=0$ is the source of the x-ray, $s=S$ is the receiver location, $x(\cdot)$ indicates the line position, and $\alpha$ is the absorption function. The observed data are a transformation of the intensities, yielding the Radon transform model for CT \cite{KaiSom05, Bardsley2011},
%\begin{equation}\label{eqn:RadonCT}
  $z(\omega, y) = \int_{L(\omega, y)} \alpha(x(s)) \, \dx s$,
%\end{equation}
where $L(\omega, y)$ is the line along which the x-ray passes through the body. The inverse problem is to reconstruct the absorption function, which provides an image of the scanned body. Discretization of the integral yields the model in~\eqref{EQ:StochasticModel}. This is typically an underdetermined system with infinitely many solutions, resulting in an ill-posed inverse problem.

Our target image is the Shepp-Logan phantom \cite{SheppLogan74}. The forward model is implemented in \texttt{MATLAB} on the same computer as in Subsection \ref{subsec:2d} with code available online \cite{BardsleyCode}. The data are simulated by adding Gaussian noise with variance $0.01^2\|\B{Ax}\|_{\infty}^2$. The target $\alpha$ is discretized to a relatively fine grid of size $128 \times 128$ so that $\dim(\B{x}) = 16,384$. We suppose that the data are observed over lines and angles such that $\dim(\B{b}) = 5000$. Thus, $\text{rank}(\B{A}) = 5000 \ll \dim(\B{x})$, guiding our choice of eigenvalue truncation in the low-rank approximation to $\B{H}$. An approximate eigendecomposition of the prior preconditioned Hessian is computed using Randomized SVD with $\ell = 5000$, as discussed in Section \ref{SEC:Proposal}, since computing the ``exact" SVD is considerably more expensive. As in Subsection \ref{subsec:2d}, we take $\B{L} = -\Delta + \delta \B{I}$.

For the nuisance parameters, we use a weakly informative prior \cite{Gelman06}, namely the proper Jeffreys prior proposed by Scott and Berger \cite{ScottBerger06}. For convenience, we parameterize the model in terms of variance components instead of precisions, $\tau^2 := \sigma^{-1}$ and $\kappa^2 := \mu^{-1}$. Then the proper Jeffreys prior on $(\kappa^2, \tau^2)$ is\footnote{We write `proportional to' ($\pi(\theta) \propto g(\theta)$) for proper priors to indicate that $\pi(\theta) = cg(\theta)$, where $c^{-1} = \int g(\theta) d\theta < \infty$ uniquely determines the density. However, the normalizing constant for an improper prior does not exist, so the prior is not unique. In the scale invariant case, any prior $\pi(\kappa^2) = a/ \kappa^2$ for $a \neq 0$ works. Since $a$ is arbitrary, we simply set it equal to $1$ for convenience and take $\pi(\kappa^2) = 1/\kappa^2$. See \cite[Chapter~3]{Berger85}.}
\begin{eqnarray}\label{EQN:Scott-Berger}
\begin{aligned}
    \pi_{SB}(\kappa^2, \tau^2) &= (\kappa^2 + \tau^2)^{-2}\\
        &= (\kappa^2)^{-1}(1 + \tau^2/\kappa^2)^{-2}\times (\kappa^2)^{-1}\\
        &\equiv \pi(\tau^2 \mid \kappa^2)\pi(\kappa^2),
        \end{aligned}
\end{eqnarray}
so that the scale invariant prior is used for $\kappa^2$ while scaling $\tau^2$ by the data level variance, as advocated by Jeffreys \cite{Jeffreys61}. The implementation of this prior as a modification to Algorithm \ref{ALG:Basic-BlockGibbs}, presented in Appendix \ref{app:pj_prior}, is similar to the approach of \cite{BrownEtAl17}. Section \ref{sec:Priors} of the Supplementary Material contains further discussion of prior specification for the nuisance parameters.

We simulate three Markov chains using our proposed LRIS approach for $40,000$ iterations (average acceptance rate $\sim100\%$). Each chain is initialized with values drawn randomly from the prior. We thin the chains by retaining every $50\text{th}$ draw to reduce the autocorrelation, making it easier to diagnose convergence. We discard the first $400$ draws of the thinned chains as a burn-in period. Trace plots and autocorrelation plots are used to verify approximate convergence of the chains. Relevant diagnostic plots are displayed in Supplementary Figures \ref{FIG:CTtrace}, \ref{FIG:CTautocorr}, and \ref{fig:CTcumavgs}. The total computation time for our sampling approach is $197,517$ seconds, or about $55$ hours. This is noteworthy since the algorithm repeatedly updates a large, nontrivial covariance matrix and samples an approximately sixteen-thousand dimensional Gaussian distribution 40,000 times. An ordinary block Gibbs sampler is simply not feasible for this problem.

%\begin{figure}[tb]
%\centering
%\includegraphics[scale=0.25]{CT-trace-kappa2-THIN.eps}
%\includegraphics[scale=0.25]{CT-trace-upsilon-THIN.eps}
%\includegraphics[scale=0.25]{CT-trace-xpixel-THIN.eps}
%\caption{Trace plots for the thinned chains in the CT image reconstruction example. Left: noise variance $\kappa^2$. Center: Variance ratio $\upsilon = \tau^2 / \kappa^2$. Right: A randomly chosen pixel of the image $\B{x}$.}
%\label{FIG:CTtrace}
%\end{figure}
%
%\begin{figure}[tb]
%\centering
%\includegraphics[scale=0.3]{CT-autocorr-kappa2-THIN.eps}
%\includegraphics[scale=0.3]{CT-autocorr-upsilon-THIN.eps}
%\caption{Autocorrelation plots for the thinned chains in the CT image reconstruction example. Left: noise variance $\kappa^2$. Right: Variance ratio $\upsilon = \tau^2 / \kappa^2$.}
%\label{FIG:CTautocorr}
%\end{figure}

\begin{figure}[tb]
        \centering
    \includegraphics[scale = 0.35]{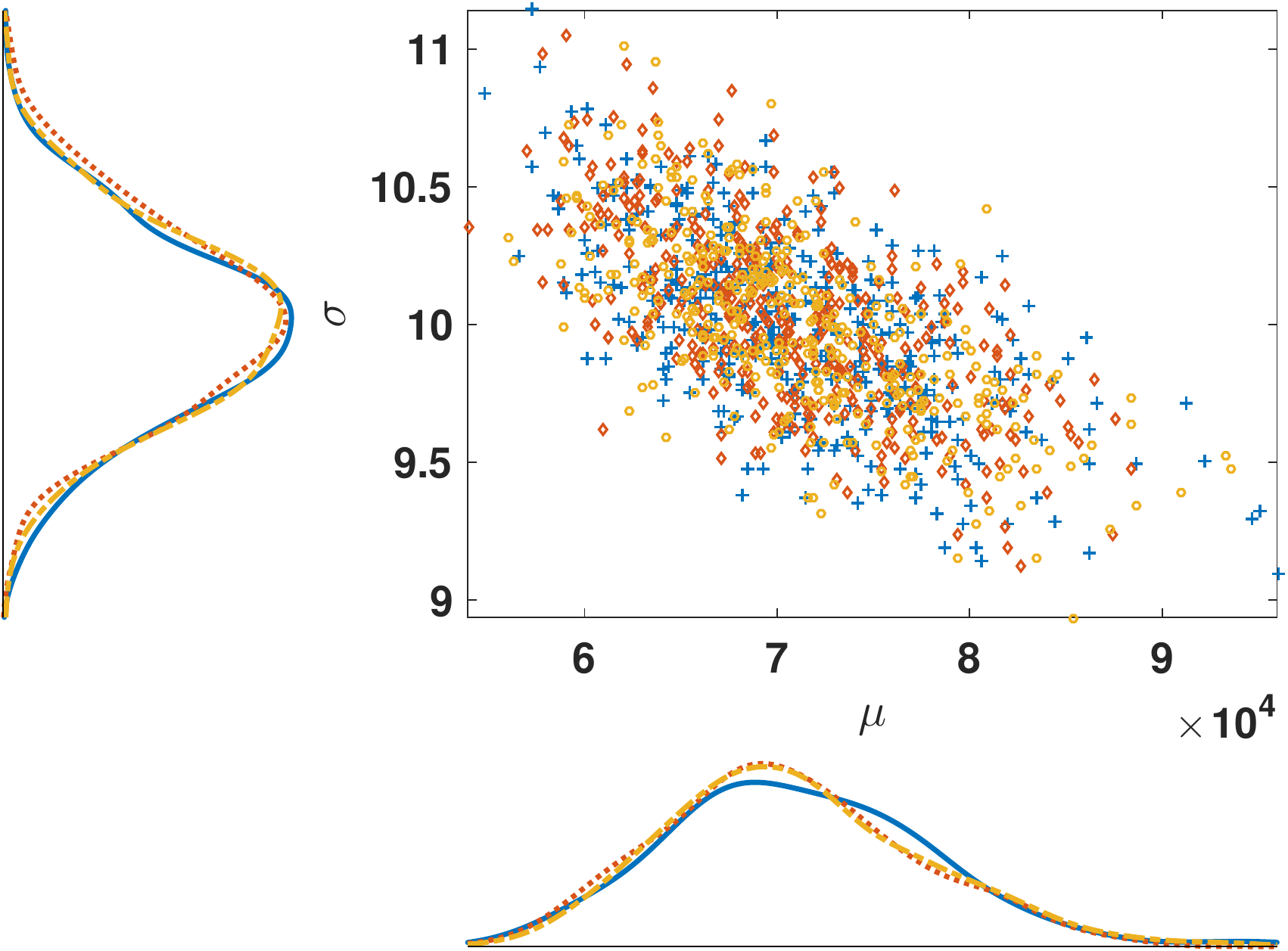}\hspace{8pt}
    \includegraphics[scale = 0.35]{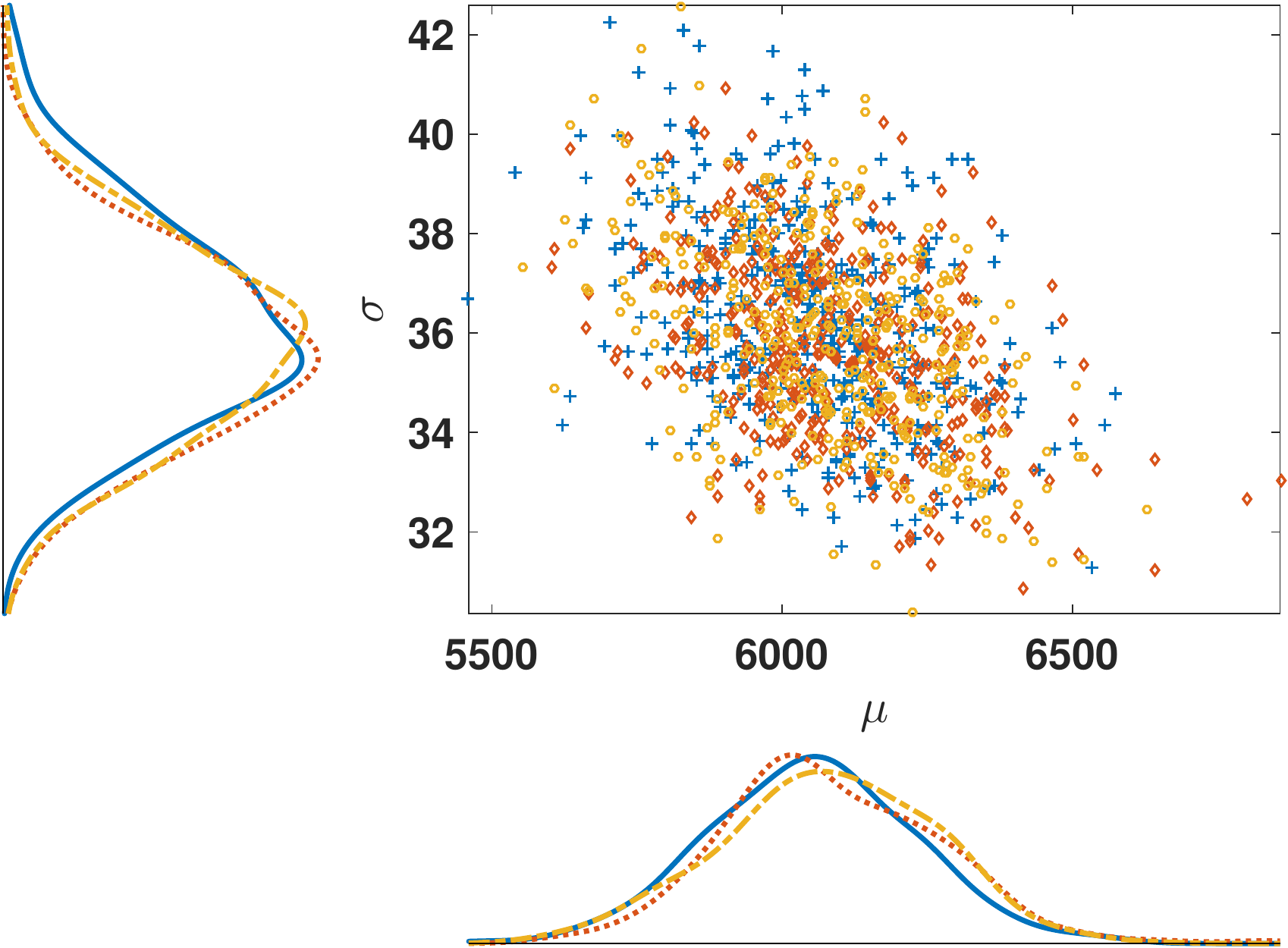}
        \caption{Estimated joint densities for the precision parameters using the proper Jeffreys prior (left) and the conjugate Gamma priors (right).}
        \label{FIG:CTscatters}
 \end{figure}

\begin{figure}[tb]
        \centering
    \includegraphics[scale = 1]{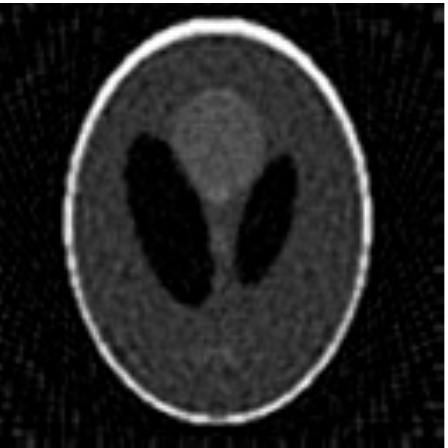}\hspace{8pt}
    \includegraphics[scale = 1]{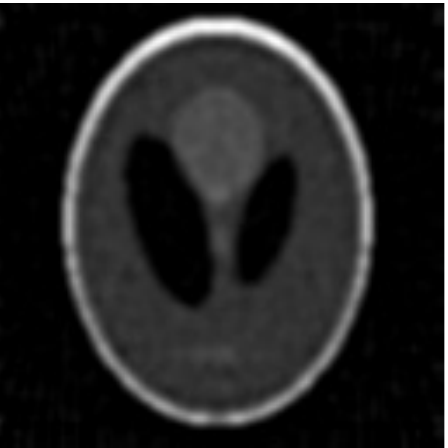}
        \caption{Posterior mean estimators of the true image in the CT image reconstruction example, using the proper Jeffreys prior (left) and the conjugate Gamma priors (right).}
        \label{FIG:CTMean}
 \end{figure}

We compare the results with samples obtained using the conjugate Gamma model with the same vague priors on $\mu$ and $\sigma$ as in Subsection \ref{subsec:2d}. Convergence diagnostics are displayed in Supplementary Figures \ref{fig:CTtracegamma} and \ref{fig:CTautocorrgamma}. Figure \ref{FIG:CTscatters} compares the approximate joint distributions for the precision parameters $(\mu, \sigma)$ under the conjugate model to the distribution based on the proper Jeffreys model, after back-transforming $\kappa^2$ and $\tau^2$. Here we see the effect of prior selection in that both $\mu$ and $\sigma$ tend to concentrate around different values, with much greater uncertainty in $\mu$ in the proper Jeffreys case. The differences between the two marginal posteriors of $(\mu, \sigma)$ affect the quality of the reconstructed images, displayed in Figure \ref{FIG:CTMean} and quantified in Table \ref{tab:MapsCTerrors}. This echoes Gelman's observation \cite{Gelman06} that even a supposedly noninformative prior on the hyper-precision can have a disproportionate influence on the results. In this case, using a weakly informative prior for $\sigma$ that depends on $\mu$ results in a higher quality reconstruction.

%We can see qualitatively the effect of both priors on the posterior mean estimator. The errors for the reconstructions are displayed in Table \ref{tab:MapsCTerrors}.

\begin{table}[th]
    \centering
    \begin{tabular}{l | c c}
        Scale parameter prior & RE \\ %& RMSE \\
        \hline
        Proper Jeffreys & 0.3270 \\ %& 0.0801 \\
        Conjugate Gamma & 0.4229 \\ %& 0.1037
    \end{tabular}
    \caption{Relative error (RE) of the estimates for the CT image reconstruction example.}
    \label{tab:MapsCTerrors}
\end{table}

%When approximating the posterior distribution via Markov chain Monte Carlo, the bottleneck is in repeatedly sampling high-dimensional Gaussian random variables.
These simulations show that by exploiting the low-rank structure of the preconditioned Hessian of the forward model, we are able to substantially reduce the computational burden compared to block Gibbs sampling. Even when the forward model is of full row rank, the results illustrate the potential for efficiency gains using our proposed LRIS approach, provided the system is underdetermined. Our approach using either proper Jeffreys or conjugate priors is much more feasible than a block Gibbs sampler, for which the computational demands can be prohibitively expensive. A comparison of the conventional conjugate Gamma priors to a weakly informative prior suggests that even a ``noninformative" prior may exert considerable influence on the results, despite strong Bayesian learning in the posterior.

In the Supplementary Material, we further consider the challenging problem of nuclear magnetic resonance (NMR) relaxometry. There, we demonstrate that our approach still produces within reasonable computation time a solution that is comparable to those obtained from deterministic iterative procedures such as conjugate gradient least squares. This is achieved with randomized SVD and without explicitly forming the forward operator $\B{A}$ in the LRIS algorithm.

%These simulation results illustrate the ability of fully Bayesian solutions to the inverse problem to perform competitively with deterministic solutions. Access to the full posterior distribution makes available more point estimators of the solution along with a set of plausible values of the regularization parameter. When approximating the posterior distribution via Markov chain Monte Carlo, the bottleneck is in repeatedly sampling high-dimensional Gaussian random variables. By exploiting the low-rank structure of the preconditioned Hessian of the forward model, we are able to substantially reduce the computational burden compared to block Gibbs sampling. Even when the forward model is of full row rank, we illustrate the potential for efficiency gains using the proposed low-rank proposal approach, provided the system is underdetermined.

\section{Discussion}\label{SEC:Discussion}
When approximating the posterior distribution via Markov chain Monte Carlo in the Bayesian Gaussian linear inverse problem, the bottleneck is in repeatedly sampling high-dimensional Gaussian random variables. Sampling from the joint posterior with standard MCMC is challenging due to the high dimensionality of the estimand, since drawing from the full conditional involves expensive operations with the covariance matrix.

In this work we propose a computationally efficient sampling algorithm which is well suited for a fully Bayesian approach in which the noise precision and the prior precision parameters are unknown and assigned prior distributions. Our proposed low rank independence sampler uses a proposal distribution constructed via low-rank approximation to the preconditioned Hessian. We show that the acceptance rate is high when the magnitudes of the discarded eigenvalues of the Hessian are small, a feature of severely ill-posed problems. When it is not obvious how to determine an appropriate truncation of rank due to the dependence of the known acceptance rates on other parameters in the model, we discuss how to adaptively determine the truncation level as part of the MCMC algorithm to find the minimal rank with high acceptance rates. We demonstrate both theoretically and empirically that the quality of the approximation is directly related to the acceptance rate of the sampler, as intuition would suggest. %Our approach is not limited by the choice of prior distribution for $\mu$ and $\sigma$, and the posterior simulation allows an exploration of plausible values for these marginalized over all plausible reconstructions $\B{x}$ given the data.
We illustrate our approach on several examples, demonstrating convergence as function of rank, as well as computational improvements to accessing the full posterior distribution.
%
%Our comparison of conventional conjugate Gamma priors on the nuisance parameters versus a more formally justified weakly informative prior shows that they can lead to substantially different marginal posterior distributions. This does not always have a dramatic effect on the reconstructed solution of interest \cite{CarlinLouis09}, but our experiment did produce a noticeable difference in the reconstruction. While this is likely due to the fact that the proper Jeffreys prior scales the prior variance by the data level variance instead of treating them independently, the precise nature of how such priors behave in linear inverse problems is, as far as we know, not yet fully understood. This is possibly an interesting avenue for future research.

%One limitation of our approach is that the prior precision parameter is correlated with the spatial unknowns $\B{x}$ in the MCMC algorithm, and this correlation increases as the mesh refinement increases \cite{agapiou2014analysis}. Methods to address this issue are proposed in~\cite{FoxNorton16,bardsley2016partially}. Combining our efficient sampler with these approaches is possibly an interesting avenue of future research.
A known issue with block Gibbs sampling in Bayesian inverse problems is the deterioration of the chains due to correlation between the hyperparameters and the estimand $\B{x}$ as the dimension of the problem increases \cite{bardsley2012mcmc}. Several approaches have been proposed to ameliorate this by breaking the dependence between the hyperparameters and $\B{x}$ in the algorithm. These include the one-block algorithm \cite{RueHeld05}, partially collapsed samplers \cite{VanDykPark08, bardsley2016partially}, noncentered parameterization (NCP) \cite{Papaspiliopoulos2003, Papaspiliopoulos2007}, and marginal then conditional (MTC) sampling \cite{FoxNorton16}. Noncentered parameterization is easily incorporated into our proposed approach without sacrificing gains in efficiency. (See the Supplementary Material for discussion of NCP and illustration of combining it with our low-rank sampler.) One-block sampling and MTC, on the other hand, require expressions for marginal densities that no longer hold when substituting the true full conditional of $\B{x}$ with an approximation, as well as approximation of determinants of large covariance matrices to which the results in this work are not directly applicable. Combining our proposed low-rank sampling approach with these algorithms is the subject of ongoing research, to appear in future work.

\section{Acknowledgements}
DAB is partially supported by Grants CMMI-1563435 and EEC-1744497 from the National Science Foundation (NSF). This material is based upon work partially supported by the NSF under Grant DMS-1127914 to the Statistical and Applied Mathematical Sciences Institute (SAMSI). The authors thank the Editors, an Associate Editor, and anonymous referees for comments and suggestions that improved this manuscript. The authors also would like to thank Duy Thai, Vered Madar, Johnathan Bardsley, and Ray Falk for useful conversations. Much of this work was done while the first author was a Visiting Research Fellow at SAMSI.

\appendix
\section{Proofs}\label{a_proofs}
%\textit{Proof of Proposition~\ref{p_mh_ratio1}}
%
%Let us consider
%\[ \frac{h(\B{x})}{g(\B{x})}  = \sqrt{\frac{\det(\condh)}{\det({\cond})}} \hat{w}(\B{x}),\]
%where we define
%\begin{equation*} \hat{w}(\B{x})  \> \equiv \> \exp\left( - \frac{1}{2}\| \B{x}-\xcond\|_{\cond^{-1}}^2 + \frac{1}{2}\| \B{x}-\xcondh\|_{\condh^{-1}}^2\right).
%\end{equation*}
%We have used the notation $\| \B{w}\|_{\B{M}} = \sqrt{\B{w}^\top \B{Mw}}$ which is a valid vector norm when $\B{M}$ is positive definite. With this definition, the  acceptance ratio simplifies to $\eta(\B{z} \mid \B{x}) = \hat{w}(\B{z})/\hat{w}(\B{x})$. We can further simplify $\hat{w}(\B{x})$
%\begin{align*}
%\log \hat{w}(\B{x})  =&  \> -\frac{1}{2}(\B{x}- \mu\cond\B{A}^\top\B{b})^\top\cond^{-1}(\B{x}- \mu\cond\B{A}^\top\B{b})  \\
%  & + \> \frac{1}{2}(\B{x}- \mu\condh\B{A}^\top\B{b})^\top\condh^{-1}(\B{x}- \condh\B{A}^\top\B{b})\\
%= & \>  -\frac{1}{2}\left( \B{x}^\top\cond^{-1}\B{x}  + \mu^2\B{b}^\top\B{A}\cond\B{Ab} \right)  \\
%&\> \qquad + \> \frac{1}{2}\left( \B{x}^\top\condh^{-1}\B{x}  + \mu^2\B{b}^\top\B{A}\condh\B{Ab} \right).
%\end{align*}
%A little bit of algebra shows that $\hat{w}(\B{z})/\hat{w}(\B{x}) = w(\B{z})/w(\B{x})$ and therefore, the  MH acceptance ratio simplifies to the desired result. $\square$

\textit{Proof of Proposition~\ref{p_mh_ratio2}:} The difference between the true and the approximate covariance matrices can be expressed as
\begin{align*}
 \cond^{-1} -\condh^{-1} = & \> \mu \B{A}^\top\B{A} + \sigma \B{L}^\top\B{L} -  \left( \mu \B{L}^\top\B{V}_k \B{\Lambda}_k \B{V}_k^\top\B{L} + \sigma \B{L}^\top\B{L}\right) ,\\
= & \> \mu\B{L}^\top\left( \B{L}^{-\top}\B{A}^\top\B{A} \B{L}^{-1} -  \B{V}_k \B{\Lambda}_k \B{V}_k^\top\right) \B{L},\\
= & \> \mu \B{L}^\top\left(\sum_{j=k+1}^n\lambda_j\B{v}_j\B{v}_j^\top\right)\B{L},
\end{align*}
giving that $\log w(\B{x})  =  -\frac{\mu}{2}\sum_{j=k+1}^n \lambda_j \left(\B{v}_j^\top\B{Lx}\right)^2 \leq 0,$ and hence the acceptance ratio is given by~\eqref{e_mh_ratio2}. $\square$

\begin{lemma}
\label{l_gauss}
Suppose $\B{M}$ is symmetric positive definite. Then,
\[ \int_{\mathbb{R}^n} \exp\left(-\frac{1}{2}\B{z}^T\B{Mz} + \B{J}^\top\B{z}\right) d\B{z} \> = \> \frac{(2\pi)^{n/2}}{\det(\B{M})^{1/2}} \exp\left( \frac{1}{2}\B{J}^\top \B{M}^{-1}\B{J}\right).  \]
\end{lemma}
\begin{proof} See~\cite[Lemma B.1.1]{SantnerEtAl03}.
%Completing the squares, we have
%\[-\frac{1}{2}\B{z}^T\B{Mz} + \B{J}^\top\B{z} = -\frac{1}{2}(\B{z}-\B{M}^{-1}\B{J})^\top\B{M} (\B{z}-\B{M}^{-1}\B{J}) + \frac{1}{2}\B{J}^\top \B{M}^{-1}\B{J}. \]
%Plugging in this expression and making a change of variables $\B{y} \leftarrow \B{z}- \B{M}^{-1}\B{J}$, the integral reduces to
%\[ \exp\left( \frac{1}{2}\B{J}^\top \B{M}^{-1}\B{J} \right)\int_{\mathbb{R}^n} \exp\left(-\frac{1}{2}\B{y}^\top\B{M}\B{y}\right) d\B{y}.\]
%The integral can be computed using elementary calculus and the answer is given in the statement of the lemma.
\end{proof}

\begin{lemma}\label{l_moments}
The moments of the acceptance ratio are
\[ \mathbb{E}_{\B{z}|\B{x}} [\eta^m(\B{z},\B{x})] \>= \> \frac{1}{N_mw^m(\B{x})}\,, \]
where $N_m$ is defined in~\eqref{e_constant}.
\end{lemma}
\begin{proof}
The proof proceeds in four steps.

\textit{1. Simplifying $\mathbb{E}_{\B{z}|\B{x}} [w^m(\B{z})]$.} We focus on $\mathbb{E}_{\B{z}|\B{x}} [w^m(\B{z})]$. By definition, this is
$$\int_{\mathbb{R}^n} w^m(\B{z}) g(\B{z})d\B{z} = \frac{\exp\left(-\frac{1}{2}(\xcondh)^\top\condh^{-1}\xcondh\right)}{(2\pi)^{n/2}\det(\condh)^{1/2}} \int \exp\left(-\frac{1}{2}\B{z}^T\B{Mz} + \B{J}^\top\B{z}\right) d\B{z}, $$
where, by using $\xcondh = \mu\condh\B{A}^\top\B{b}$, we get,
\[
    \B{M} = m(\cond^{-1}-\condh^{-1}) + \condh^{-1}, ~~\text{ and } ~~\B{J} = \condh^{-1}\xcondh= \mu\B{A}^\top\B{b}.
\]
Applying Lemma~\ref{l_gauss} and rearranging, we have
\begin{equation}\label{e_inter1}
 \mathbb{E}_{\B{z}|\B{x}} [w^m(\B{z})] = \frac{\exp\left(\frac{\mu^2}{2}(\B{A}^\top\B{b})^\top (\B{M}^{-1}-\condh)\B{A}^\top\B{b}\right)}{\det(\B{M})^{1/2} \det(\condh)^{1/2}}.
\end{equation}
We focus on the numerator and denominator of \eqref{e_inter1} separately.

\textit{2. Denominator of~\eqref{e_inter1}.} Note that $$\B{M} = m(\cond^{-1}-\condh^{-1})  + \condh^{-1} = m\cond^{-1} + (1-m)\condh^{-1} $$ and furthermore, \[ \B{M} = \B{L}^\top \left(  \mu \sum_{j=1}^k\lambda_j \B{v}_j\B{v}_j^\top + m\mu\sum_{j=k+1}^n\lambda_j \B{v}_j\B{v}_j^\top + \sigma \B{I}\right) \B{L}.\]
%Following the steps of Proposition~\ref{p_mh_bound}, we can show that
Using the properties of determinants, it can be shown that
\[
    \det(\B{M}) = \sigma^n \det(\B{L})^2 \prod_{j=1}^k\left(1 + \frac{\mu}{\sigma}\lambda_j\right)\prod_{j=k+1}^n\left(1 + \frac{m\mu}{\sigma}\lambda_j\right).
\]
Similarly,
\[
    \det(\condh^{-1}) = \sigma^n \det(\B{L})^2 \prod_{j=1}^k\left(1 + \frac{\mu}{\sigma}\lambda_j\right).
\]
Combining these results, the denominator of~\eqref{e_inter1} becomes
 $$ \det(\B{M})^{1/2}\det(\condh)^{1/2} = \sqrt{\frac{\det(\B{M})}{\det(\condh)}} =  \prod_{j>k}\left(1 + \frac{m\mu}{\sigma}\lambda_j\right)^{1/2}.$$

\textit{3. Numerator of~\eqref{e_inter1}.} Consider $\B{M}^{-1} - \condh$. By the Woodbury matrix identity,
\begin{align*}
\B{M}^{-1} - \condh = & \left(m \cond^{-1} + (1-m)\condh^{-1}\right)^{-1} - \condh \\
= &\frac{1}{\sigma} \B{L}^{-1}\left(\B{I} - \sum_{j=1}^k\frac{\mu\lambda_j}{\mu\lambda_j + \sigma}\B{v}_j\B{v}_j^\top - \sum_{j=k+1}^n \frac{m\mu\lambda_j}{m\mu\lambda_j+\sigma} \B{v}_j\B{v}_j^\top \right)\B{L}^{-\top} \\
& \qquad -\frac{1}{\sigma}\left( \B{I} - \sum_{j=1}^k\frac{\mu\lambda_j}{\mu\lambda_j + \sigma}\B{v}_j\B{v}_j^\top\right)\B{L}^{-\top}\\
 =& -\frac{1}{\sigma}\sum_{j=k+1}^n \frac{m\mu\lambda_j}{m\mu\lambda_j+\sigma} \B{L}^{-1}\B{v}_j\B{v}_j^\top\B{L}^{-\top}.
\end{align*}
The numerator is therefore,
\[
    \exp\left(\frac{\mu^2}{2}(\B{A}^\top\B{b})^\top (\B{M}^{-1}-\condh)\B{A}^\top\B{b}\right)= \exp\left(-\frac{\mu^2}{2\sigma}\sum_{j=k+1}^n \frac{m\mu\lambda_j}{m\mu\lambda_j+\sigma} (\B{b}^\top\B{A}\B{L}^{-1}\B{v}_j)^2\right).
\]
\textit{4. Combining intermediate results.} Plugging the results of Steps 2 and 3 into~\eqref{e_inter1} gives us $\mathbb{E}_{\B{z}|\B{x}} [w^m(\B{z})] = \frac{1}{N_m}$, where $N_m$ is defined in~\eqref{e_constant}. The proof readily follows because $\mathbb{E}_{\B{z}|\B{x}} [\eta^m(\B{z},\B{x})] = \mathbb{E}_{\B{z}|\B{x}} [w^m(\B{z})]/w^m(\B{x})$.
\end{proof}

\textit{Proof of Theorem~\ref{p_expect}:} From Lemma~\ref{l_moments}, we have $\mathbb{E}_{\B{z}|\B{x}} [\eta^m(\B{z},\B{x})] \>= \> \frac{1}{N_mw^m(\B{x})}$. The first result follows immediately by plugging in $m=1$. % For the second result, applying Chebyshev's inequality
%\[ \text{Prob} \left\{ \left|  \eta(\B{z},\B{x}) - \bbE_{\B{z}| \B{x}} [\eta(\B{z}, \B{x})]  \right| \geq \epsilon  \right\} \leq \frac{\text{Var}_{\B{z}|\B{x}} [\eta(\B{z},\B{x})]}{\epsilon^2 }.\]
For the second result, we use the fact that for a random variable $X$ with $\mathbb{E}(X^2) < \infty$, $\mathbb{V}[X] = \mathbb{E}[(X-\mathbb{E}[X])^2] = \mathbb{E}[X^2] - (\mathbb{E}[X])^2$. The result follows from~\eqref{e_expect} and applying Lemma~\ref{l_moments} with $m=2$. $\square$

\textit{Proof of Proposition~\ref{p_mh_bound}:} Using $h$ as defined in \eqref{EQ:Target-Density} and the definition of the proposal density $g$,
\[\frac{h(\B{x})}{g(\B{x})} =  \sqrt{\frac{\det( \condh)}{\det(\cond)}} \ w(\B{x}) \exp \left( -\frac{\mu^2}{2} \B{b^\top A} (\cond - \condh) \B{A^\top b} \right) .\]

%Note that $\sqrt{ {\det(\condh)}/{\det(\cond)} } = \sqrt{ {\det( \cond^{-1})}/{\det(\condh^{-1})} }$. Next,
%\[ \det( \cond ^{-1}) = \det(\mu \B{A}^T\B{A} + \sigma \B{L}^T\B{L}) = \sigma^n \det(\B{L})^2  \det( \B{I} + \frac{\mu}{\sigma} \B{V\Lambda V}^T )\,. \]

%Since $\B{VV}^T = \B{I}$, we can simplify $\det( \B{I} + \frac{\mu}{\sigma} \B{V\Lambda V}^T ) = \prod_{j=1}^n (1+\frac{\mu}{\sigma}\lambda_j)$; similarly we can show
%$\det( \B{I} + \frac{\mu}{\sigma} \B{V}_k\B{\Lambda}_k \B{V}^T_k ) =  \prod_{j=1}^k (1+\frac{\mu}{\sigma}\lambda_j)$ and therefore,

%$$\sqrt{{\det(\condh) }/{{\det(\cond)}}} = \left[\prod_{j>k} \left(1+\frac{\mu}{\sigma}\lambda_j\right) \right]^{1/2}.$$
From the proof of Lemma~\ref{l_moments}, $\B{M}^{-1} = \cond$ when $m=1$. Comparing terms with~\eqref{e_inter1} this  gives us  $h(\B{x}) = N_1g(\B{x})w(\B{x})$, where $N_1$ is defined in~\eqref{e_constant}. From the proof of Proposition~\ref{p_mh_ratio2}, $\log w(\B{x})\leq 0$, and therefore $ w(\B{x}) \leq 1$. The desired result follows. Note that the bound is tight because $w(\B{0}) = 1 $. $\square$

\textit{Proof of Theorem~\ref{p_randsvd}:} It is easy to show that $w(\B{x}) \leq 1$, therefore, $\eta(\B{z},\B{x}) \geq w(\B{z})$. From Proposition~\ref{p_mh_ratio1}, we need to consider the quadratic form %involving the difference in the precision matrices
\[ \frac{1}{2}\B{z}^\top (\cond^{-1} -\condh^{-1}) \B{z} = \frac{\mu}{2}\B{z}^\top\B{L}^\top (\B{H}-\widehat{\B{H}}) \B{Lz}, \]
where $\widehat{\B{H}} $ is the low-rank approximation, see~\eqref{eqn:randlowrank}. This follows from $\condh^{-1} = (\mu \B{L}^\top\widehat{\B{H}}\B{L} + \sigma \B{L}^\top\B{L})$. Using the Cauchy-Schwartz inequality, we can bound (in the spectral norm)
\[ \B{z}^\top\B{L}^\top (\B{H}-\widehat{\B{H}}) \B{Lz} \leq \|\B{Lz}\|^2_2 \| \B{H}-\widehat{\B{H}}\|_2.\]
Arguing as in~\cite[Section 5.3]{halko2011finding}, we have $\| \B{H}-\widehat{\B{H}}\|_2 \leq 2 \| \B{H} - \B{QQ}^\top \B{H}\|_2$.
Applying~\cite[Theorem 10.6]{halko2011finding} we have
\begin{equation}\label{e_expect_rand} \mathbb{E}_{\B{\Omega}}\, \| \B{H}-\widehat{\B{H}}\|_2 \> \leq \> 2 \left[\alpha \lambda_{k+1} + \beta\left(\sum_{j=k+1}^{n}\lambda_j^2\right)^{1/2}\right],\end{equation}
with constants $\alpha$ and $\beta$ given in the statement of the result. Applying Jensen's inequality
\[ \mathbb{E}_{\B{\Omega}}\,  [\eta(\B{z},\B{x})] \geq \exp\left(-\mu\|\B{Lz}\|^2 \mathbb{E}_{\B{\Omega}}\,\| \B{H}-\widehat{\B{H}}\|_2 \right).\]
Plug in~\eqref{e_expect_rand} into the above equation to complete the proof. $\square$

\iffalse
\textit{Proof of Proposition~\ref{p_gkbd}}
As in the proof of  Proposition~\ref{p_randsvd},
\[ w(\B{x}) \geq \exp (-\frac{\mu}{2}\B{x}^\top\B{L}^\top (\B{H}-\widehat{\B{H}}) \B{Lx}, \]
where $\widehat{\B{H}} = \B{K}_k\B{B}_k^\top\B{K}_k \B{K}_k^\top$ is the low-rank approximation.
From~\eqref{e-uk} and~\eqref{e_vk}, we have
\[ \B{HK}_k = \B{K}_k\B{B}_k^\top\B{K}_k + \alpha_{k+1}\beta_{k+1}\B{k}_k\B{e}_k.\]
As before,  $ \| \B{H}-\widehat{\B{H}}\| \leq 2\|\B{H}-\B{H}\\|$

$\square$
\fi

\section{Implementing the Proper Jeffreys Prior}\label{app:pj_prior}
Here we briefly discuss an implementation of the LRIS algorithm when the Scott-Berger prior~\eqref{EQN:Scott-Berger} is used instead of the independent conjugate Gamma priors.

In this case, the model~\eqref{EQN:Scott-Berger} becomes
\begin{eqnarray}\label{EQN:BayesModel-SB}
    \begin{aligned}
        \B{y} \mid \B{x}, \kappa^2 &\sim N(\B{Ax}, \kappa^2\B{I})\\
        \B{x} \mid \tau^2 &\sim N(\B{0}, \tau^2\B{\Gamma})\\
        \pi(\tau^2 \mid \kappa^2) &= \kappa^{-2}(1 + \tau^2/\kappa^2)^{-2}\\
        \pi(\kappa^2) &= \kappa^{-2}.
    \end{aligned}
\end{eqnarray}

To obtain the full conditional densities necessary for Gibbs sampling, it is convenient to reparameterize the model with $\upsilon = \tau^2/\kappa^2$.  %\sarah{Discuss whether these should be moved to the Appendix or Supplementary Materials.}
After the change of variables, the joint posterior \eqref{EQ:Joint-Conditional-Explicit} becomes
\begin{eqnarray}\label{EQN:JointPosterior-SB}
\begin{aligned}
    \pi(\B{x}, \kappa^2, \upsilon \mid \B{b}) &\propto (\kappa^2)^{-\frac{m+n}{2}-1}  \upsilon^{-n/2}(1+\upsilon)^{-2} \\
        & \quad \times \exp\left(-\frac{1}{2\kappa^2}\left[ \| \B{Ax} - \B{b} \|^2_{2} + \frac{1}{\upsilon}\| \B{Lx} \|^2_{2} \right]\right).
\end{aligned}
\end{eqnarray}
The full conditional on $\B{x}$ is a standard result for the Normal-Normal model \eqref{EQ:XFC}; i.e.
\[
	\B{x} \mid \kappa^2, \upsilon, \B{b} \sim \mc{N}(\xcond, \cond),
\]
where $\cond = \kappa^2 \left( \B{A}^\top \B{A} + \upsilon^{-1} \prior^{-1} \right)^{-1}$ and $\xcond = \left( \B{A}^\top \B{A} + \upsilon^{-1} \prior^{-1} \right)^{-1} \B{A}^\top \B{b}$. We can sample from this density using our proposed low-rank approximation approach. (See Section~\ref{SEC:Proposal}.) The full conditional on $\kappa^2$ is
\[
	\kappa^2 \mid \B{x}, \upsilon, \B{y} \sim \text{InvGamma} \left(\frac{m+n}{2}, \frac{1}{2}\| \B{Ax} - \B{b}\|^2_{2} + \frac{1}{2\upsilon}\| \B{Lx}\|^2_{2} \right).
\]
To draw from this density, draw $W \sim \text{Gamma}\left(\frac{m+n}{2}, \frac{1}{2}\| \B{Ax} - \B{b}\|^2_{2} + \frac{1}{2\upsilon}\| \B{Lx}\|^2_{2} \right)$ and set $\kappa^2 = 1/W$. The full conditional on $\upsilon$ is
\begin{eqnarray}\label{EQN:UpsilonDensity}
    \pi(\upsilon \mid \B{x}, \kappa^2, \B{b}) &\propto \upsilon^{-(n/2 + 1)-1}\exp\left[-\frac{1}{2 \kappa^2 \upsilon} \| \B{Lx}\|^2_2 \right]\left(\frac{\upsilon}{1+\upsilon}\right)^2\\
        &\equiv h(\upsilon)\left(\frac{\upsilon}{1+\upsilon}\right)^2,
\end{eqnarray}
where $h(\upsilon)$ is the density of an $\text{InvGamma}((n+2)/2, \| \B{Lx}\|^2_{2} /(2\kappa^2))$ distribution. Thus we can use an independence Metropolis-Hastings algorithm with candidate distribution $$\text{InvGamma}((n+2)/2, \| \B{Lx}\|^2_2 /(2\kappa^2)).$$
%This density is not available directly. However, since
%\[
%    \frac{\pi(\upsilon \mid \B{x}, \kappa^2, \B{b})}{h(\upsilon)} = \left(\frac{\upsilon}{1+\upsilon}\right)^2 \leq 1 ~\forall \upsilon,
%\]
%In this case, the M-H ratio is
%\begin{eqnarray*}
%    \alpha_{\upsilon} &:=& \min\left(1, \frac{\pi(\upsilon^{\ast} \mid \B{x}, \kappa^2, \B{b})h(\upsilon^{(t-1)})}{\pi(\upsilon^{(t-1)} \mid \B{x}, \kappa^2, \B{b})h(\upsilon^{\ast})}\right)\\
%        &=& \min\left(1, \left[\frac{\upsilon^{\ast}(1 + \upsilon^{(t-1)})}{\upsilon^{(t-1)}(1 + \upsilon^{\ast})}\right]^2\right),
%\end{eqnarray*}
%where $\upsilon^{\ast}$ is the proposed value and $\upsilon^{(t-1)}$ is the current state at iteration $t$.

\bibliography{mcmc_paper-new} % Using Arvind's updated bib file
\bibliographystyle{plain}

%%%% Merge with supplementary materials
\clearpage
%\documentclass[10pt,review]{siamonline0516}
%
%
%\usepackage{amsmath,amssymb,graphicx, mathtools}
%
%
%\usepackage{graphics,graphicx,color}         	% Allows inclusion of eps files.
%\usepackage{epstopdf}  % Added by Andrew to make it easy to compile directly to PDF with eps figures
%%\graphicspath{{/EPSF/}{../Figures/}{Figures/}}
%\usepackage{url}		% Allows good typesetting of web URLs.
%\usepackage{hyperref}
%\usepackage[normalem]{ulem}
%
%\newcommand{\arvind}[1]{{\color{red} #1}}
%\newcommand{\sarah}[1]{{\color{purple} #1}}
%\newcommand{\andrew}[1]{{\color{blue} #1}}
%
%% Added by Sarah
%\usepackage{floatrow}
%% Table float box with bottom caption, box width fixed to half page
%\newfloatcommand{capbtabbox}{table}[][0.5\textwidth]
%
%\usepackage{xr} % cross references
%\externaldocument[I-]{3proposal_R2}
%\externaldocument{mcmc_paper_R2}
%%\externaldocument{5_deblur.tex}
%
%\usepackage{color}  % Added by Andrew
%\frenchspacing  % Added by Andrew
%
%\usepackage{fancyhdr}
%
%
%\input{setupmacros.sty} % Added by Sarah
%\usepackage[boxed]{algorithm2e} % for algorithms
%
%\newtheorem{prop}{Proposition}

%\newcommand{\xcondh}{\hat{\boldsymbol{x}}_{\text{cond}}}
%\newcommand{\condh}{\widehat{\boldsymbol{\Gamma}}_{\text{cond}}}
%\newcommand{\xcond}{\boldsymbol{x}_{\text{cond}}}
%\newcommand{\cond}{\boldsymbol{\Gamma}_{\text{cond}}}
\newcommand{\zcond}{\boldsymbol{z}_{\text{cond}}}
%\newcommand{\prior}{\boldsymbol{\Gamma}_{\text{pr}}}
%\newcommand{\normtwo}[1]{\| #1\|_2}

%%%%%%%%%%
% TITLE AND AUTHOR
%
%\title{Supplementary Material for ``Low Rank Independence Samplers in Bayesian Inverse Problems"}
%
%\author{
%	D. Andrew Brown\thanks{
%		Department of Mathematical Sciences,
%		Clemson University, Clemson, SC 29634;
%		\href{mailto:ab7@clemson.edu}{ab7@clemson.edu};
%	}
%	\and
%	Arvind Saibaba\thanks{
%		Department of Mathematics,
%		North Carolina State University, Raleigh, NC 27695;
%		\href{mailto:asaibab@ncsu.edu}{asaibab@ncsu.edu};
%	}
%	\and
%	Sarah Vall\'elian\thanks{
%		Statistical and Applied Mathematical Sciences Institute,
%		Research Triangle Park, NC 27709;
%		\href{mailto:svallelian@samsi.info}{svallelian@samsi.info};
%	}
%}

%%%%%%%%%%
% BEGIN DOCUMENT
%
%\begin{document}
%
%\maketitle

\begin{center}
{\LARGE Supplementary Material for {Low Rank Independence Samplers in Bayesian Inverse Problems}}
\vspace{1em}

{\large D. Andrew Brown \qquad Arvind Saibaba \qquad Sarah Vall\'elian}
\end{center}
\setcounter{section}{0}
\setcounter{equation}{0}
\setcounter{figure}{0}
\setcounter{table}{0}
\setcounter{page}{1}
\renewcommand{\theequation}{S\arabic{equation}}
\renewcommand{\thetable}{S\arabic{table}}
\renewcommand\thefigure{S\arabic{figure}}
\renewcommand{\theHtable}{Supplement.\thetable}
\renewcommand{\theHfigure}{Supplement.\thefigure}

% Things which can go in this document:
%	Non-centered parameterization
%	Other numerical examples, e.g. running 2D deblurring or CT again with different priors/different hyperparameters, the 1D deblurring stuff Andrew did on all types of priors
%
% Things which probably belong in the Appendix instead:
%	Pseudocode for the sampling algorithms using centered Gamma-Gamma and proper Jeffreys priors since we use those in the paper
\vspace{2em}
This Supplementary Material contains additional material referenced in the manuscript that could aid the reader, but is not essential. We empirically demonstrate mixing behavior and efficiency of an LRIS-based Metropolis-Hastings-within-Gibbs sampler as a function of rank of the approximation, elaborate on possibilities for prior modeling of the precision (variance) components in the hierarchical model for the Bayesian inverse problem, and discuss the use of noncentered parameterization in our proposed LRIS algorithm, supported by application to the 2D deblurring example. We further demonstrate computational feasibility afforded by our proposed approach by considering MCMC/LRIS-based reconstruction of a distribution of relaxation times in nuclear magnetic resonance (NMR) relaxometry. Supplementary Figures are at the end.
%\section{Additional Proofs of Results}
%{\bf Put non-essential / less interesting proofs here to save space in the main manuscript.}

\section{Effect of the Proposal Rank on Chain Mixing}
\begin{figure}[tb]
    \centering
    \includegraphics[scale = 0.5, clip= TRUE]{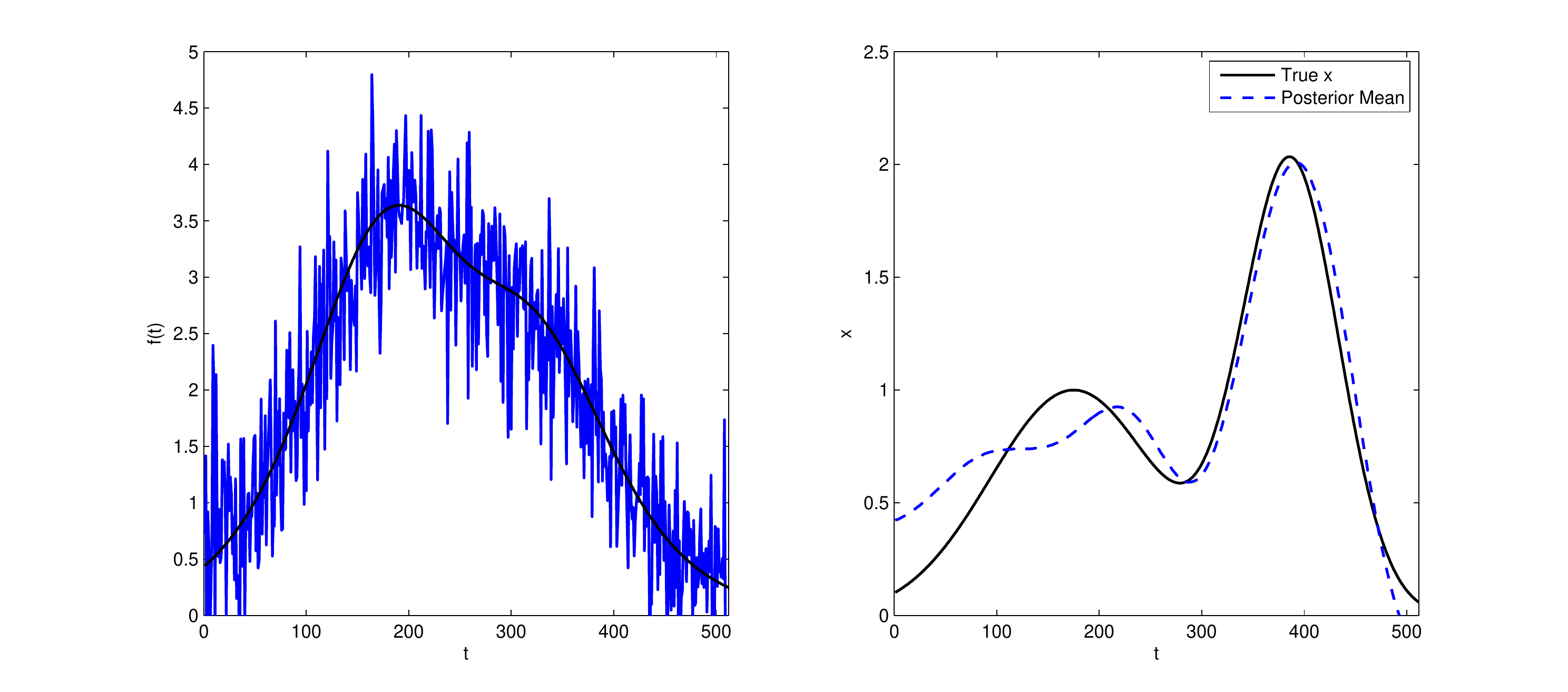}
    \caption{Observed data (left panel) and true solution (right panel) for the one-dimensional image restoration example. The solid black line in the left panel represents the noise-free observations $\B{Ax}$. The dashed blue line is the true posterior mean (approximated from a very long block Gibbs Markov chain).}
    \label{FIG:1DTruePostMean}
\end{figure}

To study how the quality of the low-rank approximation affects the mixing of an MCMC algorithm, we consider the image reconstruction problem originally appearing in \cite{shaw1972}. This is a one-dimensional image restoration problem in which the blurred image can be expressed as a Fredholm integral equation of the first kind, $g(s) = \int_{-\pi/2}^{\pi/2} K(s,t)f(t) \, \dx t,$ where $K(s,t) = (\cos(s) + \cos(t))^2\left(\sin(u)/u\right)^2$ with $u = \pi (\sin(s) + \sin(t))$, and $f$ is the true one-dimensional image given by $f(t) = 2 \exp(-6(t - 0.8)^2) + \exp(-2(t+0.5)^2)$. The problem is to reconstruct $f$ given $g$ and $K$. Using the implementation in the \texttt{Matlab} package \texttt{Regularization Tools}
 \cite{Hansen2007}, we discretize the integral via quadrature over $n=512$ points and corrupt the observations with one percent noise. The resulting model is $\B{b}= \B{Ax} + \B{\epsilon}, ~\B{\epsilon} \sim \cN(\B{0}, \B{\Gamma}),$ where $\B{b}$ is the observed data, $\B{A} \in \mathbb{R}^{n\times n}$ is the discretized forward model, $\B{x}$ is the discretized solution, and $\B{\Gamma} = 0.01^2\|\B{b}\|^2\B{I}$. The observed data $\B{b}$ and solution $f$ are displayed in Figure \ref{FIG:1DTruePostMean}.

In the hierarchical Bayesian model, we use a zero mean Gaussian process (GP) prior \cite{SantnerEtAl03, RasWilliams06}, $f(\cdot) \sim \mathcal{GP}(0, \sigma^{-1}R(\cdot, \cdot))$. We take the correlation function to be in the power exponential family, $R(t_i, t_j) = \exp\left(-|t_i - t_j|/l\right)$, with correlation length parameter $l = \pi/2$. We use vague Gamma priors about the noise precision $\mu$ and prior precision $\sigma$ with $a_{\mu} = b_{\mu} = a_{\sigma} = b_{\sigma} = 0.1$. We remark that for this example, the forward model is fast-running so that our computationally-cheap LRIS is actually not necessary. However, the fast-running model makes it feasible to run a large number of replicate MCMC chains in a reasonable amount of time, thus allowing us to empirically assess the mixing behavior of our proposed approach as a function of rank of the proposal distribution.

% AB: I have no idea why it's changing the color in this paragraph the way that it is. I didn't look at it very long, though, and moved on.
Similar to the analysis in \cite{JohnsonEtAl13}, we consider two empirical measures to assess the mixing behavior of the Markov chains as a function of rank of the approximation in the LRIS algorithm. To assess the statistical efficiency of using the mean of the MCMC output to estimate $\mathbb{E}(\B{x} \mid \B{b})$, we use the mean squared error $MSE(\overline{\B{x}}_n) := \mathbb{E}[(\overline{\B{x}}_n - \mathbb{E}(\B{x} \mid \B{b}))^2]$, where $\overline{\B{x}}_n$ is the sample mean of a chain of length $n$, $\B{x}_{(1)}, \ldots, \B{x}_{(n)}$, obtained from an MCMC run. To approximate the MSE, we find the sample mean, $\overline{\B{x}}^\ast$, from a very long run of an ordinary block Gibbs sampler and treat this as the true posterior mean $\mathbb{E}(\B{x} \mid \B{b})$. We run an additional $m$ independent Markov chains using LRIS-based Metropolis-Hastings-within-Gibbs, each of length $n$, whence we can approximate $MSE(\overline{\B{x}}_n)$ with $\widehat{MSE}(\overline{\B{x}}_n) = m^{-1}\sum_{i=1}^m(\overline{\B{x}}_n^{(i)} - \overline{\B{x}}^\ast)^2$, where $\overline{\B{x}}_n^{(i)}$ is the sample mean obtained from the $i^{\text{th}}$ chain. The second measure we consider is the expected squared Euclidean jump distance, defined as $ESEJD = \mathbb{E}(\|\B{x}_{(t+1)} - \B{x}_{(t)}\|_2^2)$. This quantity is indicative of how well a Markov chain is exploring the marginal posterior distribution of the estimand $\B{x}$. To approximate this expected value, we again run $m$ independent chains, each of length $n$. For each chain with rank $k$, we find the mean squared Euclidean jump distance $MSEJD_k := (n-1)^{-1}\sum_{t=1}^{n-1}\|\B{x}_{(t+1)} - \B{x}_{(t)}\|_2^2$. Then we obtain an estimate of ESEJD with $\widehat{ESEJD}_k = m^{-1}\sum_{i=1}^m MSEJD_k^{(i)}$, where $MSEJD_k^{(i)}$ is mean squared Euclidean jump distance from the $i^{\text{th}}$ chain.

\begin{figure}[tb]
    \centering
    \includegraphics[scale = 0.45, clip= TRUE]{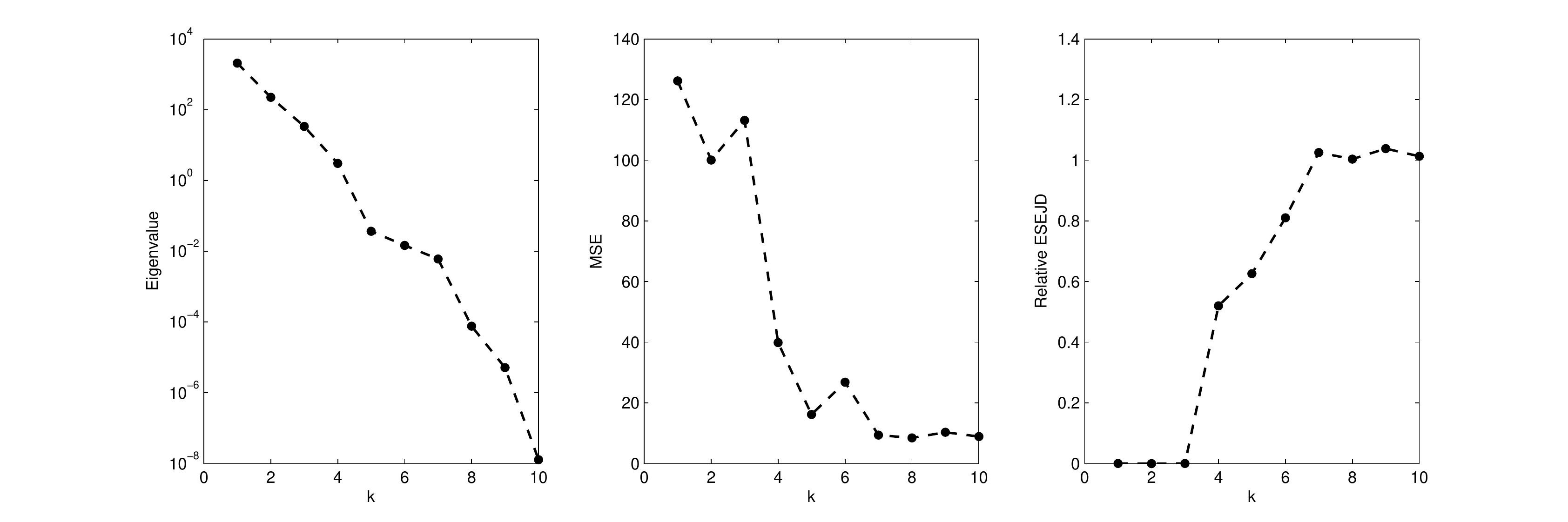}
    \caption{Left panel: Spectrum (up to the first 10 eigenvalues) of the prior-preconditioned Hessian matrix for the one-dimensional image restoration example; Middle: The approximate mean squared error (MSE) in estimating $\mathbb{E}(\B{x} \mid \B{b})$ versus the rank of the proposal distribution in the LRIS; Right: Approximate expected squared Euclidean jump distance (ESEJD) of the $\B{x}$ chain relative to the ESEJD from block Gibbs sampling at convergence.}
    \label{FIG:SpectrumAndMixing}
\end{figure}

The left panel of Figure \ref{FIG:SpectrumAndMixing} displays the first ten eigenvalues of the prior-preconditioned Hessian for the one-dimensional reconstruction example. As most of the information is captured in the first ten eigenvalues, we consider the LRIS algorithm using proposal distributions of ranks $k= 1, \ldots, 10$. For each proposal distribution, we run $m = 100$ independent Markov chains of length $n = 2000$, discarding the first $1,000$ draws as a burn-in period. To approximate the true posterior distribution, we run an ordinary block Gibbs sampler for $10^6$ iterations and approximate $\mathbb{E}(\B{x} \mid \B{b})$ with the mean of the last $5\times 10^5$ draws. This target is displayed in the right panel of Figure \ref{FIG:1DTruePostMean}. %({\bf ARVIND AND SARAH: We can put basic convergence diagnostics for the Gibbs sampler in the Supplementary Material if we want to. It may or may not be necessary given, that I ran 1 million iterations!})

\begin{figure}[tb]
    \centering
    \includegraphics[scale = 0.45, clip= TRUE]{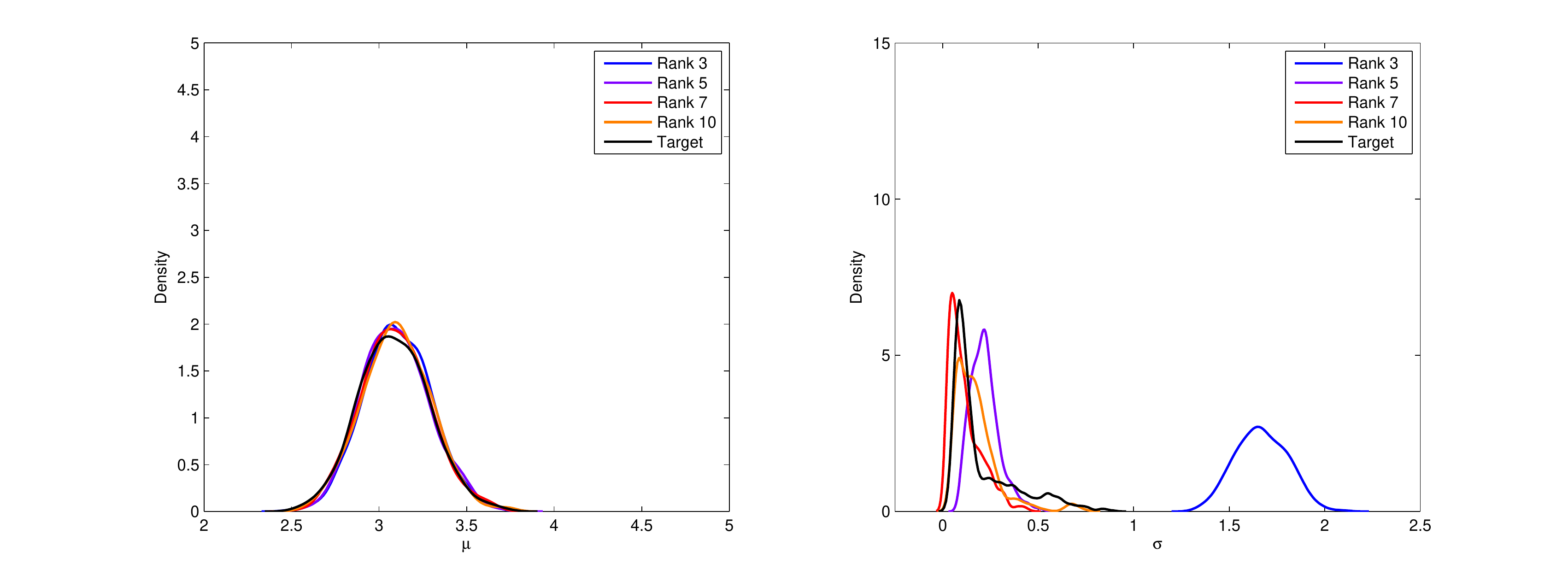}
    \caption{Smoothed approximate marginal posterior densities of the noise precision $\mu$ (left) and prior precision $\sigma$ (right) in the one-dimensional image restoration example. The target density is estimated with the last 1,000 draws from an ordinary block Gibbs sampler with chain length $10^6$.}
    \label{FIG:MarginalsAndMixing}
\end{figure}

The middle panel of Figure \ref{FIG:SpectrumAndMixing} displays the mean squared errors in the posterior mean estimates versus the rank of the LRIS proposal distribution. We see that the statistical efficiency of the sample mean increases sharply as non-negligible eigenvalues are added to the low-rank approximation, and becomes steady at $k=7$. The right panel of the Figure displays the approximate expected squared Euclidean jumping distance of the Markov chains relative to the average squared jumping distance of the block Gibbs sampler, $\widehat{ESEJD}_k/ESEJD_{\text{Gibbs}}, ~k= 1, \ldots, 10$. Similar to MSE, we can glean that only seven eigenvalues are necessary to obtain within 2,000 iterations a sample whose behavior is equivalent to a converged block Gibbs sampler. This equivalence is further supported in Figure \ref{FIG:MarginalsAndMixing}, which displays smoothed estimates of the marginal densities of $\mu$ and $\sigma$ as the rank increases. We again see the ability of the LRIS-based algorithm to closely estimate the true marginal density with only 2,000 MCMC iterates. It is interesting to observe that the noise precision $\mu$ is well identified regardless of the rank of the approximation, whereas $\sigma$ is much more difficult to estimate. This reflects identifiability issues that are characteristic of ill-posed Bayesian inverse problems \cite{bardsley2012mcmc}. %({\bf ARVIND AND SARAH: If the paper is getting too figure-heavy, we can show the MSE and ESEJD in a Table instead of a figure, and/or we can move a couple of these figures to the supplementary material.}).

%\subsection{Summary}
These results demonstrate sharp improvement in the behavior of an LRIS-based MCMC algorithm that is possible when the spectrum of the prior-preconditioned Hessian decays rapidly. Thus, for computationally expensive forward models, there is the potential to dramatically reduce the computational burden without sacrificing the convergence behavior of the Markov chain.

\section{Priors on the Nuisance Parameters}\label{sec:Priors}
Consider the Bayesian linear inverse problem with forward operator $\B{A} \in \bbR^{m\times n}$, assuming independent observations with common precision $\mu$. We assume a Gaussian prior on the target $\B{x}$ to correspond to the $L_2$ penalty in regularized inversion, with prior covariance matrix $\prior$ known up to a multiplicative constant (precision) $\sigma$. Let $\B{b}$ denote the observed data. Then the model is
\begin{eqnarray}\label{EQ:BIPMod}
    \begin{aligned}
        \B{b} \mid \B{x}, \mu &\sim \mc{N}(\B{Ax}, \mu^{-1}\B{I})\\
        \B{x} \mid \sigma &\sim \mc{N}(\B{0}, \sigma^{-1}\prior)\\
        (\mu, \sigma) &\sim \Pi,\\
    \end{aligned}
\end{eqnarray}
where $\Pi$ is some distribution with support $\mathbb{R}^+ \times \mathbb{R}^+$. In the following, we consider two different specifications of $\Pi$.

\subsection{Conditionally Conjugate Gamma Priors}\label{SEC:Gammas} The most straightforward case (and by far the most common) is to let $\mu \sim \text{Gamma}(a_{\mu}, b_{\mu})$ and $\sigma \sim \text{Gamma}(a_{\sigma}, b_{\sigma})$ independently, where we use the shape/rate parameterization of a Gamma distribution; e.g.,
\[
    \text{Gamma}( \mu \mid a, b) \propto \mu^{a-1}\exp(-b\mu), ~~\mu > 0.
\]
With these Gamma prior models on $\mu$ and $\sigma$, the joint posterior distribution of the model specified by \eqref{EQ:BIPMod} is given by
\begin{equation}\label{EQ:JointPost-Gammas}
	\pi(\B{x},\mu, \sigma \mid \B{b}) \propto \mu^{m/2 + a_{\mu}-1} \sigma^{n/2 + a_{\sigma}-1} \exp \left( -\frac{\mu}{2} \| \B{Ax} - \B{b} \|^2_{2} - \frac{\sigma}{2} \| \B{Lx} \|^2_{2} -b_{\mu} \mu - b_{\sigma} \sigma \right).
\end{equation}
Since the elements of $\B{x}$ are expected to be highly correlated in the posterior, it is desirable to update $\B{x}$ all at once in a block Gibbs sampler. The full conditional distributions in this case are
\begin{equation}\label{EQ:Conditionals-Gammas}
\begin{aligned}
	\B{x} \mid \B{b}, \mu, \sigma \>\sim & \> \mc{N}(\xcond, \cond) \\
	\mu \mid \B{b}, \B{x}, \sigma \>\sim & \> \text{Gamma}\left(m/2+a_\mu,\frac{1}{2}\normtwo{\B{Ax} - \B{b}}^2 + b_{\mu}\right) \\
	\sigma \mid \B{b}, \B{x}, \mu \> \sim & \> \text{Gamma}\left(n/2+a_\sigma,\frac{1}{2}\normtwo{\B{Lx}}^2 + b_\sigma\right),
\end{aligned}
\end{equation}
where $\cond = (\mu \B{A}^\top \B{A} + \sigma \prior^{-1})^{-1}$ and $\xcond = \mu \cond \B{A}^\top \B{b}$. The remaining question becomes specification of the hyperparameters in the Gamma priors. %The implementation of our proposed method for this choice of prior is given in Algorithm \ref{ALG:MHwG}.

%\LinesNumbered
%\begin{algorithm}[h]
%\SetAlgoLined
%\DontPrintSemicolon
%\SetKwInput{Input}{Input}
%\SetKwInput{Output}{Output}
%	\Input{Full conditional distributions \eqref{EQ:Conditionals-Gammas}, proposal density $g(\B{x})$, sample size $N$, burn-in period $N_b $.}
%	\Output{Approximate sample from the posterior distribution \eqref{EQ:JointPost-Gammas} $\{ \B{x}_{(t)}, \mu_{(t)}, \sigma_{(t)} \}_{t=N_b +1}^N$.}
%	\BlankLine
%	Initialize $\B{x}_{(0)}$, $\mu_{(0)}$, and $\sigma_{(0)}$. \;
%  	\For{$t=1$ to $N$} {
%    		Draw $\B{x}_{\ast} $ from proposal density $g(\cdot \mid \B{b}, \mu_{(t-1)}, \sigma_{(t-1)})$. \;
%		 \lIf{$u \leq \eta(\B{x}_{(t-1)}, \B{x}_{\ast})$}{ Set $\B{x}_{(t)} = \B{x}_{\ast}$.}
%		 \textbf{else} Set $\B{x}_{(t)} = \B{x}_{(t-1)}$. \;
%    		Draw $\mu_{(t)} \>\sim\>\text{Gamma}\left(m/2+a_\mu,\frac{1}{2}\normtwo{\B{Ax}_{(t)} - \B{b}}^2 + b_{\mu}\right)$. \;
%		Draw $\sigma_{(t)} \>\sim\> \text{Gamma}\left(n/2+a_\sigma,\frac{1}{2}\normtwo{\B{Lx}_{(t)}}^2 + b_\sigma\right)$. \;
%        }
%    Discard $\{ \B{x}_{(t)}, \mu_{(t)}, \sigma_{(t)} \}_{t=1}^{N_b}$ as burn-in iterates. \;
%\caption{Metropolis-Hastings-within-Gibbs (MHwG) algorithm for sampling from the posterior distribution obtained with conditionally conjugate Gamma priors.}
%\label{ALG:MHwG}
%\end{algorithm}

To impose strong prior assumptions and to stabilize the MCMC algorithm, we can rescale the observed data with $\T{\B{b}} = \B{b}/s_b$, where $s_b = \sqrt{(m-1)^{-1}\sum_{i=1}^m(b_i - \overline{b})^2}$ is the sample standard deviation.  (See, for instance, \cite{Higdon2008}.) In this case, $\mathbb{V}(\T{\B{b}})$ is expected to be reasonably close to one, though not exactly equal since correlation in the data induced by dependence on $\B{x}$ will cause $s_b$ to over- or under-estimate the true standard deviation. After rescaling, an equivalent model to \eqref{EQ:BIPMod} is
\begin{equation}\label{EQ:stdMod}
    \T{\B{b}} = \B{A}\T{\B{x}} + \T{\B{\epsilon}},
\end{equation}
where $\T{\B{x}} = \B{x}/s_b$ and $\T{\B{\epsilon}} \sim N(\B{0}, \T{\mu}^{-1}\B{I})$. This is similar to the notion of a standardized regression model. (See, e.g., \cite[Section~7.5]{Kutner2005}.) In this case, we set the hyperparameters in the prior on $\T{\mu}$ (or simply $\mu$, without loss of generality) to mildly concentrate the density about one; e.g. $a_{\mu} = b_{\mu} = 1$. By concentrating $\mu$ about 1 and allowing $\sigma$ to be vague with, say, $a_{\sigma} = b_{\sigma} = 0.1$, we do not strongly restrict values of $\lambda = \sigma/\mu$, the corresponding regularization parameter in the MAP estimator. We remark, however, that in a fully Bayesian model the primary goal is to obtain an estimate of $\B{x}$, so we are not really interested in $\mu$ or $\sigma$ in their own right (hence the term ``nuisance parameters").

If $\mu \sim \text{Gamma}(\epsilon, \epsilon)$, then $\pi(\mu) \propto \mu^{\epsilon-1}e^{-\epsilon \mu} \rightarrow \mu^{-1}$ as $\epsilon \rightarrow 0$. But $\pi(\mu) = \mu^{-1}$ is the Jeffreys prior for a scale parameter and thus is invariant to reparameterization \cite[Ch. 3]{Berger85}. For this reason it is common practice to set $\epsilon$ to some small value in a $\text{Gamma}(\epsilon, \epsilon)$ prior, say $\epsilon = 0.1$ or $\epsilon = 0.01$, to approximate the behavior of the objective prior without sacrificing propriety or conjugacy. These are the priors we use in the 2D image deblurring example in Section 4.1 of the manuscript.

\subsection{Weakly Informative Priors}
Despite the convenience associated with the Gamma priors, it was observed by \cite{Gelman06} that there is no limiting posterior distribution associated with taking $\epsilon \rightarrow 0$ in a $\text{Gamma}(\epsilon, \epsilon)$ prior, and that using such a hyperprior on the prior-level precision $\sigma$ can sometimes yield undesirable behavior. (Although Carlin and Louis \cite{CarlinLouis09} remarked that it may not make a difference in terms of the estimand of interest, $\B{x}$.) To rectify this, Gelman \cite{Gelman06} proposed as a default prior the {\em folded-t} distribution. This prior strikes a good compromise between a completely noninformative prior, which can lead to unreasonable estimates if the data are not informative about a parameter, and a strongly informative prior which prevents the data from `speaking for themselves' in determining plausible {\em a posteriori} values. As such, it is called a ``weakly informative" prior. Scott and Berger \cite{ScottBerger06} proposed what has since become known as a ``proper Jeffreys" prior \cite{PolScott2010} on the variance components. Defining $\kappa^2 = \mu^{-1}$ and $\tau^2 = \sigma^{-1}$, the proper Jeffreys prior takes
\[
    \pi(\kappa^2, \tau^2) = (\kappa^2 + \tau^2)^{-2}, ~~\kappa^2, \tau^2 > 0.
\]
This prior approximates the improper Jeffreys prior, $\pi(\kappa^2, \tau^2) = (\kappa^2 + \tau^2)^{-1}$ \cite{TiaoTan65}. Scott and Berger \cite{ScottBerger06} observed that the proper Jeffreys prior can be written as $\pi(\kappa^2, \tau^2) = \kappa^{-2}(1 + \tau^2/\kappa^2)^{-2}\kappa^{-2} \equiv \pi(\tau^2 \mid \kappa^2)\pi(\kappa^2)$, so that this model is equivalent to using the usual objective prior on the data-level variance while scaling the prior-level variance by $\kappa^2$, following the principle originally suggested by Jeffreys \cite{Jeffreys61}. The conditional prior on $\tau^2 \mid \kappa^2$ is also proper, an important consideration in finite mixture models, or when the data contain limited information about $\tau^2$. Lastly, it was observed by \cite{BrownEtAl17} that the proper Jeffreys prior is tail-equivalent to the prior obtained by placing a folded-$t_2$ on $\tau$, and thus is suitable as a default prior choice. For these reasons, this can be an attractive alternative to conjugate Gamma priors. In the CT example in Section 4.2 of the manuscript, we use the proper Jeffreys prior on the variance components. The sampling algorithm is not quite as simple since we no longer have conjugacy (see Appendix B of the manuscript), but we are still able to use our proposed low-rank independence sampler.

\section{Non-Centered Parameterizations}
In model \eqref{EQ:BIPMod}, the distribution of $\B{x}$ depends on $\sigma$, and the distribution of $\B{b}$ depends on $\B{x}$, but $\B{b}$ is conditionally independent of $\sigma$, given $\B{x}$. Under Gamma priors on $\mu$ and $\sigma$, this yields convenient conditionally conjugate distributions for use inside a block Gibbs sampler, as discussed in Section \ref{SEC:Gammas} of the Supplementary Material. However, this also leads to high correlation between the $\sigma$ and $\B{x}$ chains as the dimension of the problem increases, as noted by Bardsley \cite{bardsley2012mcmc} and Agapiou et al. \cite{agapiou2014analysis}.

A framework for potentially reducing the dependence between parameters in an MCMC algorithm is the so-called non-centered parameterization \cite{Papaspiliopoulos2003, Papaspiliopoulos2007}. A non-centered parameterization is one such that parameters are assigned {\em independent} prior distributions but still result in a model equivalent to the usual case, called a centered parameterization. The ``centered" and ``non-centered" terminology is a reference to the parameterizations considered by \cite{Gelfand1995} for efficient sampling on normal linear mixed models, where certain parameters were centered or non-centered about other parameters.

Papaspiliopoulos et al. \cite{Papaspiliopoulos2007} argued that the best choice of parameterization depends on how well the underlying parameters are identified by the data. In our case, if the data were significantly informative about $\B{x}$ so that strong Bayesian learning occurred in the posterior, then a centered parameterization would likely be appropriate. On the other hand, if $\B{x}$ is only weakly identified by the data alone and hence more dependent on prior information, then there tends to be stronger correlation between $\sigma$ and $\B{x}$ under the centered parameterization. In this case, a non-centered parameterization is likely the better option. Similar behavior was observed by Agapiou et al. \cite{agapiou2014analysis}, where it was shown that the performance of the non-centered parameterization breaks down as the data-level variance becomes small (i.e., the data become more reliable and thus contain more information about the solution). There exist also ``partially non-centered" parameterizations, which can be estimated from the data when the appropriate parameterization to use is not clear \cite{Papaspiliopoulos2003, Papaspiliopoulos2007}.

\subsection{Implementation} To determine a non-centered parameterization for the hierarchical Bayesian inverse problem, we define a random variable $\B{z}$ independent of $\sigma$ and express the distribution of $\B{x}$ in terms of these independent random variables. In our case, we have that $\B{x} \stackrel{\text{d}}{=} \sigma^{-1/2}\B{z}$, where $\B{z} \sim \mc{N}(\B{0}, \prior)$ independent of $\sigma$. Substituting this parameterization into \eqref{EQ:BIPMod}, the model becomes
\begin{eqnarray}\label{EQ:NCP}
    \begin{aligned}
        \B{b} \mid \B{z}, \mu, \sigma &\sim \mc{N} (\sigma^{-1/2}\B{Az}, \mu^{-1}\B{I})\\
        \B{z} &\sim \mc{N} (\B{0}, \prior)\\
        \mu &\sim \text{Gamma}(a_{\mu}, b_{\mu})\\
        \sigma &\sim \text{Gamma}(a_{\sigma}, b_{\sigma}).
    \end{aligned}
\end{eqnarray}
The joint posterior density is
\begin{equation}\label{EQ:JointPosterior-NCP}
    \pi(\B{z}, \mu, \sigma \mid \B{b}) \propto f(\B{b} \mid \B{z}, \mu, \sigma)\pi(\B{z})\pi(\mu)\pi(\sigma),
\end{equation}
where $f$ is the likelihood, and we adopt the conventional ambiguous use of $\pi$, understood to be defined by its arguments.

The conditional distributions for Gibbs sampling can be derived in a similar manner to the centered case discussed in Section \ref{SEC:Gammas} of the Supplementary Material, with the exception of $\sigma$. The non-centered parameterization loses the conditional conjugacy on this parameter, making an indirect sampling approach necessary. The conditional density of $\B{z}$ can be derived as the usual Normal-Normal model $\B{b} \sim \mc{N}(\B{A}_{\sigma}\B{z}, \mu^{-1}\B{I})$ and $\B{z} \sim \mc{N}(\B{0}, \prior)$, where $\B{A}_{\sigma} = \sigma^{-1/2}\B{A}$. Thus,
\[
    \B{z} \mid \B{b}, \mu, \sigma \sim \mc{N}( \zcond, \cond),
\]
where $\zcond = (\mu \B{A}_{\sigma}^\top \B{A}_{\sigma} + \B{L}^\top \B{L})^{-1} = ((\mu/\sigma)\B{A}^\top \B{A} + \B{L}^\top \B{L})^{-1}$ and $\zcond = \mu\cond\B{A}_{\sigma}^\top \B{b} = (\mu/\sigma^{1/2})\cond\B{A}^\top \B{b}$. One can easily show the same results for this parameterization as in the manuscript with appropriate substitutions of $\mu$ and $\sigma$. In particular, the same simplification of the MH acceptance ratio for the $\B{z}$ chain holds as in Proposition 1. The full conditional of $\mu$ is also straightforward. It is derived similarly to the centered case:
\begin{eqnarray*}
    \pi(\mu \mid \B{b}, \B{z}, \sigma) &\propto& \mu^{m/2 + a_{\mu}-1}\exp\left[-\mu\left(\frac{1}{2} \normtwo{\B{A}_{\sigma}\B{z}-\B{b}}^2 + b_{\mu}\right)\right]\\
    \Rightarrow \mu \mid \B{b}, \B{z}, \sigma &\sim& \text{Gamma}\left(m/2 + a_{\mu}, \frac{1}{2} \normtwo{ \B{A}_{\sigma}\B{z} - \B{b} }^2 + b_{\mu}\right).
\end{eqnarray*}

The most substantial difference between the centered and non-centered parameterizations is the loss of conditional conjugacy on $\sigma$. The conditional density for $\sigma$ is
\begin{equation}\label{EQ:deltaDens}
    \pi(\sigma \mid \B{b}, \B{z}, \mu) \propto \exp\left[-\frac{\mu}{2} \normtwo{\sigma^{-1/2}\B{Az} - \B{b}}^2 \right]\sigma^{a_{\sigma}-1}e^{-b_{\sigma}\sigma}.
\end{equation}
While there is no obvious simplification or standard distribution for $\sigma$, we can use a random walk Metropolis step to sample from it. Simplifying \eqref{EQ:deltaDens} slightly, we have
\[
    \pi(\sigma \mid \B{b}, \B{z}, \mu) \propto \sigma^{a_{\sigma}-1}e^{-b_{\sigma}\sigma}\exp\left[-\frac{\mu}{2\sigma}\left(\B{z}^\top\B{A}^\top \B{Az} - 2\sigma^{1/2}\B{z}^\top \B{A}^\top \B{b}\right)\right], ~~\sigma > 0.
\]
To eliminate boundary constraints on $\sigma$ and thus facilitate Gaussian proposals in the Metropolis algorithm, reparameterize the model with $\omega = \log(\sigma)$. Then the density for $\omega$ becomes
\[
    \pi_{\omega}(\omega \mid \B{b}, \B{z}, \mu) \propto e^{(a_{\sigma}-1)\omega - b_{\sigma}e^{\omega}}\exp\left[-\frac{\mu}{2}\left(e^{-\omega}\B{z}^\top \B{A}^\top \B{Az} - 2e^{-\omega/2}\B{z}^\top \B{A}^\top \B{b}\right)\right].
\]
To implement the Metropolis step inside the block Gibbs algorithm, we can explicitly separate the burn-in phase from the sampling phase. During the first phase, we seek a suitable proposal distribution for $\sigma$ (i.e., a suitable Gaussian proposal for $\omega$). This is done by adaptively controlling the variance $c$, adjusting its scale based on the acceptance rate of the $\sigma$ samples. A skeleton of this procedure is as follows:

\begin{enumerate} % Try to make this look nicer?
    \item
    Initialize $c = c_0$ and $accept = 0$.

    \item
    For $i \in $ (burn-in iterations)
    \begin{enumerate}
        \item
        Draw $\omega^{\ast} \sim N(\omega^{(i-1)}, c)$

        \item
        Accept/reject $\omega^{\ast}$ according to the Metropolis ratio. If accepted, $accept \leftarrow accept + 1$.

        \item
        If $\mod(i,100) = 0$, then
        \begin{enumerate}
            \item
            If $accept/100 < 0.35, ~c \leftarrow 0.75c$,

            \item
            Else, if $accept/100 > 0.5, ~c \leftarrow 1.75c$

            \item
            $accept \leftarrow 0$
        \end{enumerate}

        \item
        Repeat
    \end{enumerate}

    \item
    For $i \in $ (sampling iterations)
    \begin{enumerate}
        \item
        Draw $\omega^{\ast} \sim N(\omega^{(i-1)}, c)$, where $c$ is fixed at the value determined from the burn-in period.

        \item
        Accept/reject $\omega^{\ast}$ according to the Metropolis ratio.

        \item
        Repeat
    \end{enumerate}
\end{enumerate}

Agaipiou et al. \cite{agapiou2014analysis} also considered a non-centered parameterization in a Bayesian Gaussian linear inverse problem, similar to that considered in this work. They relied on Metropolis sampling, but with a different proposal mechanism. Instead of sampling $\sigma$ directly, they considered $\zeta := \sigma^{-1/2}$. After a change of variables, the full conditional distribution of $\zeta$ is
\begin{eqnarray*}
    \pi(\zeta \mid \B{b}, \B{z}, \mu) &\propto& f(\B{b} \mid \B{z}, \mu, \zeta )\pi(\zeta)\\
        &\propto& \exp\left[-\frac{\mu}{2}(\B{b} - \zeta\B{Az})^\top(\B{b} - \zeta\B{Az})\right]\left(\frac{1}{\zeta^2}\right)^{a_{\sigma} + 1/2}e^{-b_{\sigma}/\zeta^2}.
\end{eqnarray*}
We see that the likelihood contribution to this density, written as a function of $\zeta$, is proportional to a Gaussian density. That is,
\begin{equation}\label{eqn:ncpProp}
    g(\zeta) := \exp\left[-\frac{\mu}{2}(\B{b} - \zeta\B{Az})^\top(\B{b} - \zeta\B{Az})\right] \propto \exp\left[-\frac{\mu\B{z}^\top\B{A}^\top\B{Az}}{2}\left(\zeta - \frac{\B{z}^\top\B{A}^\top\B{b}}{\B{z}^\top\B{A}^\top\B{Az}}\right)^2\right],
\end{equation}
which is the density of a normal distribution with mean $\zeta_c := \B{z}^\top\B{A}^\top\B{b}(\B{z}^\top\B{A}^\top\B{Az})^{-1}$ and variance $\xi_c := (\mu\B{z}^\top\B{A}^\top\B{Az})^{-1}$.  Agaipiou et al. \cite{agapiou2014analysis} use this Gaussian distribution as a proposal for an independence sampler, except that $\zeta$ is restricted to be positive. In other words, their proposal distribution is a truncated Gaussian with density
\[
    q(\zeta^{\ast}) = \cN \left(\zeta^{\ast} \mid \zeta_c, \xi_c\right)\left(1 - \Phi(-\zeta_c/\sqrt{\xi_c})\right)^{-1},
\]
where $\Phi(\cdot)$ is the cumulative distribution function of the standard normal distribution. It is important to note that this approach assumes $\B{A}$ is of full rank, so that $\B{Az} \neq \B{0}$ for all $\B{z} ~(a.e.)$. In cases where $\B{A}$ is rank deficient, one can use the random walk Metropolis step previously discussed, but other approaches are possible.

\subsection{Illustration with Image Deblurring} 
\begin{figure}[tb]
    \centering
    \includegraphics[scale= 0.4]{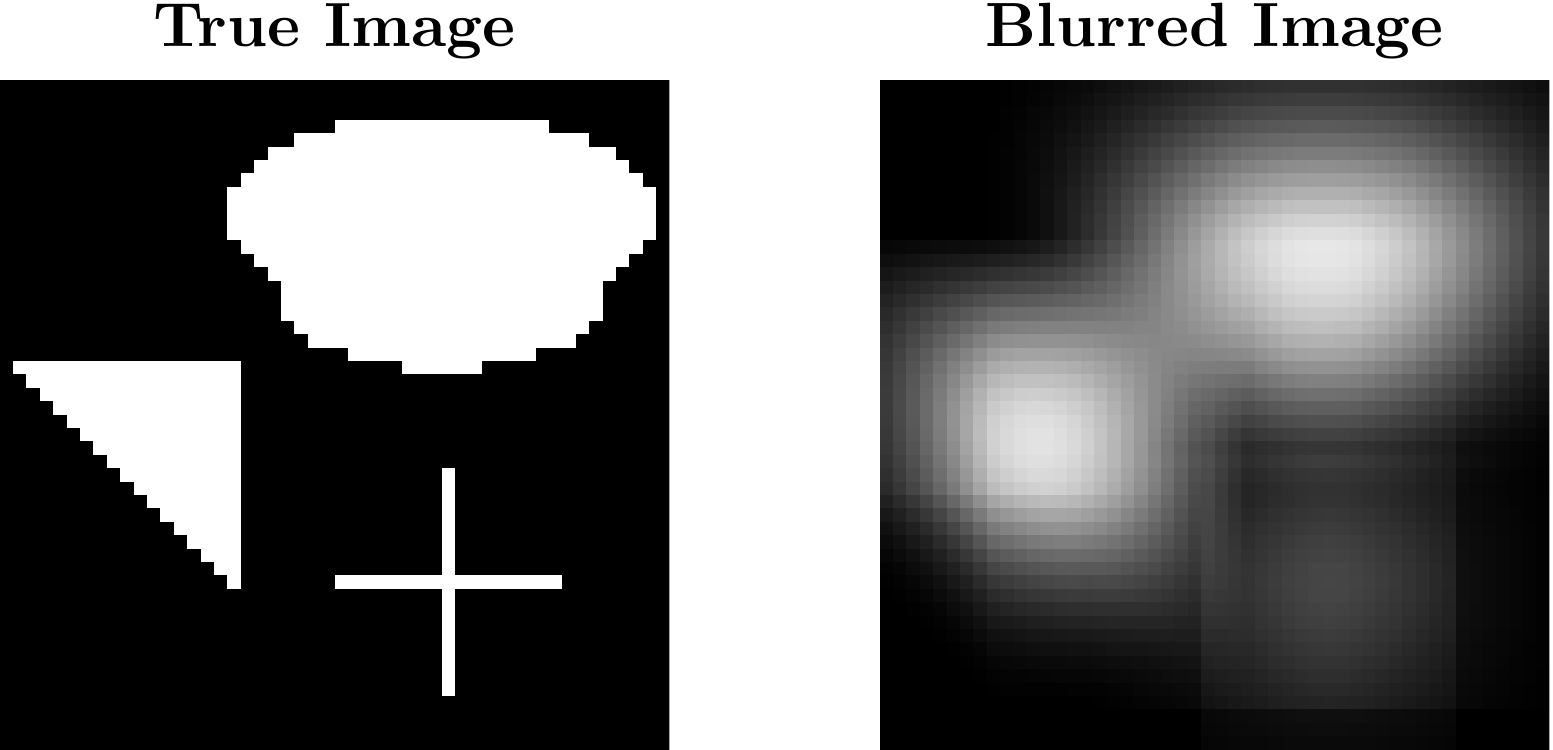} \hspace{1.5em} \includegraphics[scale= 0.4]{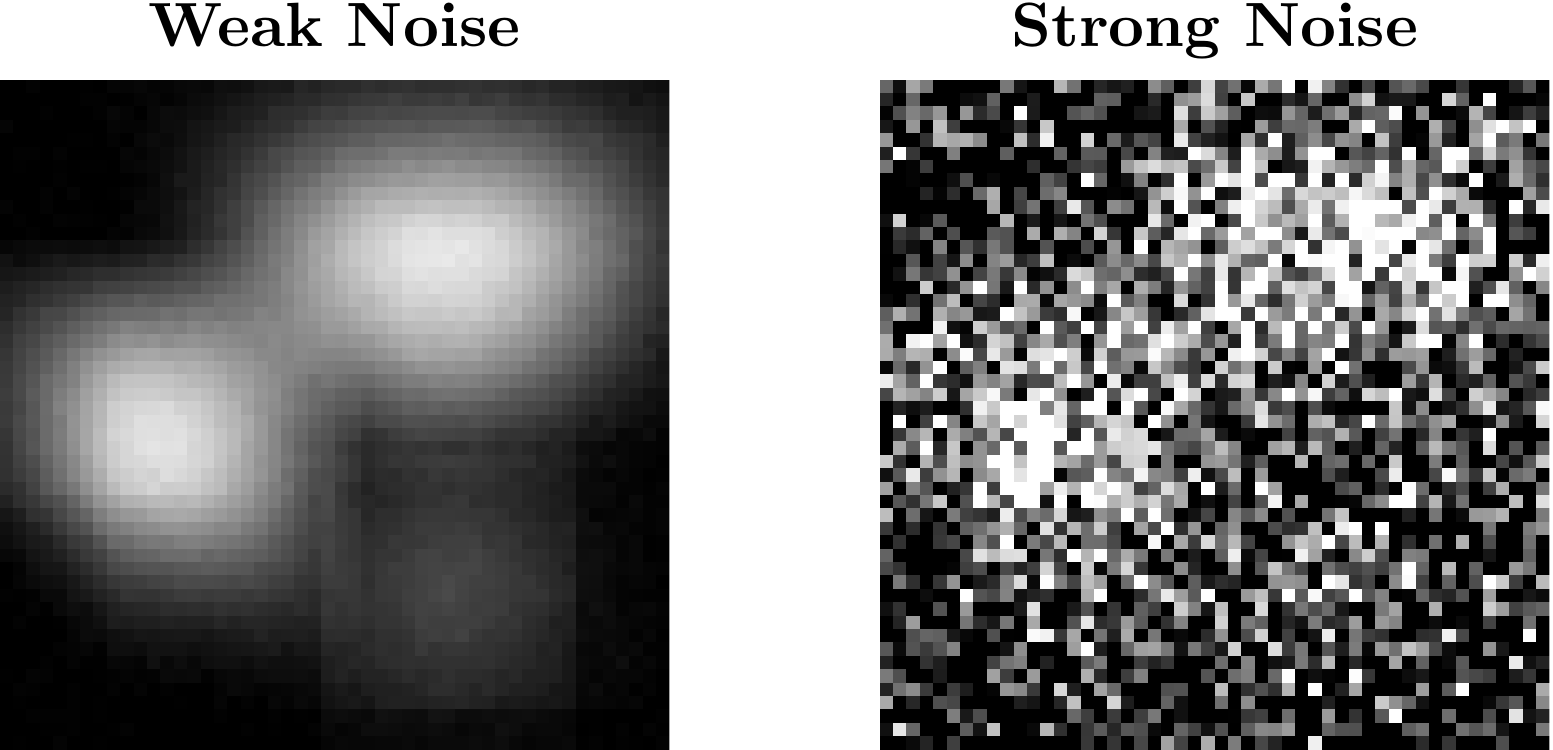}
    \caption{From left to right: target image, blurred image, reliable data, and corrupted data, used for comparing noncentered and centered parameterizations under low-rank independence sampling for 2D deblurring.}
    \label{fig:NCPvsCPTargets}
\end{figure}

To illustrate use of the non-centered parameterization in our proposed LRIS algorithm, we consider again the 2D image deblurring example from Section 4.1 of the manuscript. Here, the target image $\B{x}$ and the blurring operator $\B{A}$ are the same as before. However, we create two different observed datasets $\B{b}_i = \B{Ax} + \B{\epsilon}_i, ~~i= 1,2$, with two different levels of noise. One set is strongly corrupted with 50\% noise, $Var(\B{\epsilon}_1) = 0.5^2\|\B{b}\|_{\infty}^2\B{I}$, and the other contains much less noise, $Var(\B{\epsilon}_2) = 0.01^2\|\B{b}\|_{\infty}^2\B{I}$, and thus is more informative about the true solution $\B{x}$. The target image, blurred image, and noisy data sets are displayed in Figure \ref{fig:NCPvsCPTargets}.

For the centered parameterization, we use the same priors on $\B{x}, \mu$ and $\sigma$ in model \eqref{EQ:BIPMod} as in the manuscript, namely $\prior^{-1} = \B{L}^{\top}\B{L}$ with $\B{L} = -\Delta + \delta\B{I}$ and $\epsilon = 0.1$ in the $\text{Gamma}(\epsilon, \epsilon)$ priors on $\mu$ and $\sigma$ to approximate scale invariant objective priors. We use the same low-rank Metropolis-Hastings-within-Gibbs algorithm as in the manuscript. To implement the non-centered parameterization, we use our proposed low-rank independent sampler to update $\B{x}$ and the independent sampler proposed by \cite{agapiou2014analysis} to update $\sigma^{-1/2}$. Our interest here is not in posterior inference about the target image, but in the mixing behavior of these two parameterizations applied to data with different amounts of noise. Thus, rather than running to and diagnosing convergence, we run only 5,000 iterations of each Markov chain to study autocorrelation and how quickly the chains appear to be moving through their support. Under each configuration, we run three chains in parallel, with $\mu$ and $\B{x}$ initialized by drawing them randomly from their prior distributions. The prior precision $\sigma$ is initialized at $0.1, 6$, and $25$ for chains 1, 2, and 3, respectively.

\begin{figure}[tb]
    \centering
    \includegraphics[scale= 0.75]{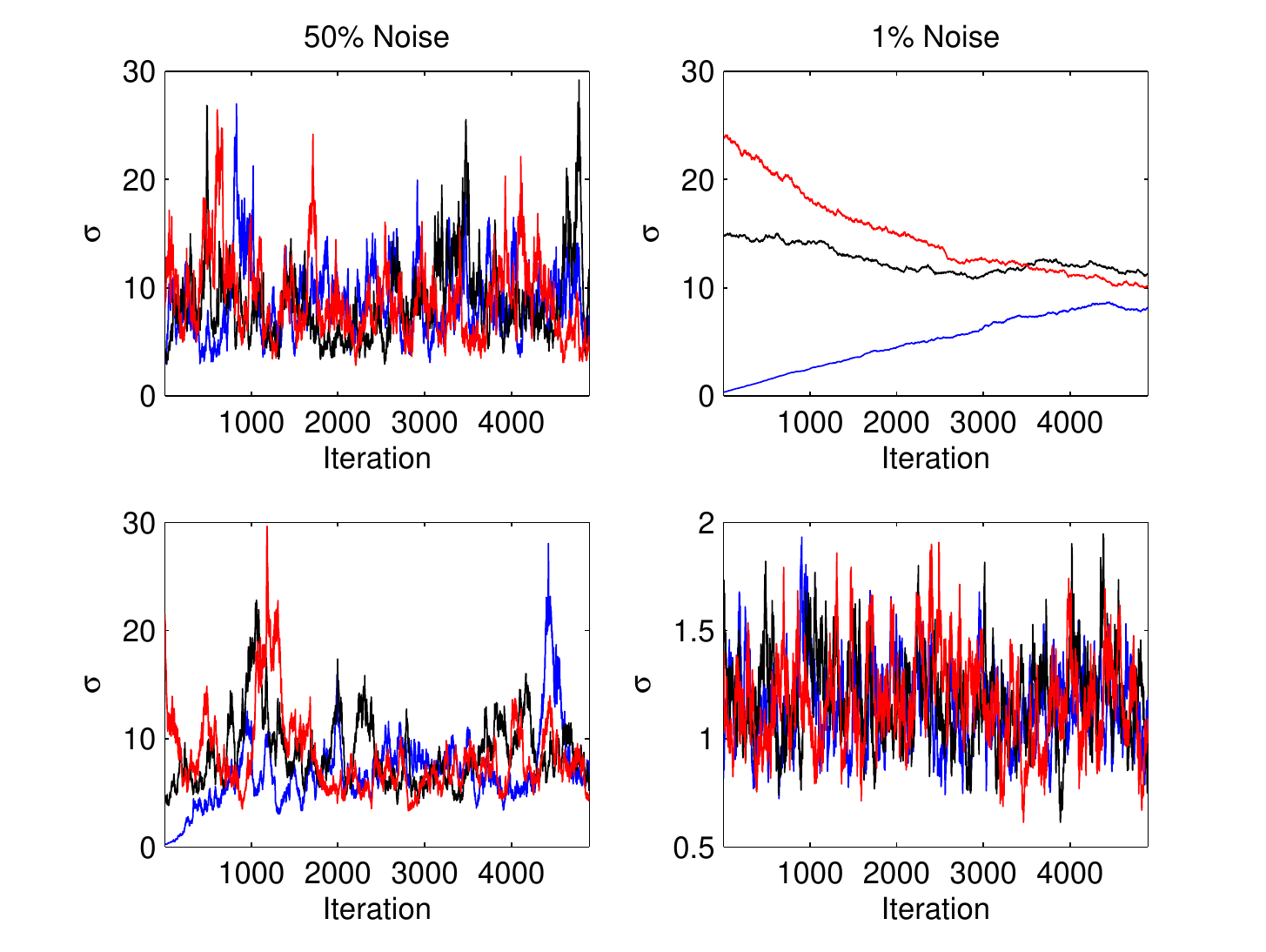}
    \caption{Trace plots of the $\sigma$ chains obtained from centered and non-centered parameterizations using both datasets displayed in Figure \ref{fig:NCPvsCPTargets}. The left column corresponds to the noisy data and right column to the reliable data. The top row are plots obtained from the noncentered parameterization, and the bottom row is obtained from the centered parameterization.}
    \label{fig:NCPvsCPTrace}
\end{figure}

\begin{figure}[!h]
    \centering
    \includegraphics[scale= 0.25]{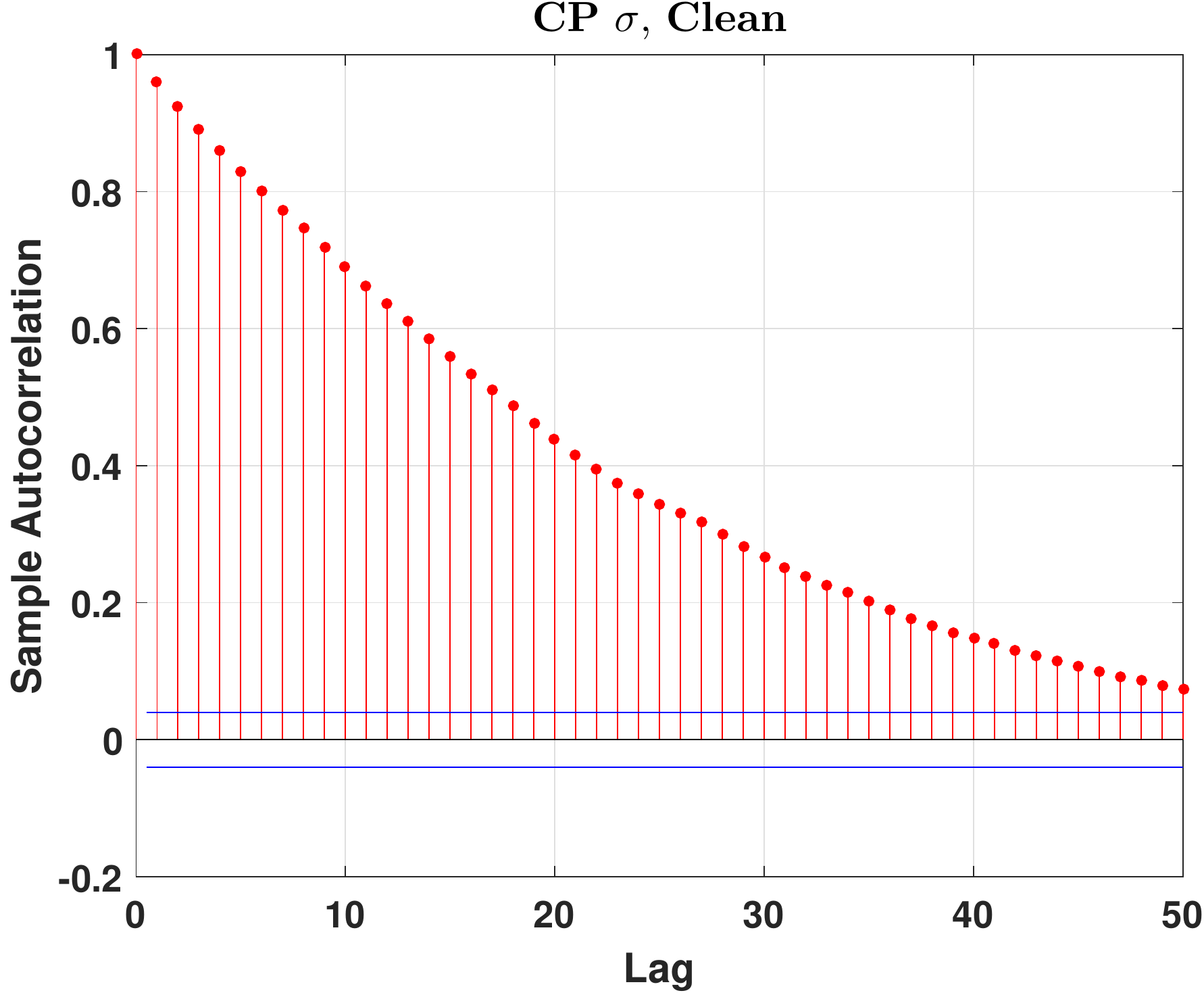}
    \includegraphics[scale= 0.25]{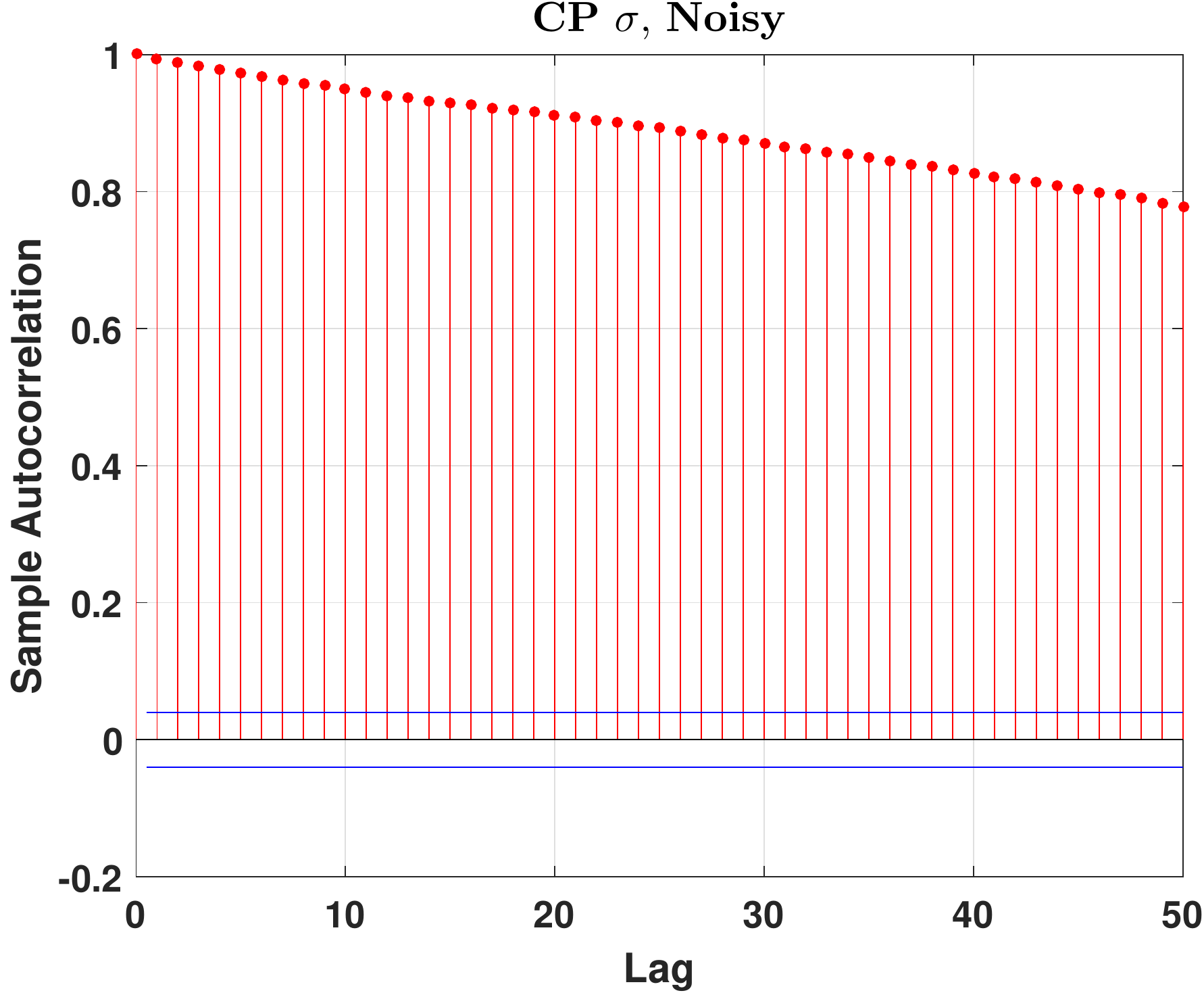}\\
    \includegraphics[scale= 0.25]{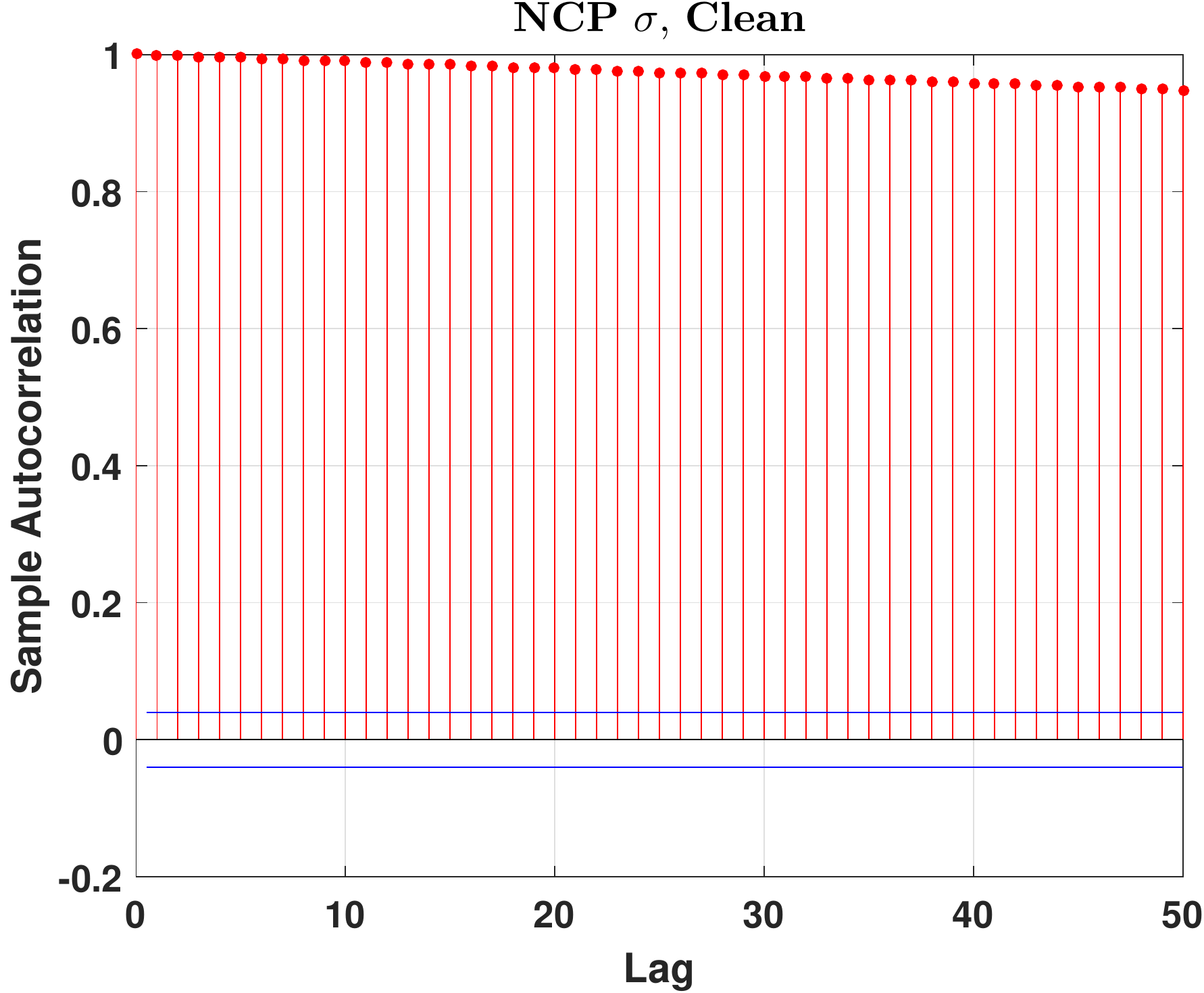}
    \includegraphics[scale= 0.25]{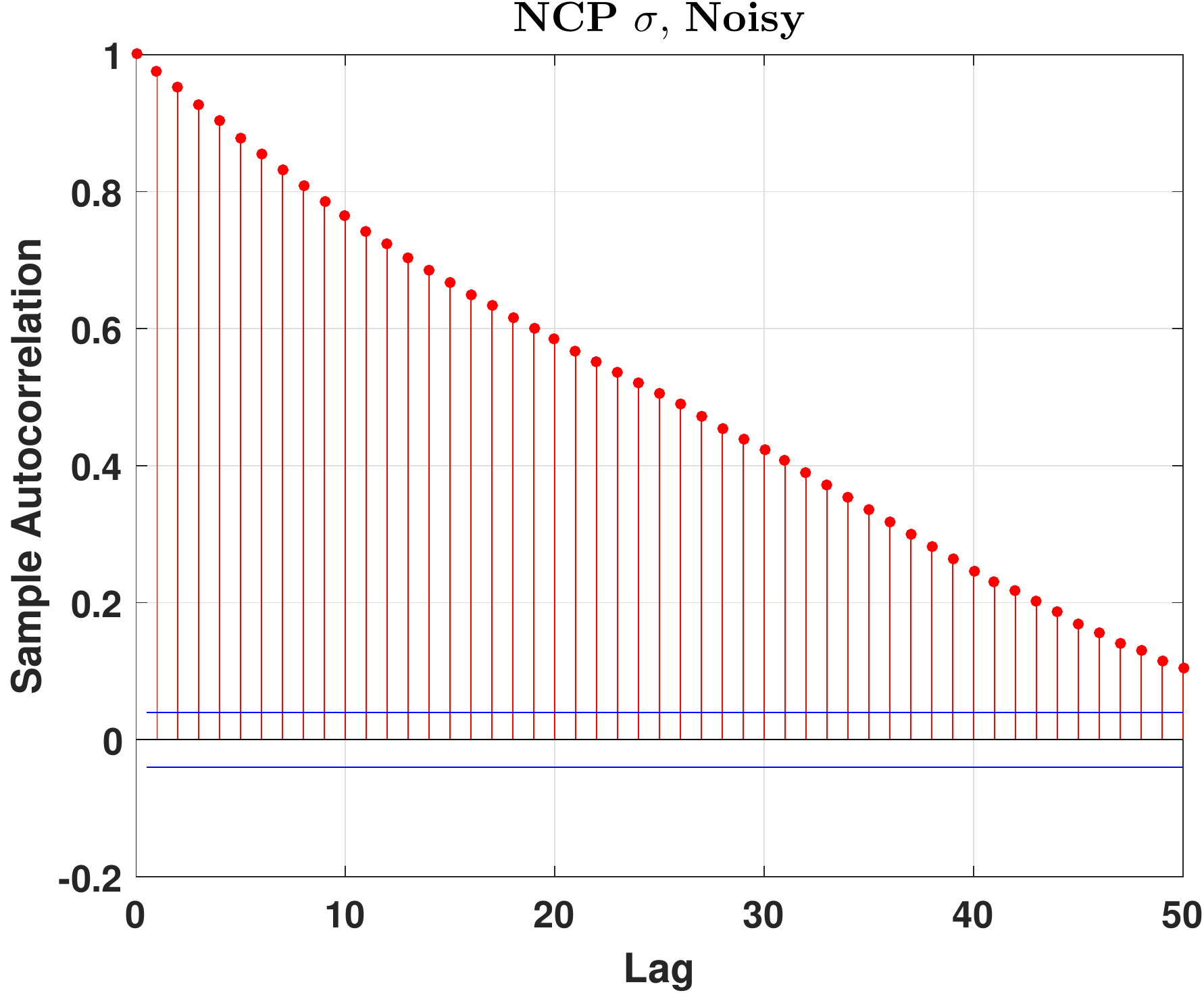}
    \caption{Estimated autocorrelation functions of one $\sigma$ chain from the Supplementary 2D deblurring simulation study, each obtained from centered (top row) and non-centered (bottom row) parameterizations using both weakly (left column) and strongly (right column) corrupted data.}
    \label{fig:NCPvsCPACF}
\end{figure}

Figure \ref{fig:NCPvsCPTrace} displays the trace plots of the three $\sigma$ chains under each of the four combinations of data and parameterization. With severely noisy data, the noncentered parameterization shows stronger mixing than the centered parameterization. The opposite is true with the reliable data containing only 1\% noise. The differences under the reliable data are particulary striking, where we see improvement in the centered paramterization but a very severe degradation in performance of the noncentered parameterization. The drift in all of the chains is indicative of considerable autocorrelation, and this is confirmed by examining the autocorrelation functions plotted in Figure \ref{fig:NCPvsCPACF} and the estimated lag 1 and lag 50 correlation coefficients in Table \ref{tab:NCPvsCPLags}. Each chain suffers from high lag 1 autocorrelation, but it decays faster under the non-centered parameterization for the noisy data, and faster for the centered parameterization under the more reliable data. The decay of the autocorrelation is particularly poor for the noncentered parameterization with the reliable data.
\begin{table}[!h]
    \centering
    \begin{tabular}{l | c c | c c}
     & \multicolumn{2}{c|}{\underline{Reliable}} & \multicolumn{2}{c}{\underline{Noisy}}\\
    & Lag 1 & Lag 50 & Lag 1 & Lag 50\\
    \hline
    CP & 0.968 & 0.174 & 0.995 & 0.774 \\
    NCP  & 0.999 & 0.974 & 0.982 & 0.348\\
    \end{tabular}
    \caption{Estimated autocorrelation coefficients of one $\sigma$ chain, for each data / parameterization combination in Supplementary simulation study.}
    \label{tab:NCPvsCPLags}
\end{table}

The results of this illustration demonstrate the ease with which either the non-centered or centered parameterization can be used in combination with our proposed LRIS. Further, the relative performance of non-centered versus centered parameterizations previously observed in \cite{Papaspiliopoulos2003, Papaspiliopoulos2007, agapiou2014analysis} is still present when using our more computationally efficient alternative. In particular, for strongly corrupted data, $\B{x}$ is more strongly determined through the prior than the data, so that a non-centered parameterization is preferable. More reliable data impose the constraint that $\B{b} \approx \sigma^{-1/2}\B{Az}$, severely degrading the performance of the non-centered parameterizaton. Thus, when using our proposed approach, a practitioner can still appeal to the same considerations when choosing a more effective parameterization for convergence of their MCMC algorithm.

\section{Adaptive LRIS}
The target rank $k$ for the low-rank approximation may not be known in practice, or may depend explicitly on the parameters $\mu$ and $\sigma$. Here we outline a simple adaptive strategy for determining the target rank. The basic idea is to increment the rank till the acceptance ratio meets the desired tolerance.
\begin{enumerate}
\item Inputs: user defined tolerance $\delta$, initial rank $k_\text{init}$, rank increment $k_\text{inc}$
\item Start $k = k_\text{init}$
\item While not \texttt{converged}
\item Check if acceptance ratio is higher than $1-\delta$. If, yes, mark \texttt{converged}.
\item Increment the rank $k \leftarrow k + k_\text{inc}$.
\item End While
\end{enumerate}
Monitoring the acceptance ratio is important to determine a stopping criterion for the algorithm described above. We suggest several different strategies:
\begin{enumerate}
\item The acceptance ratio can be monitored empirically as the independence sampler is exploring the distribution. If the acceptance ratio is too small, the target rank may be incremented.
\item In the case that the exact eigenpairs are used, we derive a lower bound for $w(\B{x})$. From the proof of Theorem 1, we find that if $\widehat{\B{H}} = \B{V}_k\B\Lambda_k\B{V}_k^\top$
\[ w(\B{x}) \geq \exp\left(-\frac{\mu}{2}\|\B{Lx}\|_2^2 \| \B{H}- \widehat{\B{H}}\|_2\right) = \exp\left(-\frac{\mu\lambda_{k+1}}{2}\|\B{Lx}\|_2^2 \right). \]
Here $\lambda_{k+1}$ is the largest eigenvalue that is discarded. This suggests that the target rank $k$ is the minimum index which satisfies \[ \lambda_{k+1} \leq \frac{2}{\mu\|\B{Lx}\|_2^2}\log\frac{1}{1-\delta}.\]
This will ensure $w(\B{x}) \geq 1 -\delta$, where $\delta$ is a user-defined parameter. It is worth mentioning that the entire eigendecomposition need not be computed, nor recommended. In practice, the eigenpairs can be computed in an incremental fashion.

\item In the randomized low-rank approach, $\lambda_{k+1}$ may not be available. However, it may be easy to estimate $\| \B{H}- \widehat{\B{H}}\|_2$ cheaply using a randomized estimator; see~\cite[Lemma 4.1]{halko2011finding}. Then to ensure that $w(\B{x}) \geq 1 -\delta$, it is required that $\| \B{H}- \widehat{\B{H}}\|_2 \leq \varepsilon$, where $\varepsilon$ satisfies
\[ \varepsilon \leq \frac{2}{\mu\|\B{Lx}\|_2^2}\log\frac{1}{1-\delta}.\]
The low-rank approximation can be adaptively and efficiently determined using the adaptive range finding algorithm, see~\cite[Algorithm 4.2]{halko2011finding}.
\end{enumerate}

\section{NMR Relaxometry Simulation}

In this section, we investigate the performance of our approach on a large-scale ill-posed inverse problem with a different low-rank nature than the CT example. We consider the problem of nuclear magnetic resonance (NMR) relaxometry in which nuclear magnetic moments are used to infer physical or chemical properties of a medium. The deterioration (or \emph{relaxation}) of the signal over time is analyzed using various electromagnetic pulse sequences and acquisition strategies. The longitudinal and traverse relaxation times of the medium, denoted $T_1$ and $T_2$, respectively, can be recovered from these acquisitions. The goal of 2D NMR relaxometry is to reconstruct the joint distribution of the relaxation times, $f(T_1, T_2)$, from measured signals gathered at different experimental times, denoted $g(\tau_1, \tau_2)$. Mathematically, the forward model is a Fredholm integral equation of the first kind, $g(\tau_1, \tau_2) = \int_{0}^{\hat{T}_2} \int_{0}^{\hat{T}_1} K(\tau_1, \tau_2, T_1, T_2) f(T_1, T_2) \, \dx T_1 \dx T_2.$ As in the 1D image deblurring example, the measurement noise is typically assumed to be Gaussian and the problem is to reconstruct $f$ from $g$ and $K$.

\begin{figure}[tb]
	\includegraphics[scale=0.25]{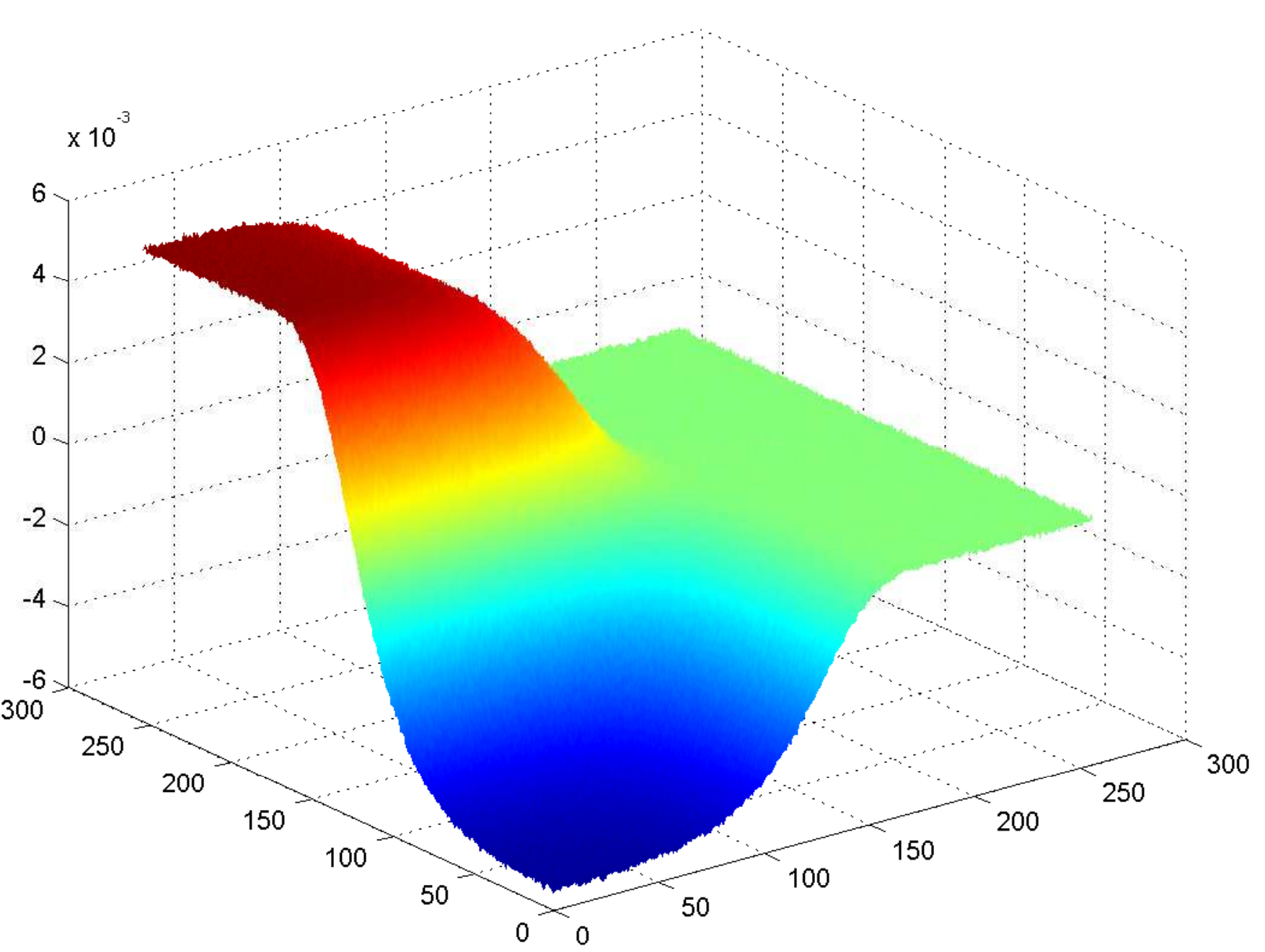}
	\includegraphics[scale=0.25]{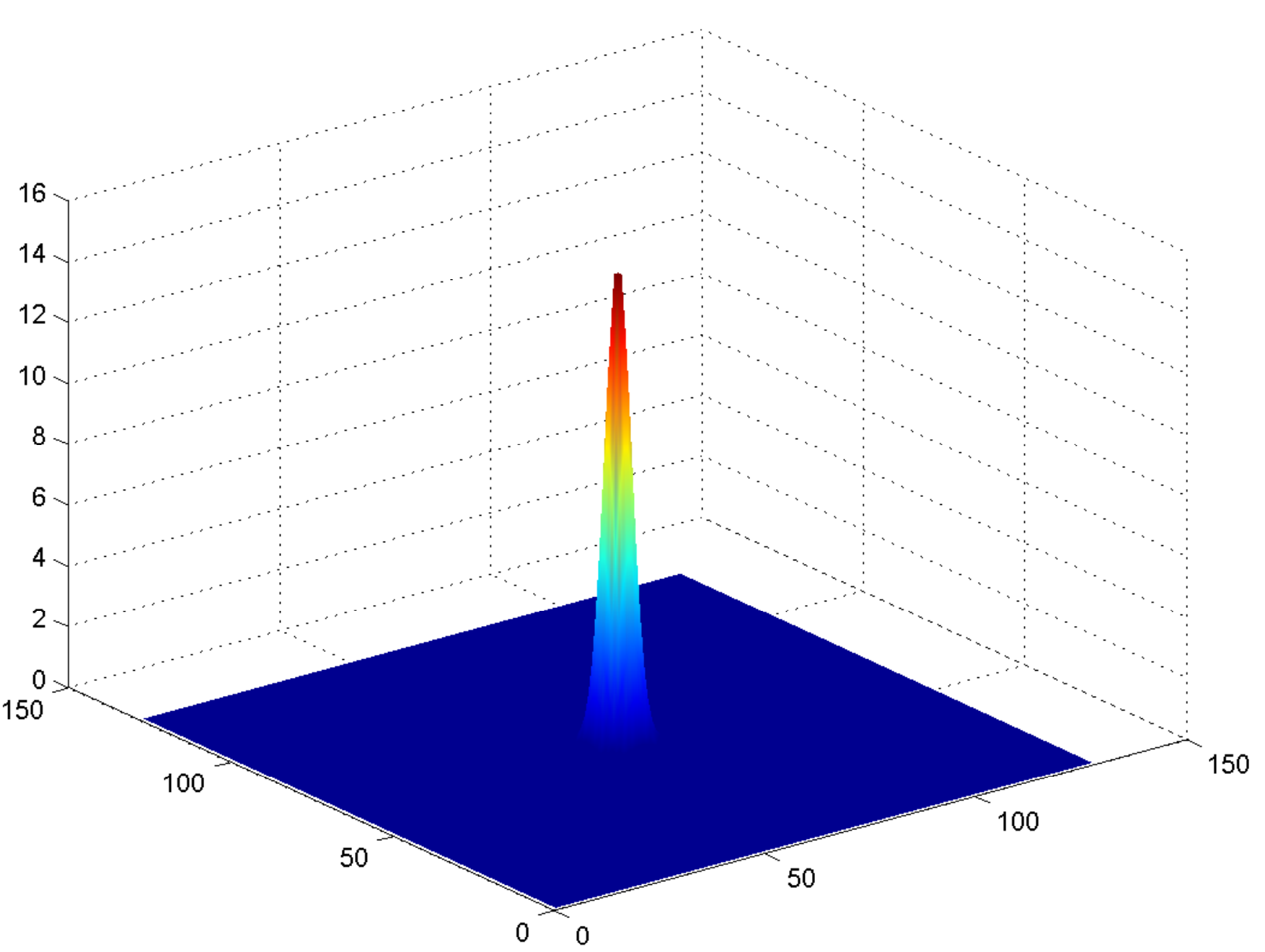}
	\includegraphics[scale=0.25]{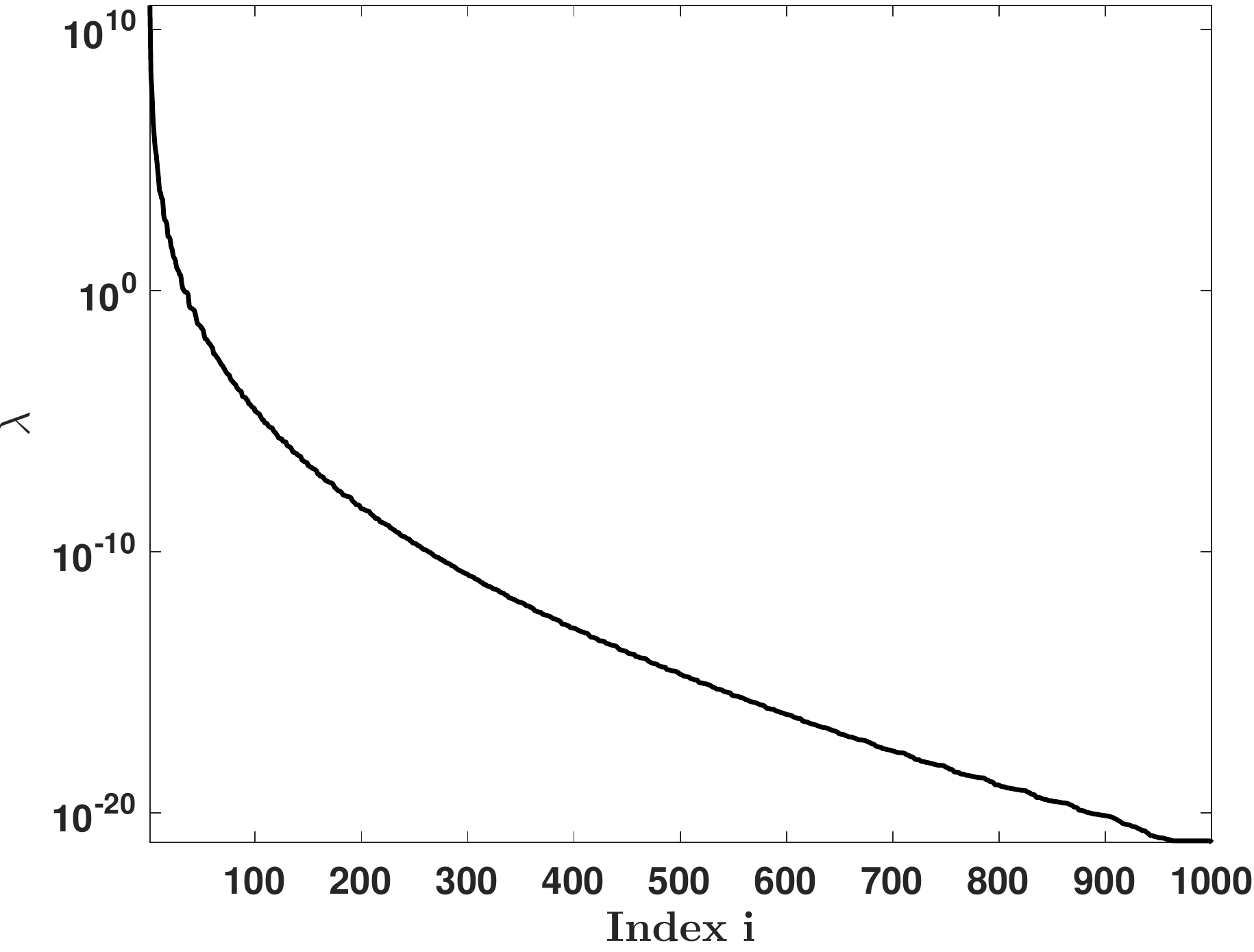}
	\caption{Noisy measured data $\B{b}$ (left) and true relaxation distribution $\B{x}$ (center) in the NMR relaxometry example. Approximate spectrum of $\B{H}$ (right) shows significant decay within the first few hundred eigenvalues.}
	\label{fig:NMRsetup}
\end{figure}

We use the implementation in the \texttt{Matlab} package \texttt{IR Tools} \cite{Gazzola17} to simulate NMR data. In this example, the kernel is separable, and after discretization we obtain the linear inverse problem $\B{Ax} = \B{b}$, where $\B{b}$ is the vectorization of the data $g$ and $\B{x}$ the vectorization of the unknown $f$. Unlike the other examples we consider in this work, the forward map $\B{A}$ and its transpose are not constructed explicitly, but their matrix-vector products are computed using available function handles. We again enforce smoothness in $\B{x}$ by specifying $\B{L} = -\Delta + \delta\B{I}$ and corrupt the data with one percent Gaussian random noise. With the default discretizations, the data have dimension $m=65,536$ and the unknown has dimension $n=16,384$. The noisy data $\B{b}$ and unknown $\B{x}$ are displayed in Figure~\ref{fig:NMRsetup}.  We compute approximate eigenvalues of the prior preconditioned Hessian using Randomized SVD with $\ell = 1000$. The sharp decay in the eigenvalues, also displayed in Figure~\ref{fig:NMRsetup}, illustrates the severe ill-posedness of the NMR relaxometry inverse problem.

\begin{figure}[tb]
	\includegraphics[scale=0.25]{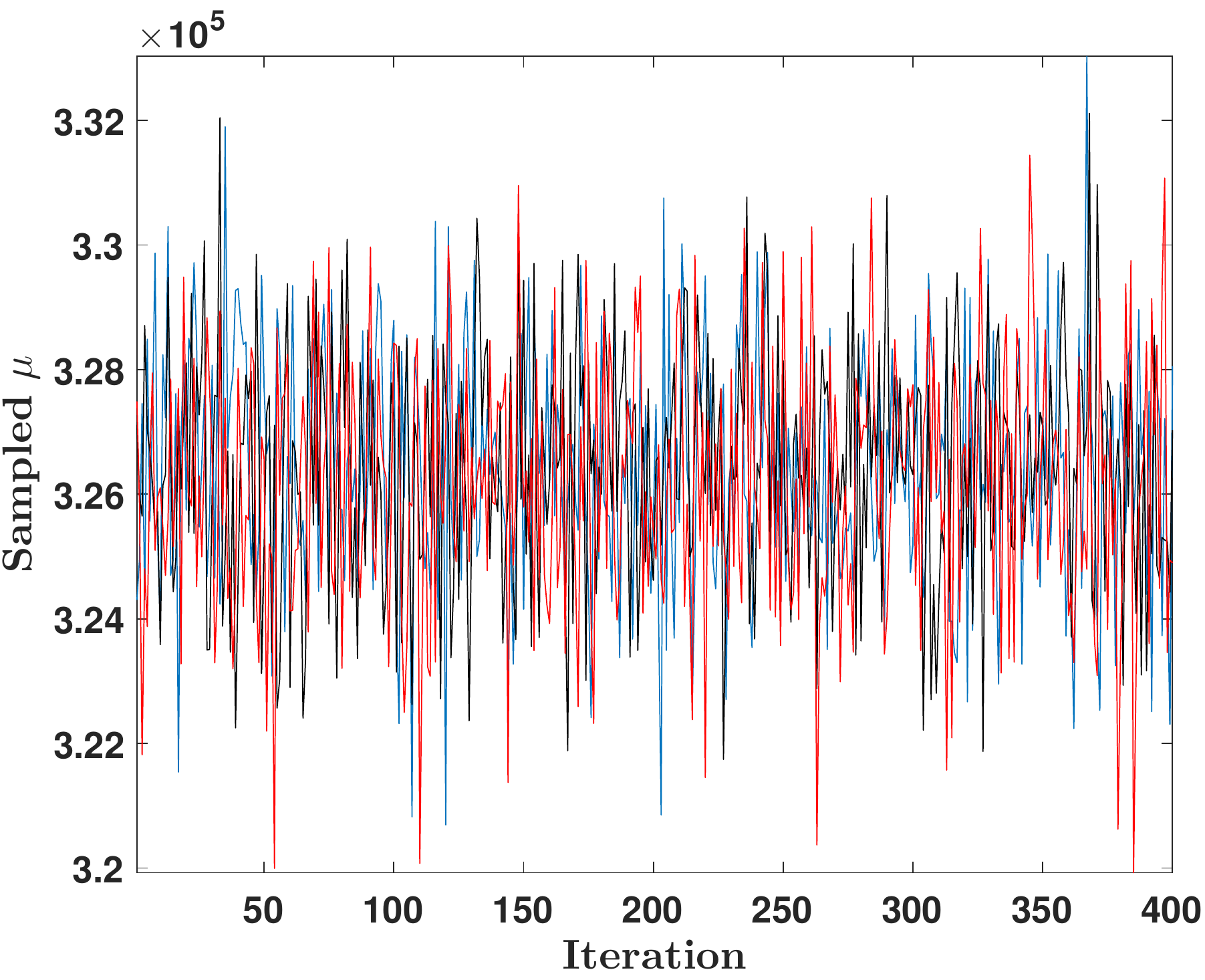}
	\includegraphics[scale=0.25]{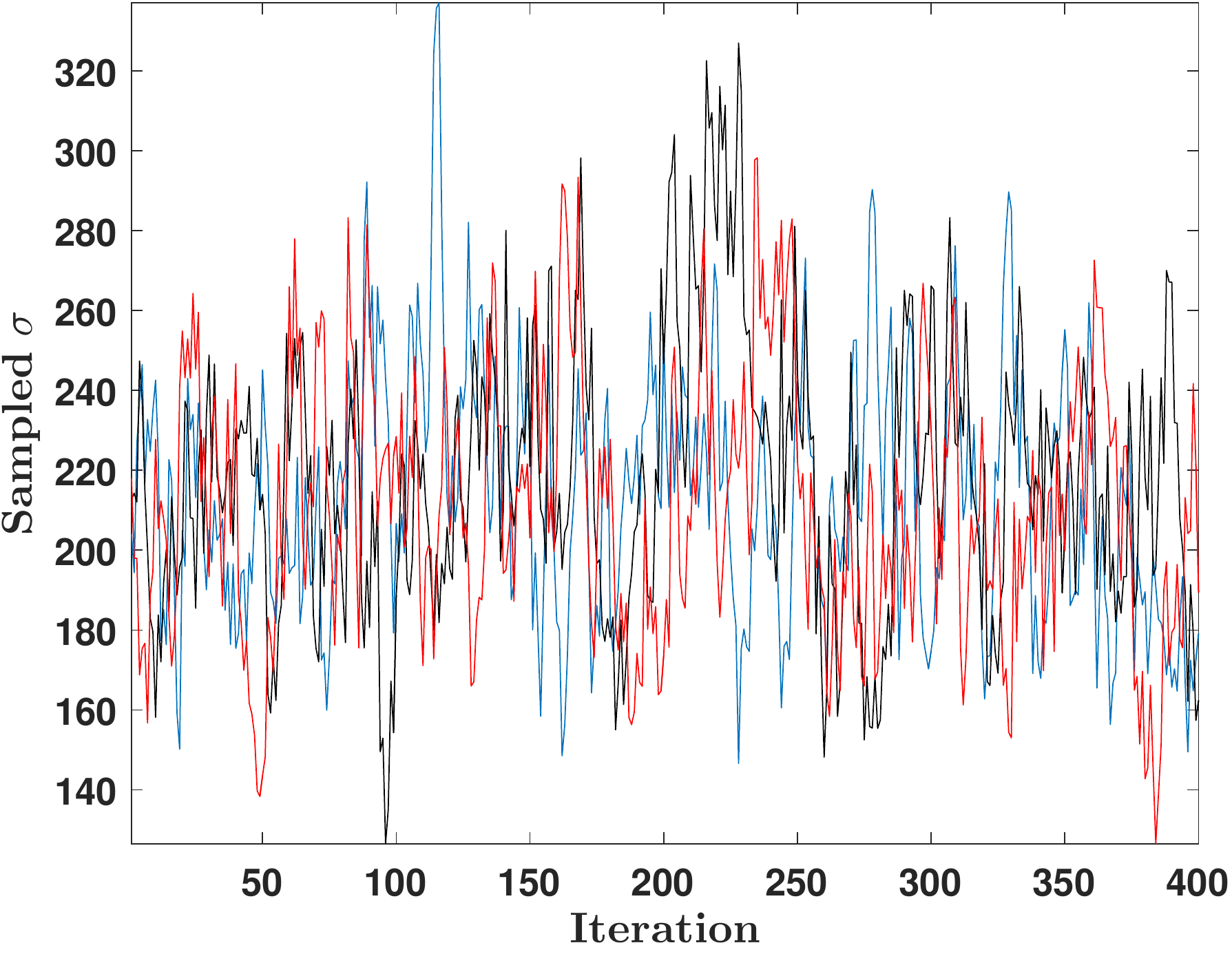}
	\includegraphics[scale=0.25]{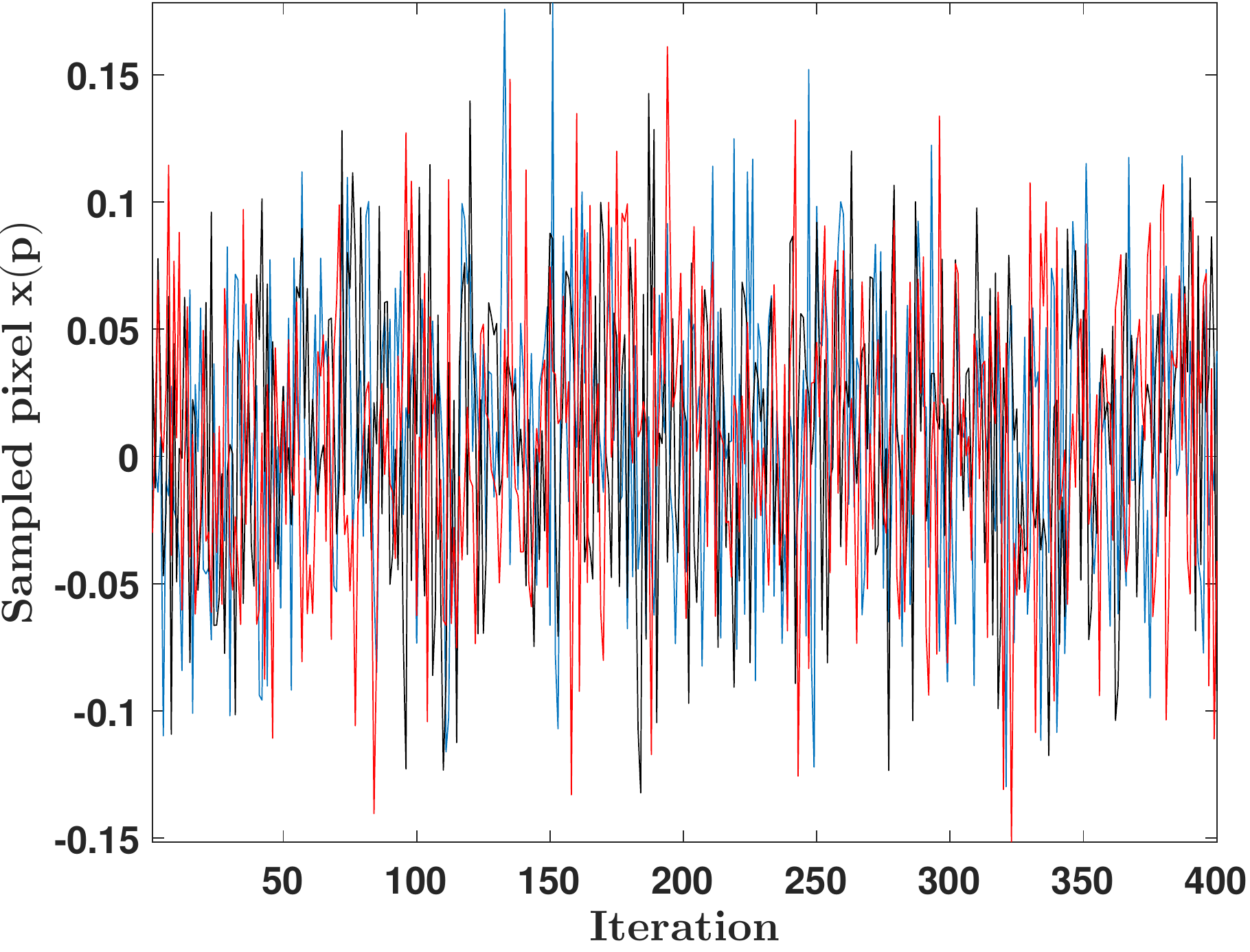}
	\caption{Trace plots for the three thinned chains using Gamma priors in the NMR relaxometry example. Left: noise precision $\mu$. Center: prior precision $\sigma$. Right: a randomly chosen pixel of the image $\B{x}$.}
	\label{fig:NMRtrace}
\end{figure}

\begin{figure}[tb]
	\includegraphics[scale=0.3]{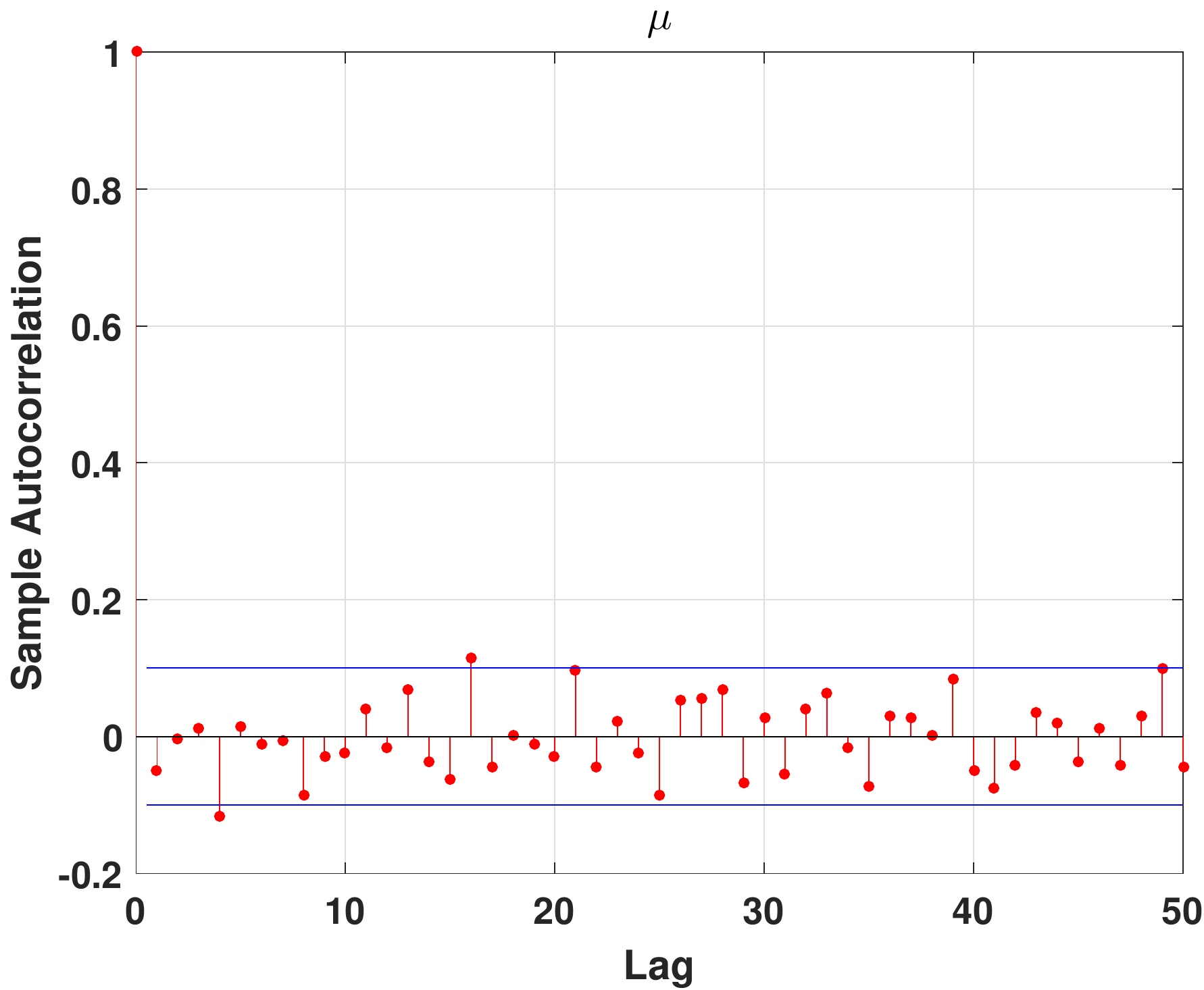}
	\includegraphics[scale=0.3]{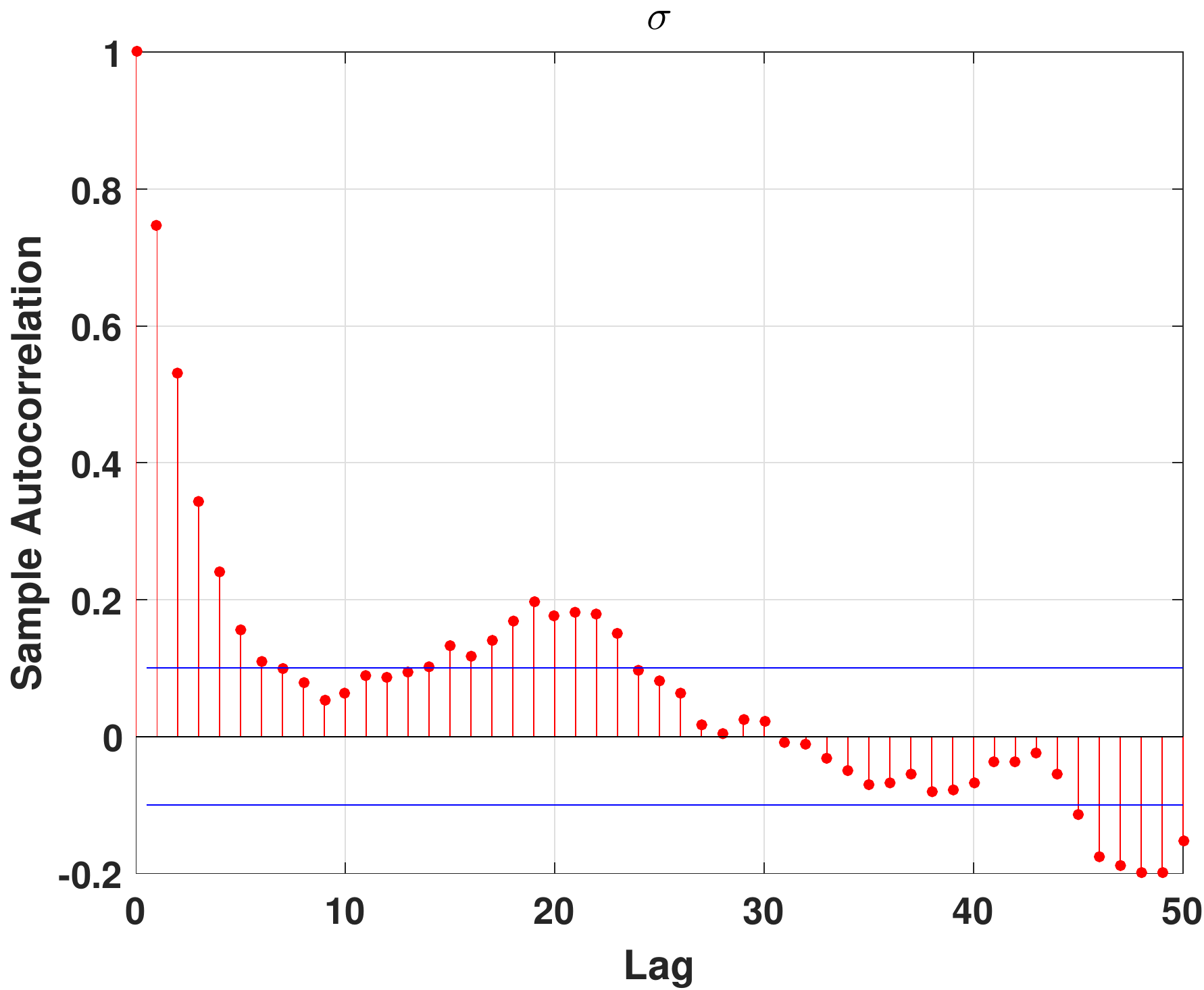}
	\caption{Autocorrelation plots for one thinned chain using Gamma priors in the NMR relaxometry example. Left: noise precision $\mu$. Right: prior precision $\sigma$.}
	\label{fig:NMRautocorr}
\end{figure}

For implementation, we perform the same experiment as in Section 4.2 of the manuscript with the same vague conjugate Gamma priors on $\mu$ and $\sigma$. Analyzing the acceptance rates for different eigenvalue truncation levels as well as the spectral decay in Figure~\ref{fig:NMRsetup}, we specify the truncation level $k=300$ (average acceptance rate $\approx 100\%$). We simulate three Markov chains for 40,000 iterations each, initializing each chain with random draws from the priors. To reduce the autocorrelation in the chains, they are thinned to retain every 50th draw, after which the first 400 draws are discarded as the burn-in period. Approximate convergence of the chains is diagnosed via trace plots and autocorrelation plots, displayed in Figures~\ref{fig:NMRtrace} and~\ref{fig:NMRautocorr}. The PSRFs for $\mu$ and $\sigma$ are 0.99 and 1.01, respectively. The computations are carried out in \texttt{MATLAB 2013a} on a Dell Precision T3600 desktop PC running Windows 7 Enterprise with an Intel Xeon E5-1660 3.30GHz CPU and 64GB RAM. The total computation time is 13,263 seconds, or about 3.7 hours. %This is significantly faster than the CT example in Section 4.2 of the manuscript, which has the same dimensionality but requires a higher truncation level $k$.

\begin{figure}[tb]
	\includegraphics[scale=0.4]{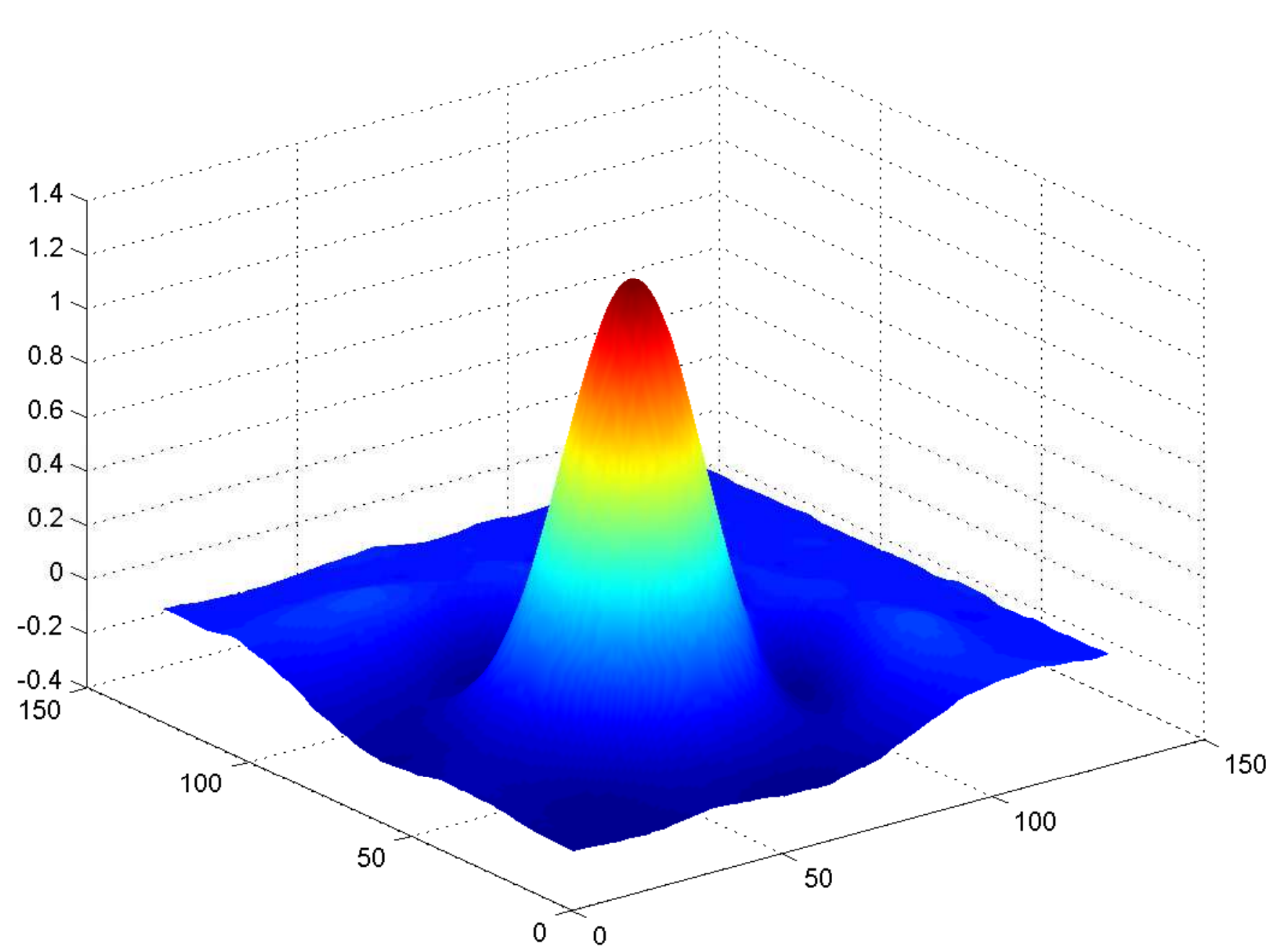}\hspace{18pt}
    \includegraphics[scale=0.35]{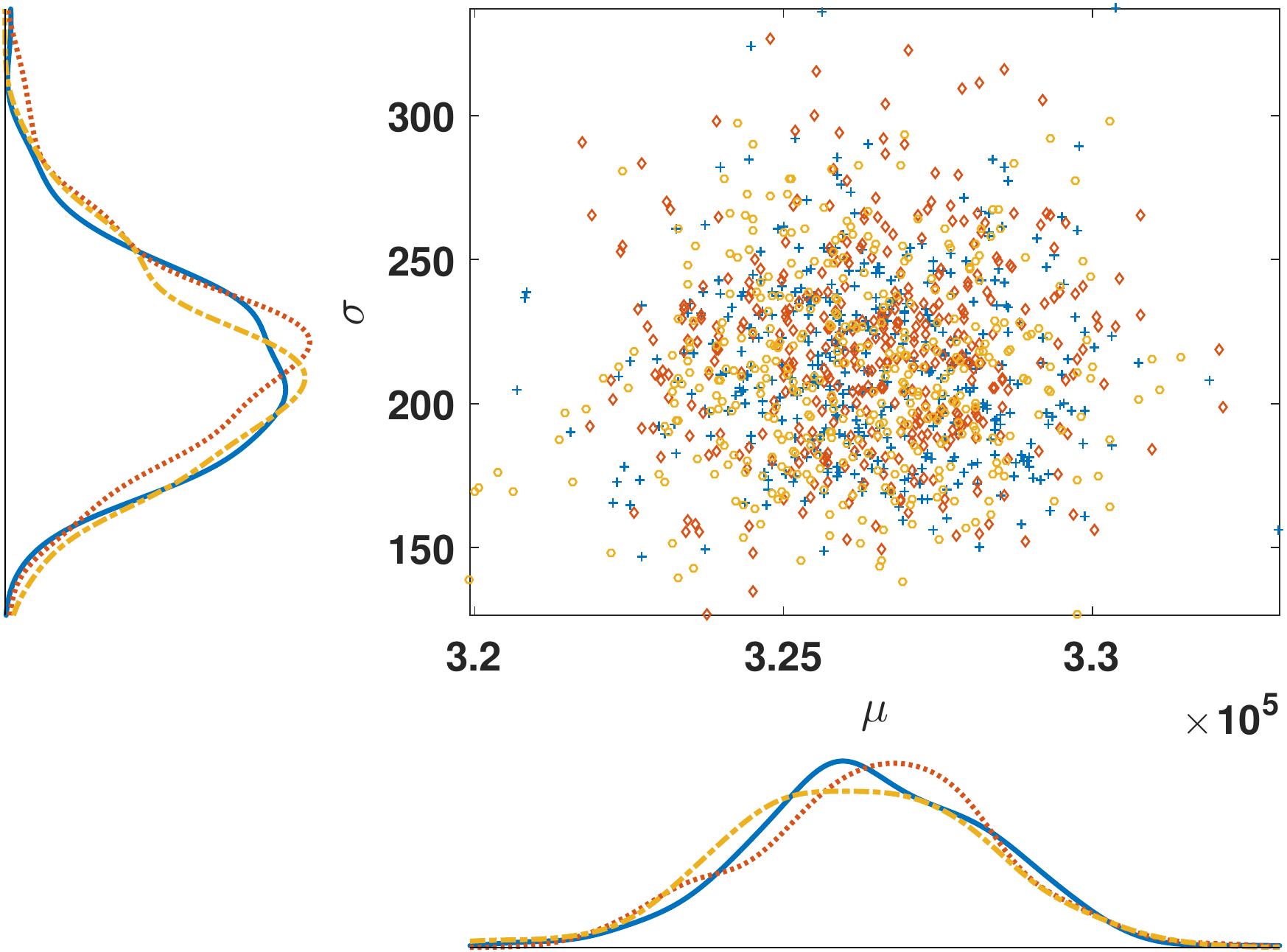}
	\caption{Posterior mean estimate of $\B{x}$ (left) and smoothed histograms of $\mu$ and $\sigma$ (right) in the NMR relaxometry example.}
	\label{fig:NMRpostmean}
\end{figure}

Figure~\ref{fig:NMRpostmean} displays the posterior mean of $\B{x}$ based upon the output of the MCMC at approximate convergence, as well as the approximate densities of $\mu$ and $\sigma$. While we obtain a reasonable estimate of the distribution of relaxation times, it is strongly attenuated compared to the true solution and exhibits background fluctuations (including negative values). The attenuation, though, is similar to the solution obtained in \cite{Gazzola17} via conjugate gradient least squares, though our posterior mean produces a smoother estimate. Indeed, \cite{Gazzola17} remark that the NMR problem is extremely challenging and that even the iterative methods they consider require many thousands of iterations to produce interpretable estimates of $\B{x}$.

Our low-rank approach offers considerable computational advantages over standard MCMC, thus making it feasible to reconstruct a competitive solution along with appropriate measures of uncertainty. It is demonstrated in \cite{Gazzola17} that non-negativity constraints on $\B{x}$ can yield an improved solution. In our framework, this suggests that a truncated Gaussian or lognormal prior on $\B{x}$ may produce solutions superior to that produced by the Gaussian prior. Such an investigation is beyond the scope of this work and left for future research.

%Once approximate convergence has been attained, we contextualize our results by comparing with known behavior of deterministic solutions to the problem. The inverse problem in 2D NMR relaxometry is quite challenging to solve, and box constraints are commonly used to improve the quality of the reconstruction. Figure 10 in \cite{Gazzola17} illustrates this difficulty: the unconstrained solution contains fluctuations in the background (including negative values) and does not attain the correct values at the peak. The posterior mean estimate obtained with our sampling approach, displayed in Figure~\ref{fig:NMRpostmean}, suffers similarly from background fluctuations and loss of peak amplitude. Our low-rank approximation approach offers significant computational advantages, even for such challenging ill-posed Bayesian inverse problems. The nonnegativity constraints applied in \cite{Gazzola17} improve reconstruction significantly in this example. This suggests that a lognormal prior on $\B{x}$ may perform better than the Gaussian prior necessary for our approach. However, pursuit in this direction is beyond the scope of this work.

\section{Supplementary Figures}

%\subsection{2D Deblurring Example}
In Supplementary Figure \ref{fig:2Dmarginals}, we compare the (vague) Gamma priors on $\mu$ and $\sigma$ with the marginal posterior distributions estimated from the MCMC samples of the joint posterior.  The prior for $\mu$ has been scaled by a factor of 10,000 for visualization. Note the significant difference between the distributions, indicating that strong Bayesian learning has occurred.
\begin{figure}[!htb]
    \centering
    \includegraphics[scale= 0.3]{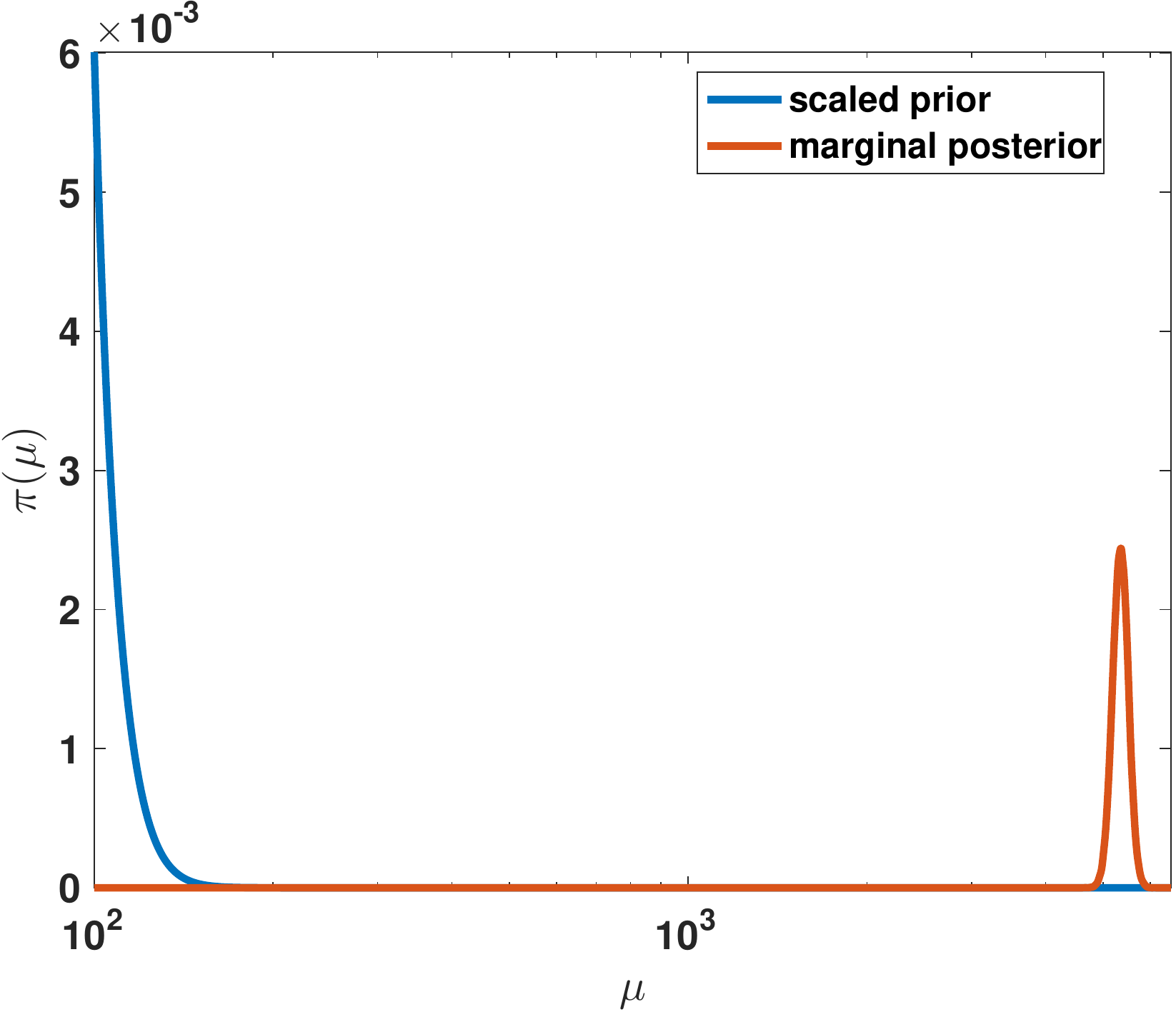}
    \includegraphics[scale= 0.3]{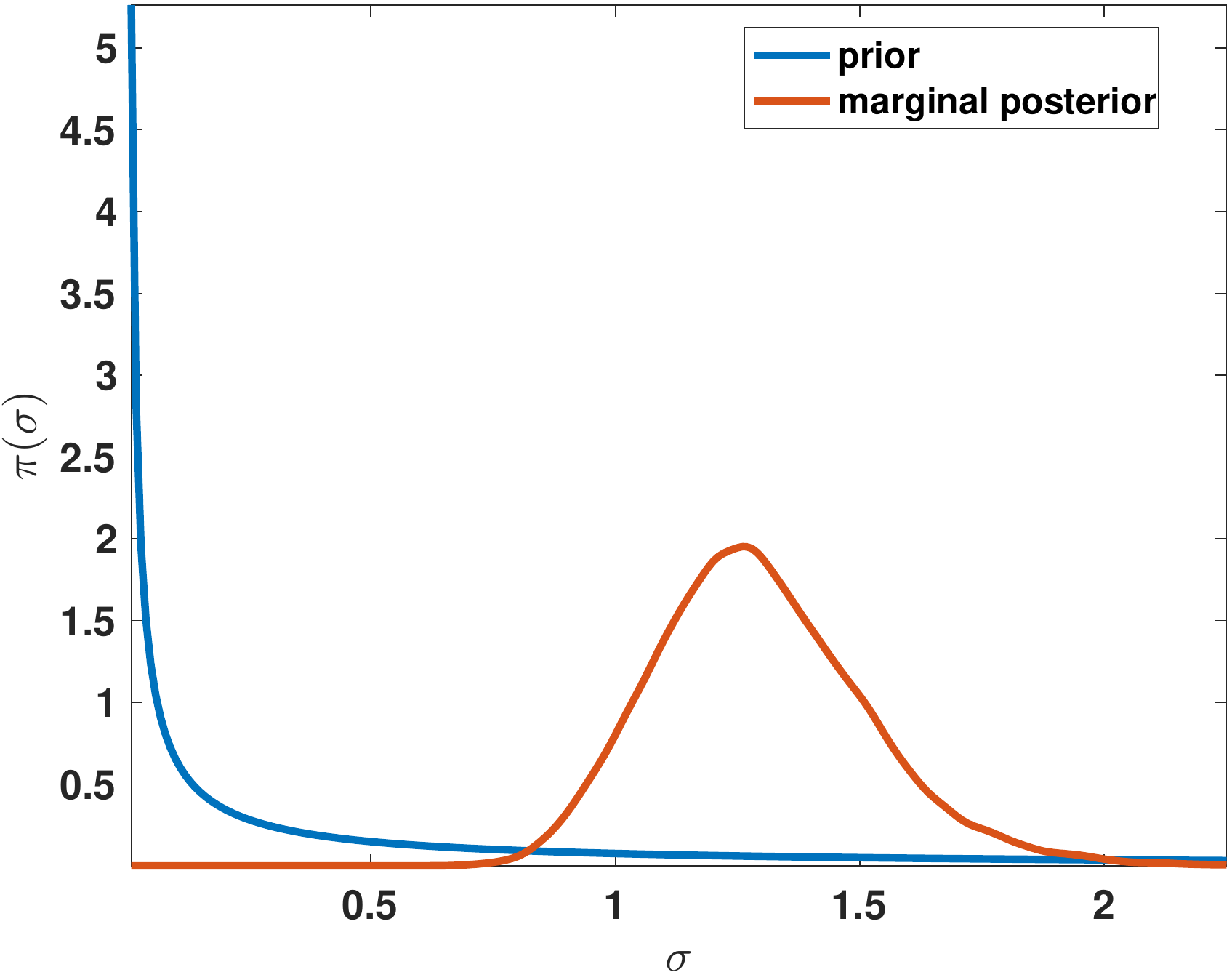}
    \caption{Priors and marginal posterior distributions for $\mu$ (left) and $\sigma$ (right) in the 2D deblurring example.}
    \label{fig:2Dmarginals}
\end{figure}

%\subsection{CT Example}
Trace plots of the samples from the thinned MCMC chains using weakly informative priors in the CT example (Section 4.2) are displayed in Supplementary Figure \ref{FIG:CTtrace}, and autocorrelation plots are displayed in Supplementary Figure \ref{FIG:CTautocorr}. In Supplementary Figure \ref{fig:CTcumavgs}, we plot the cumulative averages of the parameters $\kappa^2$ and $\upsilon$ of these samples. The PSRFs for these parameters were both 1.00. Even after thinning, we see some correlation between the samples.
\begin{figure}[tb]
\centering
\includegraphics[scale=0.25]{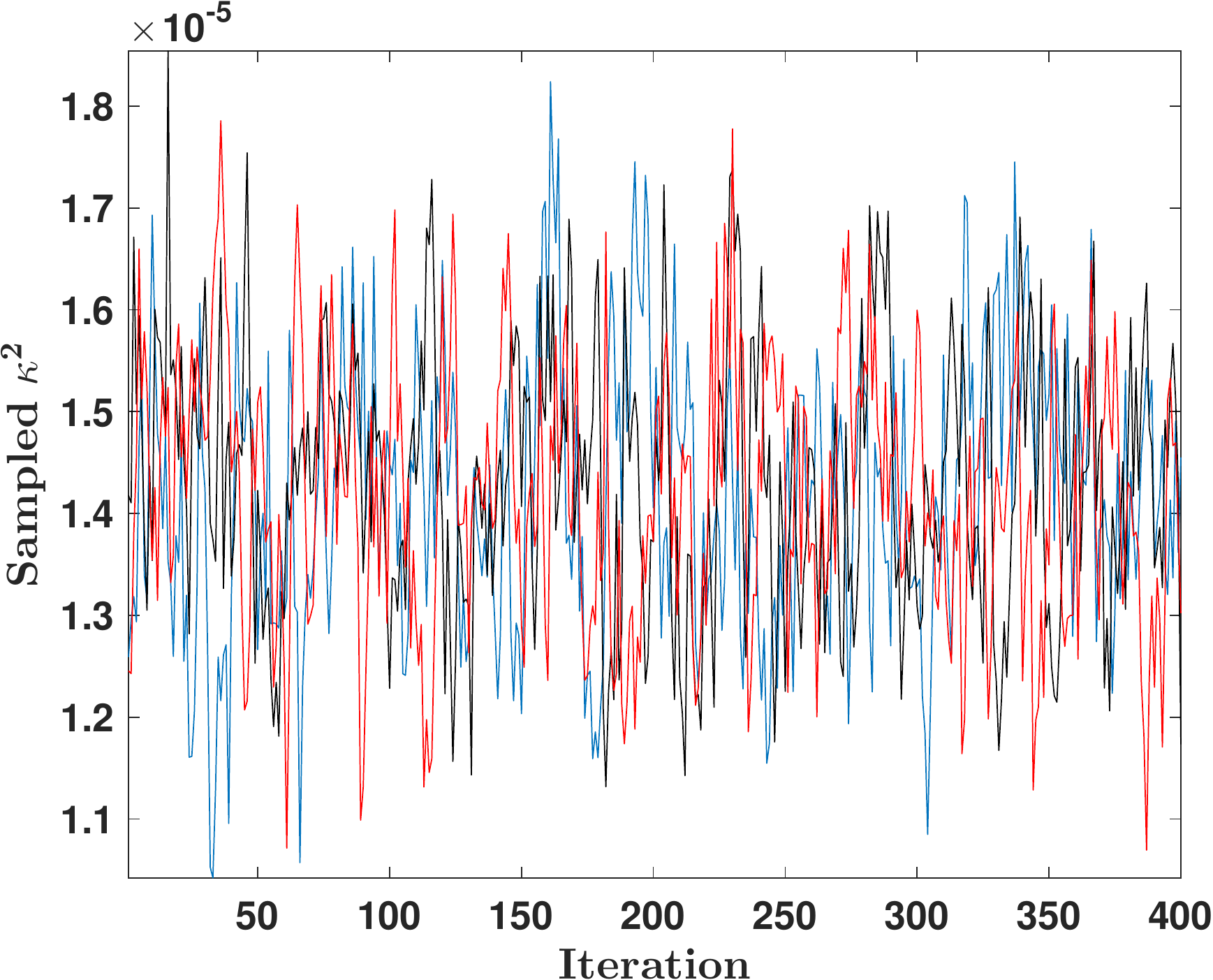}
\includegraphics[scale=0.25]{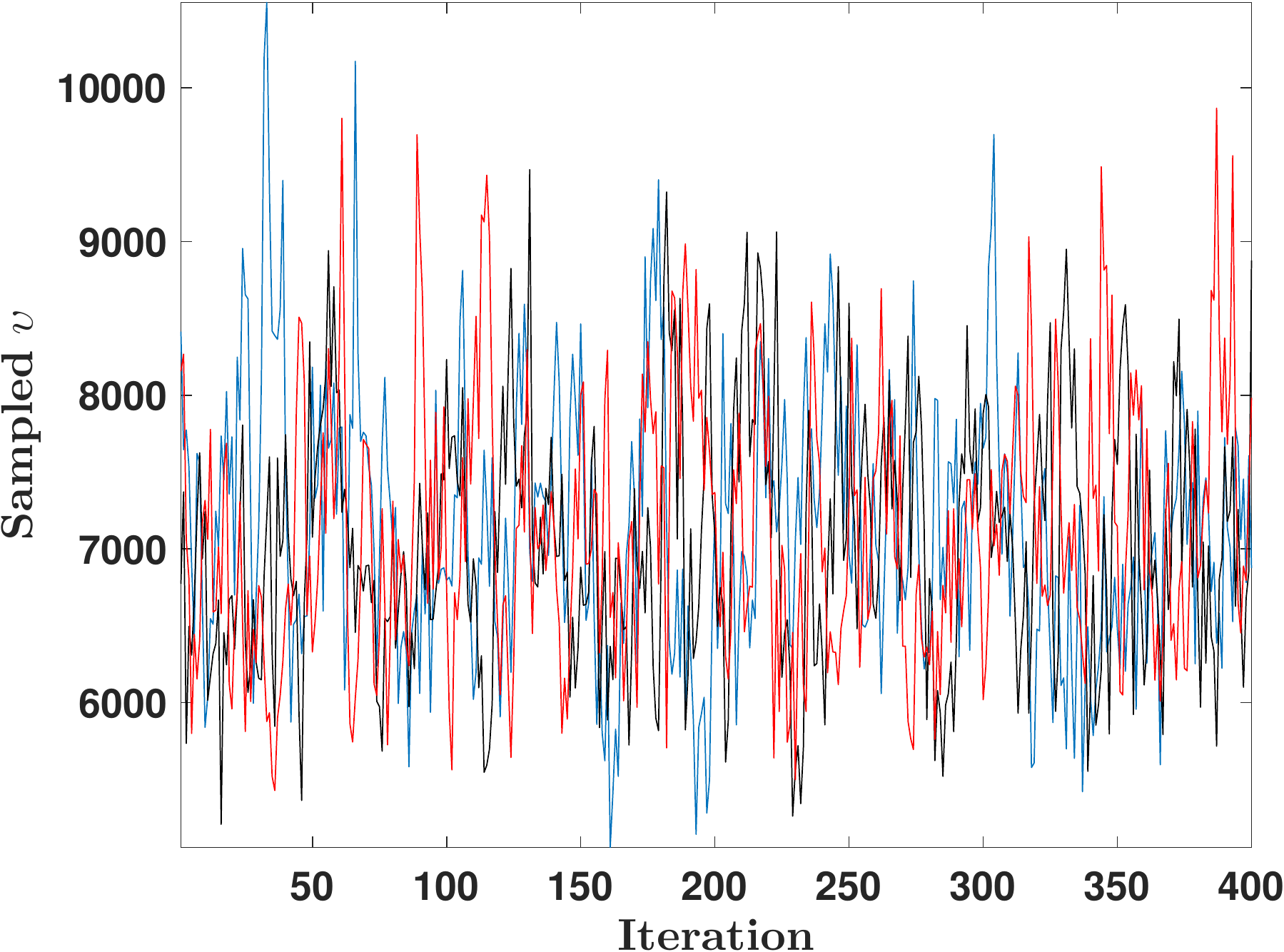}
\includegraphics[scale=0.25]{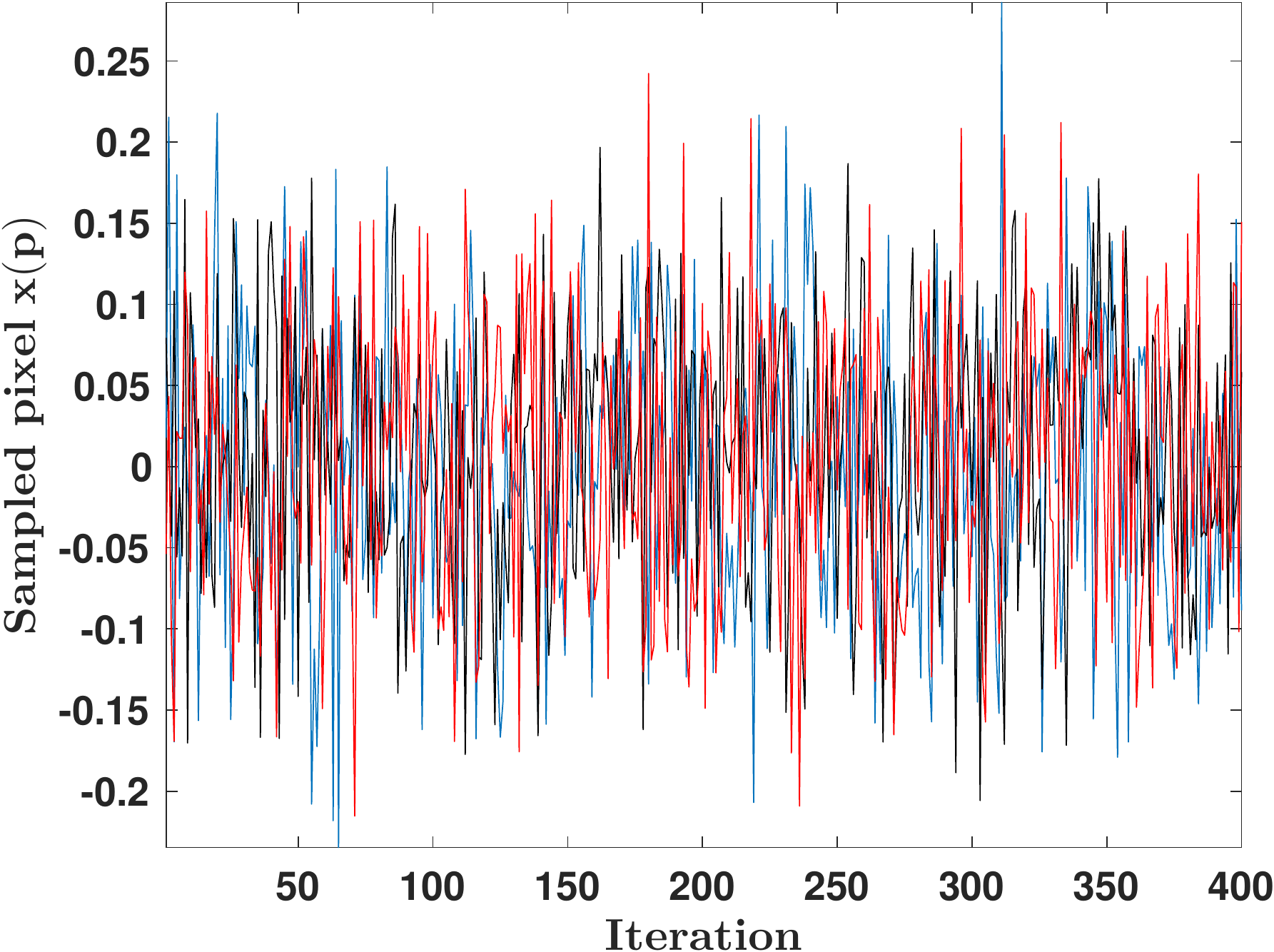}
\caption{Trace plots for the thinned chains in the CT image reconstruction example. Left: noise variance $\kappa^2$. Center: Variance ratio $\upsilon = \tau^2 / \kappa^2$. Right: A randomly chosen pixel of the image $\B{x}$.}
\label{FIG:CTtrace}
\end{figure}

\begin{figure}[tb]
\centering
\includegraphics[scale=0.3]{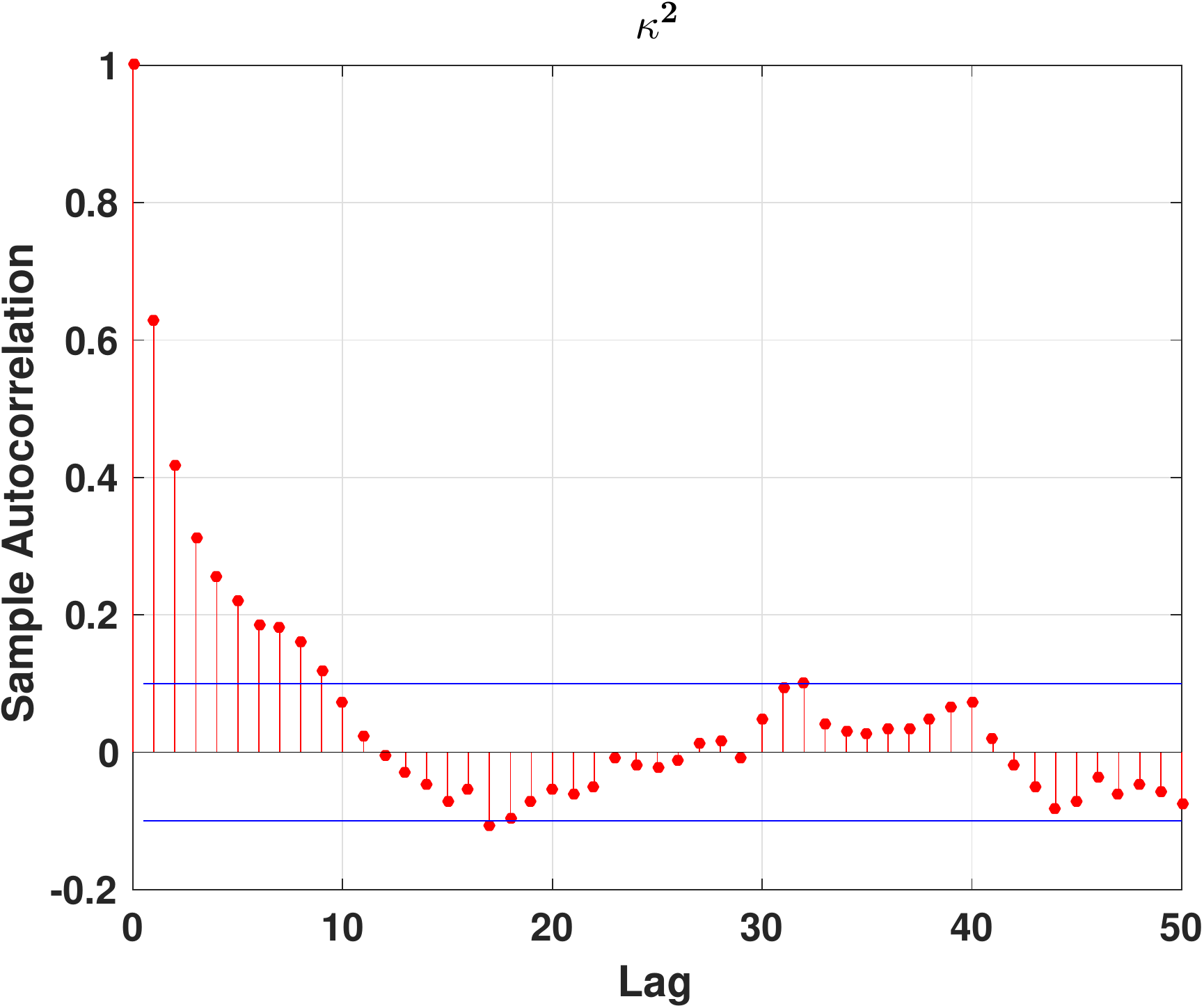}
\includegraphics[scale=0.3]{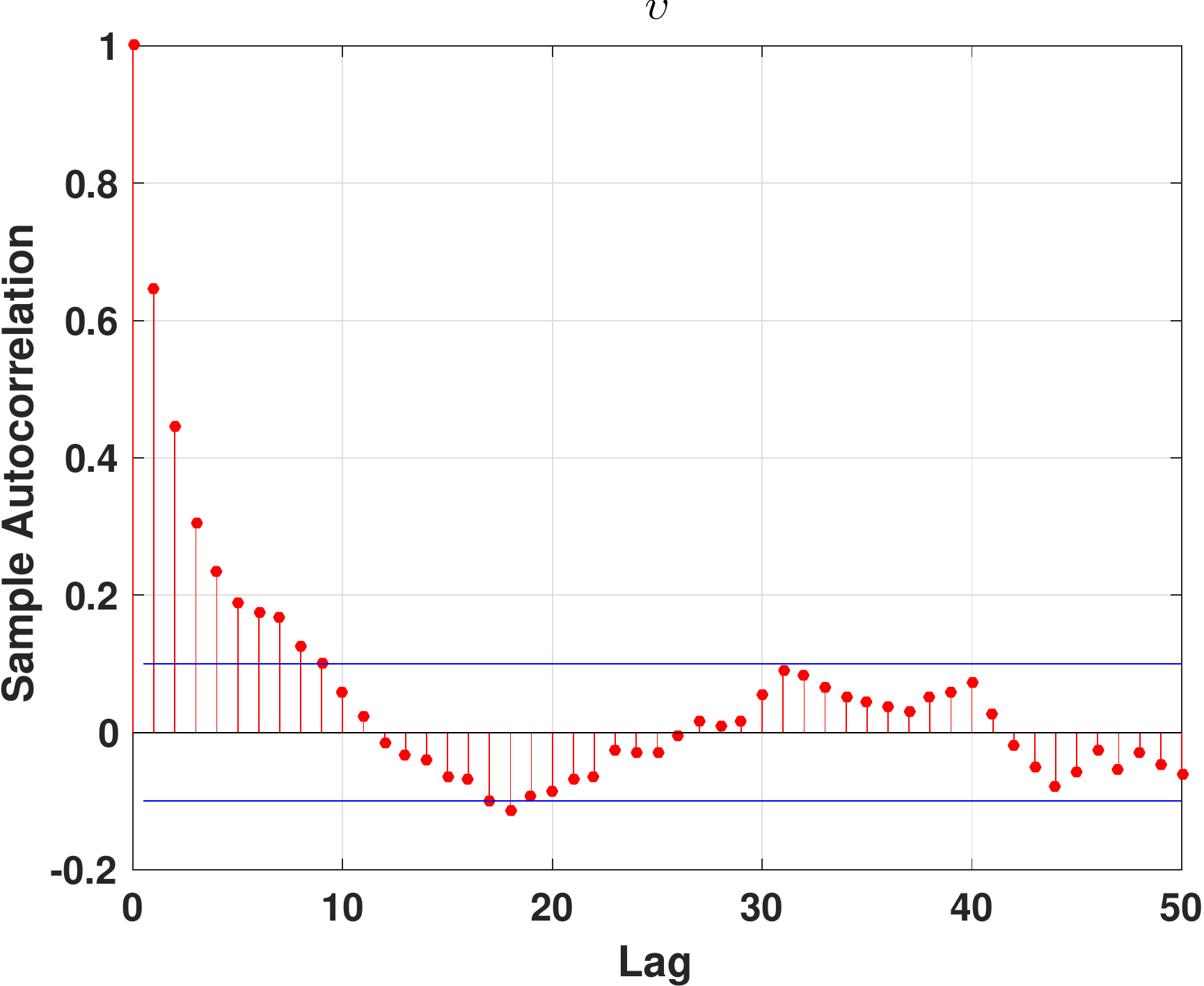}
\caption{Autocorrelation plots for the thinned chains in the CT image reconstruction example. Left: noise variance $\kappa^2$. Right: Variance ratio $\upsilon = \tau^2 / \kappa^2$.}
\label{FIG:CTautocorr}
\end{figure}

\begin{figure}[tb]
\centering
\includegraphics[scale=0.3]{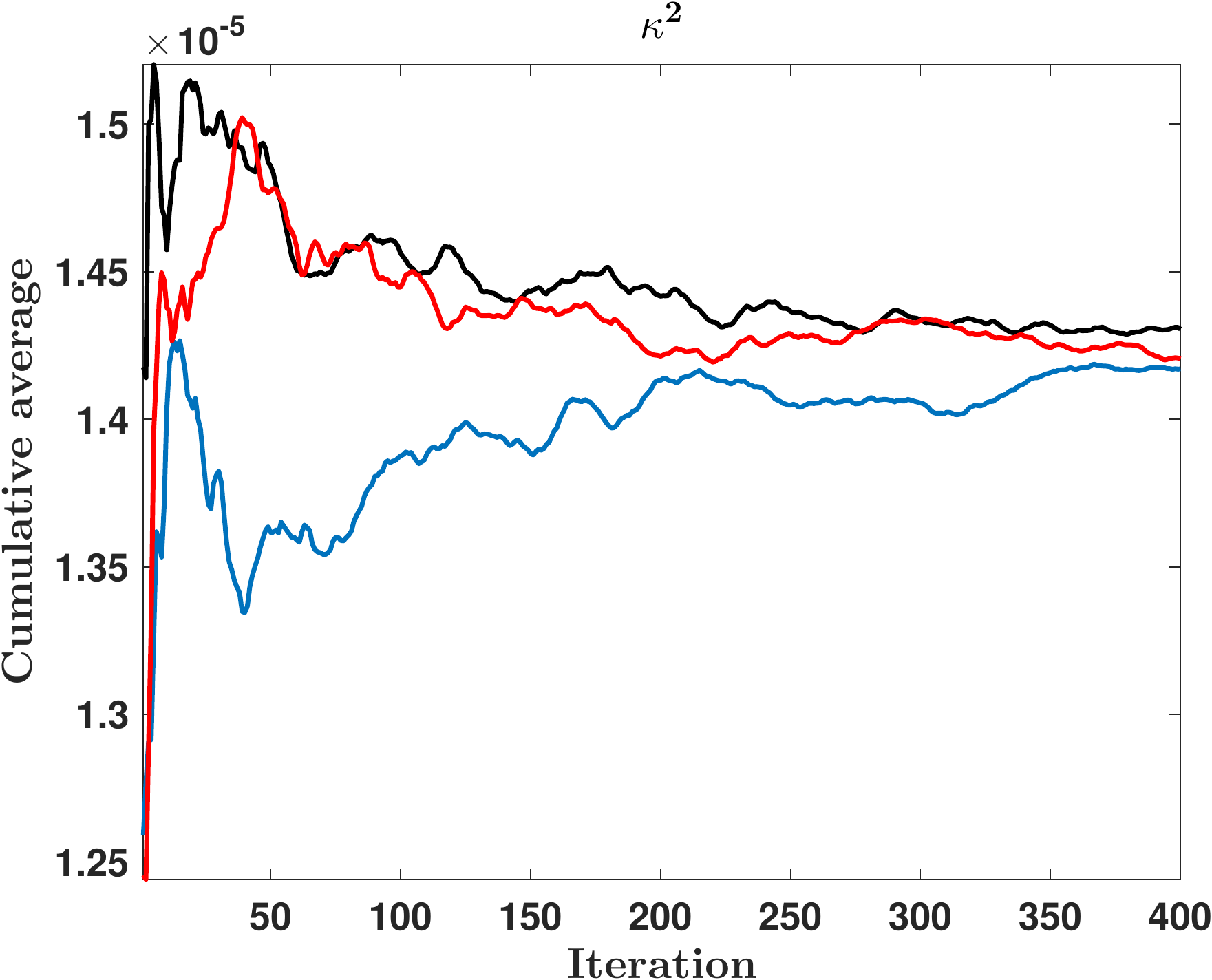}
\includegraphics[scale=0.3]{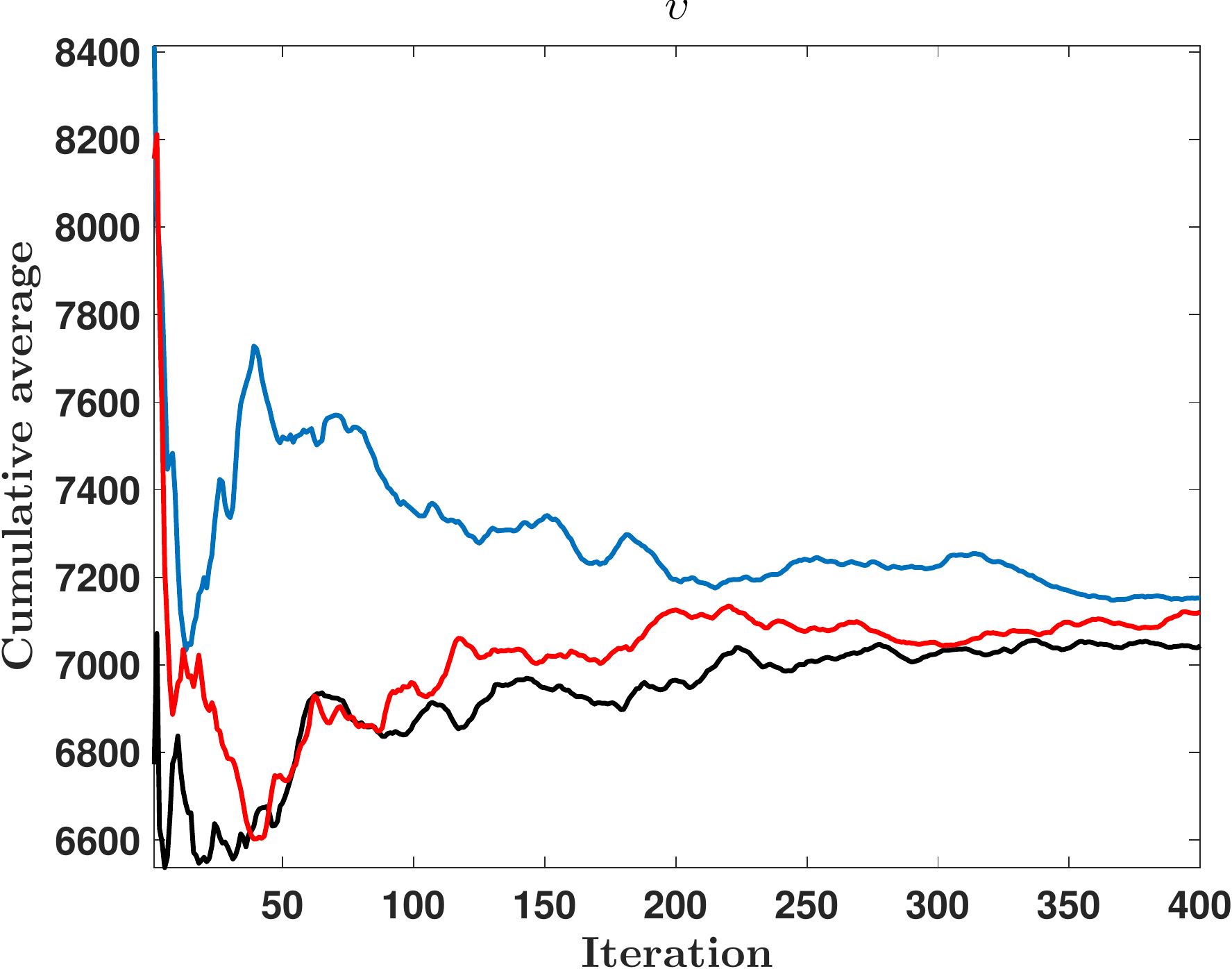}
    \caption{Cumulative averages of $\kappa^2$ and $\upsilon$ for each of the thinned MCMC chains in the CT example.}
    \label{fig:CTcumavgs}
\end{figure}

We next present convergence diagnostics for the CT example with the conjugate Gamma prior instead of proper Jeffreys to justify the comparison in the manuscript. We use vague Gamma priors for $\mu$ and $\sigma$ as in the 2D example, and the same prior for $\B{x}$ as in the original CT experiment. We again use Randomized SVD with target rank $\ell = 5000$ in the low-rank proposal distribution for $\B{x}$. We simulate three (randomly initialized) Markov chains using the MCMC algorithm with our proposed approach for 20,000 iterations. These chains are then thinned and the burn-in period discarded to produce an equal Monte Carlo sample size as in the previous experiment. Trace plots, autocorrelation plots, and PSRFs are used to determine approximate convergence of the chains. The trace plots are displayed in Supplementary Figure~\ref{fig:CTtracegamma} and the autocorrelation plots for $\mu$ and $\sigma$ are in Supplementary Figure~\ref{fig:CTautocorrgamma}. The PSRFs for $\mu$ and $\sigma$ were 1.00 and 1.01 respectively.% In Supplementary Figure~\ref{fig:CTcumavggamma}, we plot the cumulative averages for $\mu$ and $\sigma$.

\begin{figure}[tb]
	\includegraphics[scale=0.25]{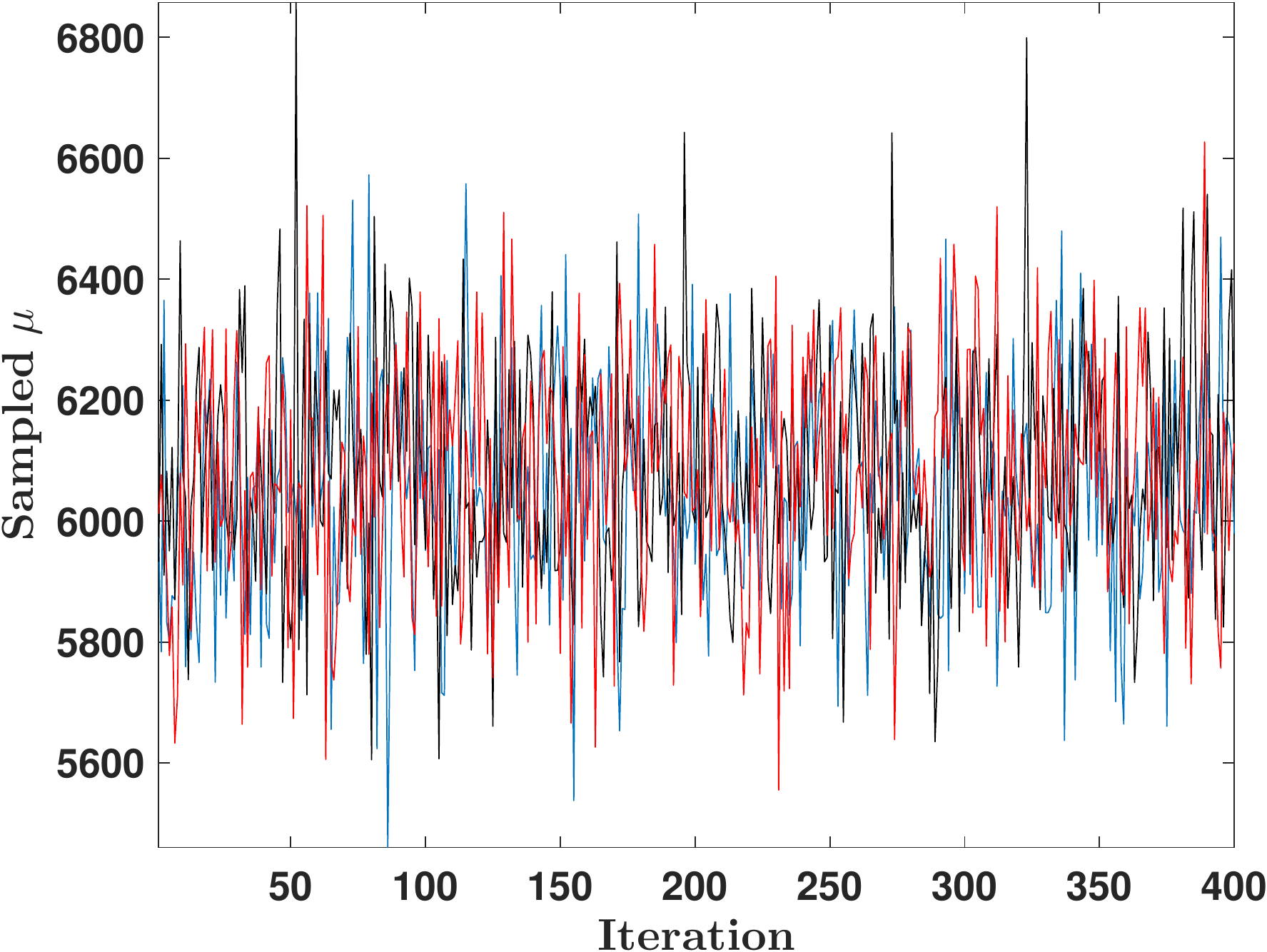}
	\includegraphics[scale=0.25]{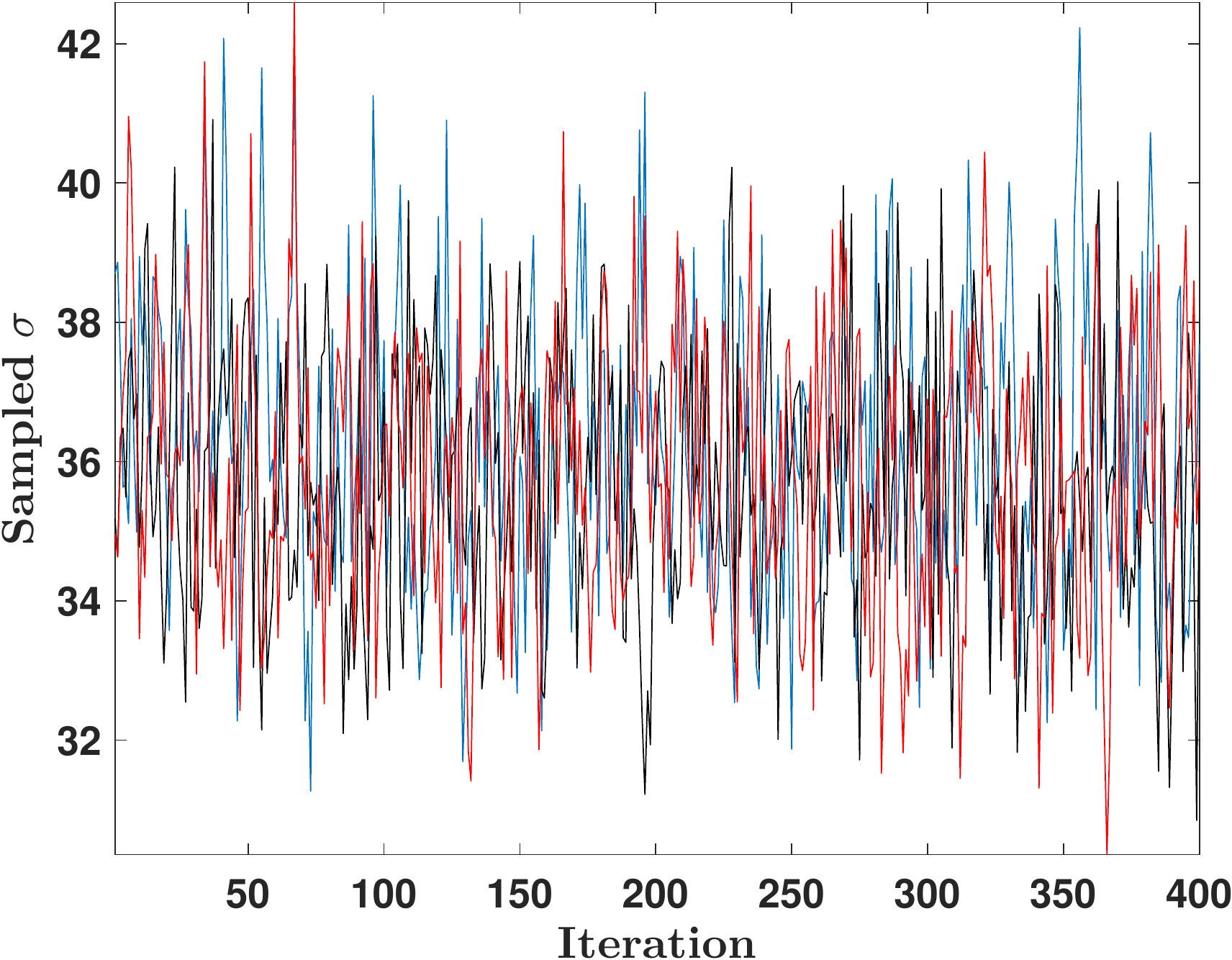}
	\includegraphics[scale=0.25]{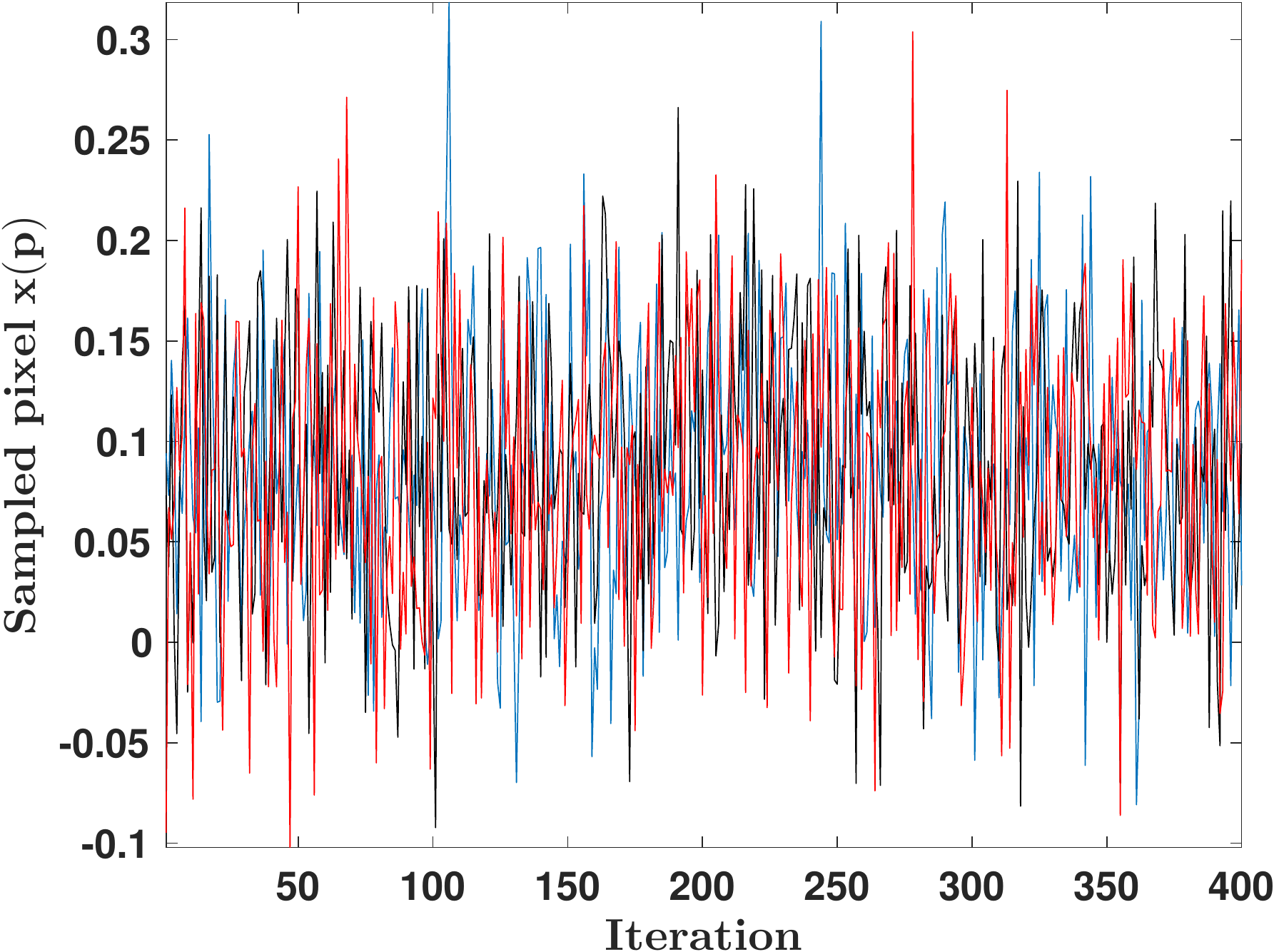}
	\caption{Trace plots for the three thinned chains using Gamma priors in the CT image reconstruction example. Left: noise precision $\mu$. Center: prior precision $\sigma$. Right: a randomly chosen pixel of the image $\B{x}$.}
	\label{fig:CTtracegamma}
\end{figure}

\begin{figure}[tb]
	\includegraphics[scale=0.3]{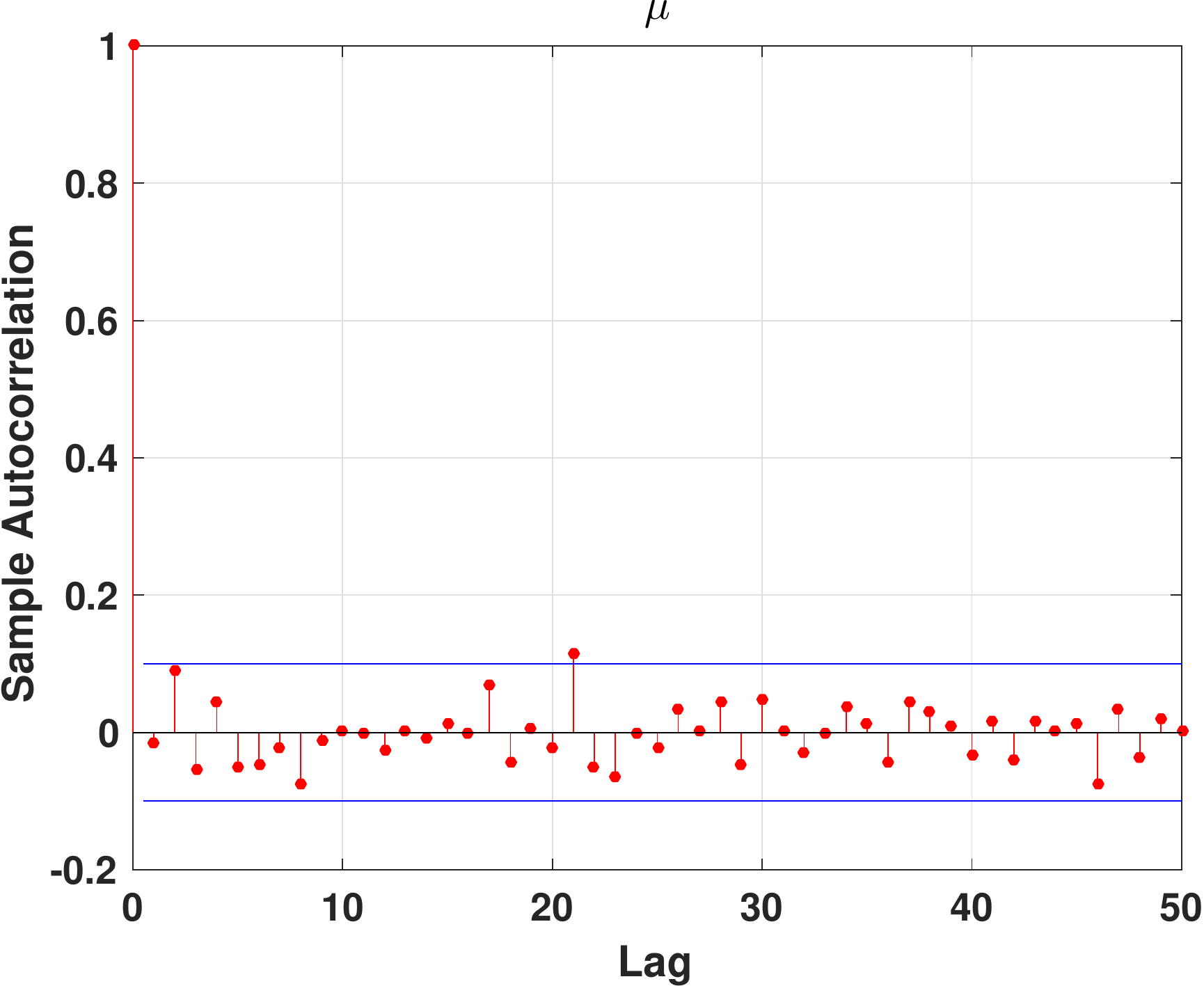}
	\includegraphics[scale=0.3]{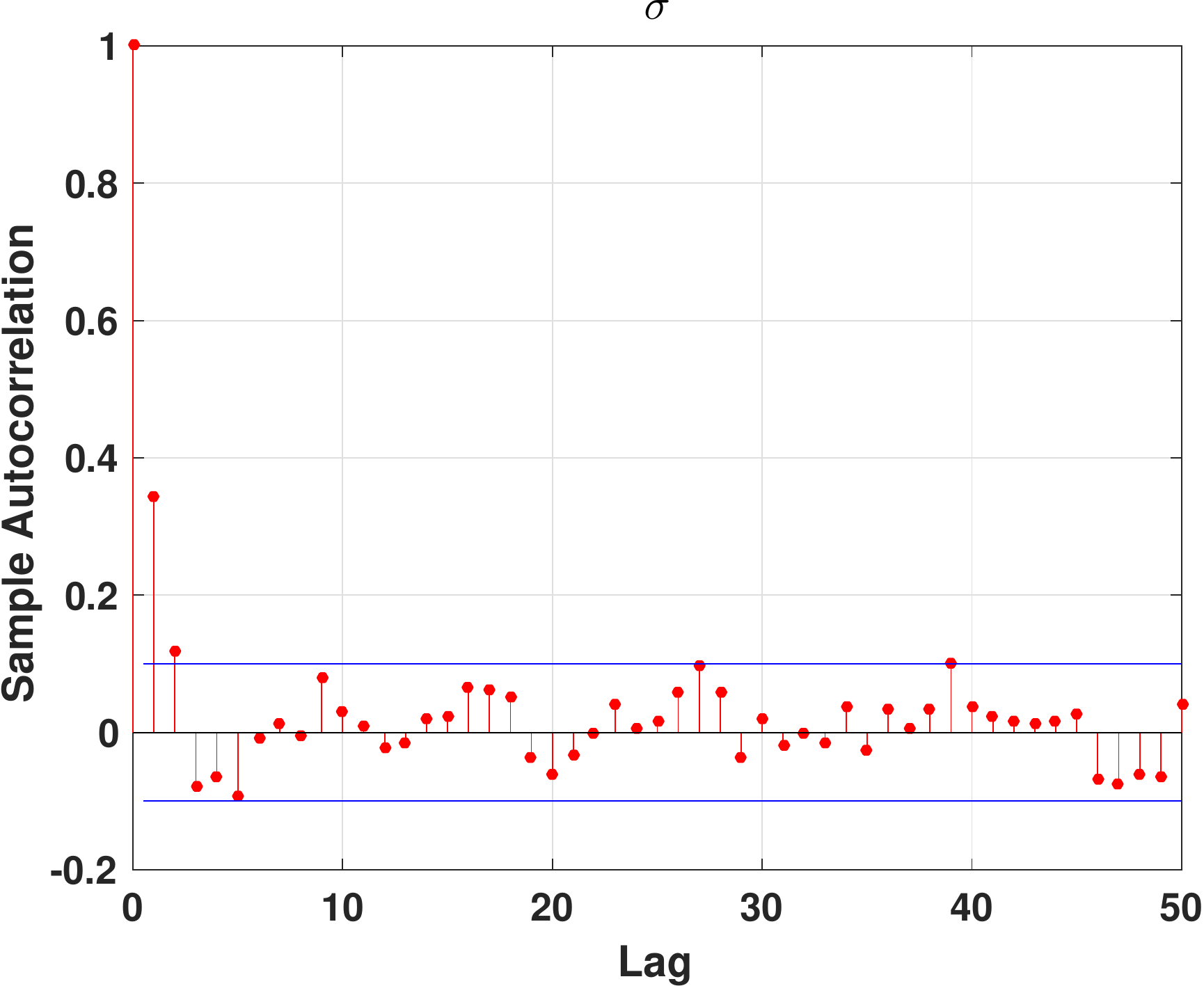}
	\caption{Autocorrelation plots for one thinned chain using Gamma priors in the CT image reconstruction example. Left: noise precision $\mu$. Right: prior precision $\sigma$.}
	\label{fig:CTautocorrgamma}
\end{figure}

%\begin{figure}[!htb]\label{fig:CTcumavggamma}
%	\includegraphics[scale=0.3]{Figures/CT-cumavg-mu-Gamma-THIN}
%	\includegraphics[scale=0.3]{Figures/CT-cumavg-sigma-Gamma-THIN}
%	\caption{Cumulative averages for the three thinned chains using Gamma priors in the CT image reconstruction example. Left: noise precision $\mu$. Right: prior precision $\sigma$.}
%\end{figure}

%The total computation time for the MCMC simulations was 99,459 seconds, or about 27.6 hours. Thus collecting 40,000 samples as in the proper Jeffreys case would require a comparable computation time.

%
%\clearpage
%
%\bibliography{mcmc_paper-new} % Using Arvind's updated bib file
%\bibliographystyle{plain}
%\end{document}

\end{document}